\documentclass[12pt]{article}
\usepackage{amsfonts}
\usepackage{amsmath}
\usepackage{amssymb}

\usepackage{graphicx}
\usepackage{amsmath}
\usepackage{multirow}
\usepackage{subcaption}
\usepackage{float}

\begin{document}

\title {Bifurcations of a predator-prey system with cooperative hunting and Holling III functional response
\author{\small Yong Yao$^a$,
\ \ Teng Song$^b$,
\ \ Zuxiong Li$^c$\thanks{Corresponding author. lizx0427@126.com (Z. Li).}\\
{\small\it a. School of Mathematics and Physics, Wuhan Institute of Technology,}\\
{\small\it  Wuhan, Hubei 430205, P.R. China}\\
{\small\it b. School of Science, Wuhan University of Technology,}\\
{\small\it  Wuhan, Hubei 430070, P.R. China}\\
{\small\it c. College of Mathematics and  Statistics, Chongqing Three Gorges University,}\\
{\small\it Wanzhou, Chongqing 404100, P.R. China}\\
}
\date{}
}


\newtheorem{thm}{Theorem}
\newtheorem{defi}[thm]{Definition}
\newtheorem{lem}[thm]{Lemma}
\newtheorem{prop}[thm]{Proposition}
\newtheorem{coro}[thm]{Corollary}
\newtheorem{exam}[thm]{Example}
\newtheorem{rem}{Remark}

\def\theequation{\thesection.\arabic{equation}}
\def\qed\nopagebreak\hfill{\rule{3pt}{6pt}}
\makeatletter

\@addtoreset{equation}{section}

\makeatother

\def\qed{\hfill \rule{3pt}{6pt}}


\maketitle

\begin{center}
\begin{minipage}{15cm} \hspace{0.1cm}
\noindent{\bf Abstract:}
In this paper, we consider the dynamics of a predator-prey system of Gause type with cooperative hunting among predators and Holling III functional response.
The known work numerically shows that the system exhibits saddle-node and Hopf bifurcations except homoclinic bifurcation for some special parameter values.
Our results show that there are a weak focus of multiplicity three and a degenerate equilibrium with double zero eigenvalues (i.e., a cusp of codimension two) for general parameter conditions and the system can exhibit various bifurcations as perturbing the bifurcation parameters appropriately, such as the transcritical and the pitchfork bifurcations at the degenerate boundary equilibrium, the saddle-node and the Bogdanov-Takens bifurcations at the degenerate positive equilibrium  and the Hopf bifurcation around the weak focus.
The comparative study demonstrates that the dynamics are far richer and more complex than that of the system without cooperative hunting among predators.
The analysis results reveal that appropriate intensity of cooperative hunting among predators is beneficial for the persistence of predators and the diversity of ecosystem.
\end{minipage}
\end{center}

\parskip=0pt 


\section{Introduction}
Mathematical models have played a central role in population biology.
To understand the interactions among populations, many mathematical modeling
efforts have been made through investigating the food chain dynamics.
One of the most classical mathematical models is the Gause type predator-prey system taking the following form(\cite{Kot}),
\begin{eqnarray}
\left\{
\begin{array}{l}
\frac{dx}{dt}=rx(1-\frac{x}{K})-y\Phi(\cdot),\\
\frac{dy}{dt}=cy\Phi(\cdot)-dy,
\end{array}
\right.
\label{(i1)}
\end{eqnarray}
where $x$ and $y$ are the prey density and predator density respectively,
$r$ represents the per capita intrinsic growth rate of prey, $K$  stands for the carrying capacity of prey, $c$ is the conversion rate,
$d$ represents  the mortality rate of predators and
$\Phi(\cdot)$ is the functional response referring to the consumption rate of a predator to prey. The functional response usually depends on many factors, such as the population density, the encounter rate and the handling time. The following functional responses have been extensively used in modeling population dynamics.

(i) Lotka-Volterra functional response
\begin{eqnarray*}
\left.
\begin{array}{l}
\Phi(x)=ex,
\end{array}
\right.
\end{eqnarray*}
where $e$ is the encounter rate, which is an unbounded and linear function(\cite{Lotka,Volterra}).

(ii) Holling II functional response
\begin{eqnarray*}
\left.
\begin{array}{l}
\Phi(x)=\frac{ex}{1+ehx},
\end{array}
\right.
\end{eqnarray*}
where $e$ is the encounter rate and $h$ is the handling time, which is a bounded and monotonically increasing nonlinear function(\cite{Holling}).

(iii) Holling III functional response
\begin{eqnarray*}
\left.
\begin{array}{l}
\Phi(x)=\frac{ex^2}{1+ehx^2},
\end{array}
\right.
\end{eqnarray*}
where $e$ and $h$ have the same biological meanings,  which also is a bounded and monotonically increasing nonlinear
function like Holling II functional response.  Although both Holling II and III functional responses
are approaching an asymptote, the former is decelerating and the latter is sigmoid(\cite{Holling}).

(iv) Holling IV functional response
\begin{eqnarray*}
\left.
\begin{array}{l}
\Phi(x)=\frac{ex}{1+ehx+eh_0x^2},
\end{array}
\right.
\end{eqnarray*}
where $e$ and $h$ still have the same biological meanings and $h_0$ describes how handling time increases with prey density due to group defense, which is a bounded and nonmonotonic function(\cite{Collings}).

System (\ref{(i1)}) with various functional responses has attracted a lot of attention from mathematicians and theoretical biologists and displays interesting dynamics despite their complete dependence on prey density. System (\ref{(i1)}) with Lotka-Volterra functional response
has relatively simple dynamics except the rise of the unique positive equilibrium from a transcritical bifurcation  since the linear functional response is sufficiently simple. System (\ref{(i1)}) with nonlinear functional responses has relatively complicated dynamical properties. For example, the system with Holling II functional response called Rosenzweig-MacArthur model(\cite{Rosenzweig}) can exhibit a unique positive equilibrium and a unique stable limit cycle  around the unstable positive equilibrium, which are induced by the transcritical bifurcation and the Hopf bifurcation respectively.
Therefore, both of predators and prey coexist at either the coexistence equilibrium or the limit cycle.
The appearance of the stable oscillation demonstrates exactly the paradox of enrichment(\cite{Cheng, Hsu, Huang}).
Chen and Zhang(\cite{ChenZhang1986}) considered the system with Holling III functional response and obtained parameter conditions for the global stability of the unique positive equilibrium and the uniqueness of the limit cycle around the positive equilibrium.
Moreover, the bifurcation diagram and the corresponding phase portraits of the system with Holling III functional response also were exhibited in the book of Bazykin(\cite{Bazykin}).
Some researchers investigated the system with Holling IV functional response (the parameters chosen such that the denominator of the functional response does not vanish for non-negative prey density) and obtained that the system can display complicated bifurcation phenomena, such as the degenerate Hopf bifurcation and the Bogdanov-Takens bifurcation of codimensions two and three(\cite{HuangXiao, XiaoZhu, ZhuCampbellWolkowicz}).

The independence of functional response of predator density is not always true in the ecosystem  since the independence
means that one predator affects the growth rate of its prey independently of its conspecifics. The two inverse ecological phenomena, interference and
facilitation (or cooperative hunting) among predators, therefore need the predator-dependent functional response to reflect. The predator interference phenomenon first received the attention of scholars and was modeled by the ratio-dependent(\cite{Aguirre, ArditiGinzburg}) and the Beddington-DeAngelis(\cite{Beddington, DeAngelis, LiuBeretta})  functional responses, etc.
Similarly, the cooperative hunting among predators also is a common phenomenon in nature(\cite{Dugatkin}), which even is an effective mechanism to promote the evolution and diversity of species(\cite{Packer}).
A lot of recent studies paid attention to the cooperative hunting among predators(\cite{Alves,Berec,ChowJangWang19,ChowJangWang20,Jang,Pal,Sen,Yan}).
Alves and Hilker(\cite{Alves}) considered the following four types of predator-dependent functional response to model the cooperative hunting among predators based on the above-mentioned four types of functional response.

(a) Lotka-Volterra functional response with cooperative hunting
\begin{eqnarray*}
\left.
\begin{array}{l}
\Phi(x,y)=(e+ay)x.
\end{array}
\right.
\end{eqnarray*}

(b) Holling II functional response with cooperative hunting
\begin{eqnarray*}
\left.
\begin{array}{l}
\Phi(x,y)=\frac{(e+ay)x}{1+h(e+ay)x}.
\end{array}
\right.
\end{eqnarray*}

(c) Holling III functional response with cooperative hunting
\begin{eqnarray*}
\left.
\begin{array}{l}
\Phi(x,y)=\frac{(e+ay)x^2}{1+h(e+ay)x^2}.
\end{array}
\right.
\end{eqnarray*}

(d) Holling IV functional response with cooperative hunting
\begin{eqnarray*}
\left.
\begin{array}{l}
\Phi(x,y)=\frac{(e+ay)x}{1+h(e+ay)x+h_0(e+ay)x^2}.
\end{array}
\right.
\end{eqnarray*}
The parameter $a$ is the intensity of cooperative hunting among predators.  The modeling motivation is that the encounter rate of predators with their prey
is no longer a constant due to the cooperative hunting among predators, but rather an increasing function of predator density, which is analogous to the
encounter-driven functional response put forward by Berec(\cite{Berec}).
Alves and Hilker(\cite{Alves}) investigated the two-parameter bifurcation of system (\ref{(i1)}) with the four types of functional response (a)-(d) respectively by numerical simulations and came to the conclusion that cooperative hunting among predators is a mechanism to produce the Allee effect in predators. Their numerical results showed that system (\ref{(i1)}) with Lotka-Volterra, Holling II and Holling IV types functional response (i.e., (a), (b) and (d)) can display Bogdanov-Takens bifurcation of codimension two for some special parameters, but the system with Holling III type functional response (i.e., (c)) can not exhibit  Bogdanov-Takens bifurcation of codimension two, exactly,
there is neither a Bogdanov-Takens bifurcation point nor a homoclinic bifurcation curve in the two-parameter bifurcation diagram.
Whereafter, Zhang and Zhang(\cite{ZhangZhang}) further considered system (\ref{(i1)}) with Lotka-Volterra type functional response, and
discussed the location of equilibria qualitatively and gave analytical conditions for not only the Bogdanov-Takens bifurcation of codimension two at positive equilibrium but also the transcritical and pitchfork bifurcations at boundary equilibrium.
For system (\ref{(i1)}) with Holling II type functional response,
which actually is the  special case of Berec's system(\cite{Berec}), Yao(\cite{Yao2020}) further
gave the parameter conditions for the qualitative properties of equilibria and analysed the bifurcations at nonhyperbolic equilibria including saddle-node, transcritical, pitchfork and Hopf bifurcations and Bogdanov-Takens bifurcation of codimension two
by applying the pseudo-division reduction and complete discrimination system of parametric polynomial.
The comparative studies revealed that system (\ref{(i1)}) with  cooperative hunting has even richer and more complicated dynamics  than  the system without  cooperative hunting for both Lotka-Volterra type and Holling II type functional responses.

In this paper, we consider system (\ref{(i1)}) with  Holling III type functional response and analyse whether the system  can also exhibit complicated dynamics. The system takes the following form
\begin{eqnarray}
\left\{
\begin{array}{l}
\frac{dx}{dt}=rx(1-\frac{x}{K})-\frac{(e+ay)x^2y}{1+h(e+ay)x^2},\\
\frac{dy}{dt}=\frac{c(e+ay)x^2y}{1+h(e+ay)x^2}-dy.
\end{array}
\right.
\label{(1)}
\end{eqnarray}
By the scaling
$\tilde{x}:=\sqrt{\frac{ce}{d}}x$, $\tilde{y}:=\sqrt{\frac{e}{dc}}y$, $\alpha:=\frac{a}{e}\sqrt{\frac{dc}{e}}$, $\sigma:=\frac{r}{d}$, $\kappa:=\sqrt{\frac{ce}{d}}K$, $\tilde{h}:=\frac{hd}{c}$ and $\tau:=dt$,
system (\ref{(1)}) can be  simplified as
\begin{eqnarray}
\left\{
\begin{array}{l}
\frac{dx}{dt}=\sigma x(1-\frac{x}{\kappa})-\frac{(1+\alpha y)x^2y}{1+h(1+\alpha y)x^2},\\
\frac{dy}{dt}=\frac{(1+\alpha y)x^2y}{1+h(1+\alpha y)x^2}-y,
\end{array}
\right.
\label{(2)}
\end{eqnarray}
where we still denote $\tilde{x}$, $\tilde{y}$, $\tilde{h}$ and $\tau$ as $x$, $y$, $h$ and $t$ respectively,
and $\alpha$ represents the intensity of cooperative hunting among predators.
Alves and Hilker(\cite{Alves}) displayed the two-parameter bifurcation diagram when varying the cooperative hunting $\alpha$ and intrinsic growth rate of prey $\sigma$ by choosing $h=0.9$ and $\kappa=0.8$, in which the
Hopf bifurcation curve does not touch the saddle-node bifurcation curve and the homoclinic bifurcation curve does not appear,
which means that system (\ref{(2)}) does not undergo a Bogdanov-Takens bifurcation for the special parameter values.
Nevertheless, we can prove that the system can exhibit not only a Bogdanov-Takens bifurcation of codimension two at the degenerate positive equilibrium but also a degenerate Hopf bifurcation at the nonhyperbolic positive equilibrium.
The time-rescaling transformation shows that system (\ref{(2)}) and the following quartic polynomial system are orbitally equivalent
\begin{eqnarray}
\left\{
\begin{array}{l}
\frac{dx}{dt}=x\{\sigma(\kappa-x)(1+h(1+\alpha y)x^2)-k(1+\alpha y)xy\},\\
\frac{dy}{dt}=\kappa y\{(1-h)(1+\alpha y)x^2-1\}.
\end{array}
\right.
\label{(3)}
\end{eqnarray}
The paper is organized as follows. In section 2,
we prove qualitatively that system (\ref{(3)}) has at most four equilibria and give
the parameter conditions for the existence and their qualitative properties.
In section 3, we distinguish the topological types of the nonhyperbolic equilibria and analyse various possible bifurcations around them, such as
a transcritical bifurcation and a pitchfork bifurcation at the degenerate boundary equilibrium,
a saddle-node bifurcation and a Bogdanov-Takens bifurcation of codimension two at the degenerate positive equilibrium and a Hopf bifurcation at the weak focus of multiplicity at most three.
Finally, we make numerical simulations to demonstrate our theoretical results and
end the paper with a brief discussion in section 4.


\section{Equilibria and Their Properties}
In this section we discuss the existence of equilibria of system (\ref{(3)}) and their qualitative properties.
We first give the following partition of parameters $(\alpha,\kappa,\sigma)$ as $0<h<1$ such that we can state our results conveniently
\begin{eqnarray*}
\left.
\begin{array}{l}
\mathbb{R}_+^3:=\{(\alpha,\kappa,\sigma)\in\mathbb{R}^3:\alpha>0,\kappa>0, \sigma>0\}=\mathcal{P}_1\cup\mathcal{P}_2\cup\mathcal{P}_3\cup\mathcal{P}_4\cup\mathcal{P}_5\cup\mathcal{P}_6\cup\mathcal{P}_7,
\end{array}
\right.
\end{eqnarray*}
where
\begin{eqnarray*}
\left.
\begin{array}{l}
\mathcal{P}_1:=\{(\alpha,\kappa,\sigma)\in\mathbb{R}_+^3: \alpha_1<\alpha<\alpha_2, \kappa<\kappa_1\},\\
\mathcal{P}_2:=\{(\alpha,\kappa,\sigma)\in\mathbb{R}_+^3: \alpha=\alpha_1, \kappa\leq \kappa_1\},\\
\mathcal{P}_3:=\{(\alpha,\kappa,\sigma)\in\mathbb{R}_+^3: \alpha=\alpha_2, \kappa<\kappa_1\},\\
\mathcal{P}_4:=\{(\alpha,\kappa,\sigma)\in\mathbb{R}_+^3: \alpha<\alpha_1, \kappa\leq \kappa_1\},\\
\mathcal{P}_5:=\{(\alpha,\kappa,\sigma)\in\mathbb{R}_+^3: \kappa>\kappa_1\},\\
\mathcal{P}_6:=\{(\alpha,\kappa,\sigma)\in\mathbb{R}_+^3: \alpha>\frac{2}{\kappa\sigma}, \kappa=\kappa_1\},\\
\mathcal{P}_7:=\{(\alpha,\kappa,\sigma)\in\mathbb{R}_+^3: \alpha>\alpha_2, \kappa<\kappa_1\}
\end{array}
\right.
\end{eqnarray*}
with $\kappa_1:=\frac{1}{\sqrt{1-h}}$ and
\begin{eqnarray}
\left.
\begin{array}{l}
\alpha_1:=-\frac{16}{81}\frac{9\kappa^2(h-1)+16}{\sigma\kappa^3(h-1)}-\frac{1+\sqrt{3}i}{2}\sqrt[3]{-\frac{A}{2}+\sqrt{\Delta}}+\frac{-1+\sqrt{3}i}{2}\sqrt[3]{-\frac{A}{2}-\sqrt{\Delta}},\\
\alpha_2:=-\frac{16}{81}\frac{9\kappa^2(h-1)+16}{\sigma\kappa^3(h-1)}+\sqrt[3]{-\frac{A}{2}+\sqrt{\Delta}}+\sqrt[3]{-\frac{A}{2}-\sqrt{\Delta}},\\
A:=\frac{16}{531441 \sigma^3 \kappa^9 (h-1)^3}\{10935 \kappa^8 (h-1)^4+77760 (h-1)^3 \kappa^6\\
\phantom{A:=}
+1492992 (h-1)^2 \kappa^4+(3538944 h-3538944) \kappa^2+2097152\},\\
\Delta:=\frac{64}{14348907}\frac{(h \kappa^2-\kappa^2+1) (3 h \kappa^2-3 \kappa^2+128)^3}{\kappa^8 \sigma^6 (h-1)}.
\end{array}
\right.
\label{(2-1)}
\end{eqnarray}

The following theorem shows the various parametric conditions for the existence of equilibria of system (\ref{(3)}) and their corresponding qualitative properties.
\begin{thm}
System (\ref{(3)}) has at most four equilibria as follows.
(i) System (\ref{(3)}) always has two boundary equilibria $E_{0}(0,0)$ and $E_{\kappa}(\kappa,0)$. Moreover, $E_{0}(0,0)$ is a saddle and $E_{\kappa}(\kappa,0)$ is a stable node (or a saddle or degenerate) if $\kappa>0$ and $h\geq1$ or $0<\kappa<\kappa_1$ and $0<h<1$ (or $\kappa>\kappa_1$ and $0<h<1$ or $\kappa=\kappa_1$ and $0<h<1$).
(ii) System (\ref{(3)}) has at most two positive equilibria. More concretely,
(ii.a) system (\ref{(3)}) has no positive equilibrium if $h\geq1$ and $(\alpha,\kappa,\sigma)\in\mathbb{R}_+^3$ or $0<h<1$ and $(\alpha,\kappa,\sigma)\in\mathcal{P}_1\cup\mathcal{P}_2\cup\mathcal{P}_4$;
(ii.b) system (\ref{(3)}) has a unique positive equilibrium $E_1(x_1,y_1)$ if $0<h<1$ and $(\alpha,\kappa,\sigma)\in\mathcal{P}_5\cup\mathcal{P}_6$ and system (\ref{(3)}) has a unique positive equilibrium $E_*(x_*,y_*)$, if $0<h<1$ and $(\alpha,\kappa,\sigma)\in\mathcal{P}_3$;
(ii.c) system (\ref{(3)}) has two positive equilibria $E_1(x_1, y_1)$ and $E_2(x_2, y_2)$ with $x_1<x_2$, if $0<h<1$ and $(\alpha,\kappa,\sigma)\in\mathcal{P}_7$.
Furthermore, $E_1(x_1,y_1)$ is a stable node or focus (or an unstable node or focus or a center or weak focus) if $T(x_1)<0$ (or $>0$ or $=0$), $E_2(x_2, y_2)$ is a saddle and $E_*(x_*,y_*)$ is degenerate,
where
\begin{eqnarray*}
\left.
\begin{array}{l}
y_i=\sigma x_i(1-\frac{x_i}{\kappa})$~~~$(i=1, 2, *)
\end{array}
\right.
\end{eqnarray*}
and
\begin{eqnarray*}
\left.
\begin{array}{l}
T(x_1):=-\kappa (h-1)^2 x_1^2-2 h \sigma x_1+\kappa \{\sigma(2 h-1)-h+1\}.
\end{array}
\right.
\label{(2-2)}
\end{eqnarray*}
\label{thm1}
\end{thm}
\textbf{Proof}.
The following algebraic equations determine all the equilibria of system (\ref{(3)})
\begin{eqnarray}
\left\{
\begin{array}{l}
x\{\sigma(\kappa-x)(1+h(1+\alpha y)x^2)-k(1+\alpha y)xy\}=0,\\
\kappa y\{(1-h)(1+\alpha y)x^2-1\}=0.
\end{array}
\right.
\label{(2-3)}
\end{eqnarray}
Note that system (\ref{(3)}) always has two boundary equilibria $E_{0}(0,0)$ and $E_{\kappa}(\kappa,0)$ for all permissible parameters.
The positive equilibria of system (\ref{(3)}) lie on the curve
\begin{eqnarray}
y=\sigma x(1-\frac{x}{\kappa}), ~~~0<x<\kappa.
\label{(2-4)}
\end{eqnarray}
We obtain the following function by substituting curve (\ref{(2-4)}) into the second equation in (\ref{(2-3)})
\begin{eqnarray}
F(x):=-\alpha\sigma(h-1)x^4+\kappa\alpha\sigma(h-1)x^3+\kappa(h-1)x^2+\kappa,
\label{(2-5)}
\end{eqnarray}
whose zeros in the interval $(0,\kappa)$ determine the abscissas of  positive equilibria of system (\ref{(3)}).
It is obvious that $F(x)$ has zeros only if $h\neq1$.
Clearly, the derivative of $F(x)$ is
\begin{eqnarray}
F'(x)=(h-1)x(-4\alpha\sigma x^2+3\kappa\alpha\sigma x+2\kappa).
\label{(2-6)}
\end{eqnarray}
The unique positive zero of $F'(x)=0$ is
\begin{eqnarray}
x_*:=\frac{3\kappa\alpha\sigma+\sqrt{\kappa\alpha\sigma(9\kappa\alpha\sigma+32)}}{8\alpha\sigma}.
\label{(2-7)}
\end{eqnarray}
If $h>1$, then $F(x)$ increases first and then decreases for $x\in(0,+\infty)$ and $F(\kappa)=\kappa\{\kappa^2(h-1)+1\}>0$, which implies that $F(x)$ has no zero for $x\in(0,\kappa)$.
If $0<h<1$, then $F(x)$ decreases first and then increases for $x\in(0,+\infty)$.
In order to determine the number of zeros of $F(x)$ in the interval $(0,\kappa)$ for this case, we need
consider the signs of $F(x_*)$ and $F(\kappa)$ and compare the two numbers $x_*$ and $\kappa$ for $0<h<1$ and $(\alpha,\kappa,\sigma)\in\mathbb{R}_+^3$.
The discussion is divided into the following three subcases.
(I) If $F(x_*)>0$, i.e., $0<h<1$ and $(\alpha,\kappa,\sigma)\in\mathcal{P}_1$, then $F(x)$ has no positive zero.
(II) If $F(x_*)=0$ and $x_*\geq\kappa$, i.e., $0<h<1$ and $(\alpha,\kappa,\sigma)\in\mathcal{P}_2$, then $F(x)$ has no positive zero in the interval $(0,\kappa)$;
If $F(x_*)=0$ and $x_*<\kappa$, i.e., $0<h<1$ and $(\alpha,\kappa,\sigma)\in\mathcal{P}_3$, then $F(x)$ has a unique positive zero $x_*$ in the interval $(0,\kappa)$.
(III) If $F(x_*)<0$, $F(\kappa)\geq0$ and $x_*>\kappa$, i.e., $0<h<1$ and $(\alpha,\kappa,\sigma)\in\mathcal{P}_4$, then $F(x)$ has no positive zero in the interval $(0,\kappa)$;
If $F(x_*)<0$ and $F(\kappa)<0$, i.e., $0<h<1$ and $(\alpha,\kappa,\sigma)\in\mathcal{P}_5$, then $F(x)$ has a unique positive zero denoted as $x_1$ in the interval $(0,\kappa)$;
If $F(x_*)<0$, $F(\kappa)=0$ and $x_*<\kappa$, i.e., $0<h<1$ and $(\alpha,\kappa,\sigma)\in\mathcal{P}_6$, then $F(x)$ has a unique positive zero $x_1$ in the interval $(0,\kappa)$;
If $F(x_*)<0$, $F(\kappa)>0$ and $x_*<\kappa$, i.e., $0<h<1$ and $(\alpha,\kappa,\sigma)\in\mathcal{P}_7$, then $F(x)$ has two positive zeros denoted as $x_1$ and $x_2$$(x_1<x_2)$ in the interval $(0,\kappa)$.
Furthermore, $\alpha_1$ and $\alpha_2$$(\alpha_1<\alpha_2)$ given in (\ref{(2-1)}) are two larger positive zeros of
\begin{eqnarray*}
\left.
\begin{array}{l}
f(\alpha):=-27 \kappa^3 \sigma^3 (h-1) \alpha^3-16 \sigma^2 (9 h \kappa^2-9 \kappa^2+16) \alpha^2-4 \sigma \kappa (h-1) (h \kappa^2-\kappa^2+32) \alpha\\
\phantom{f(\alpha):=}
-16 \kappa^2 (h-1)^2,
\end{array}
\right.
\end{eqnarray*}
which determines the sign of $F(x_*)$.
Therefore, system (\ref{(3)}) has no positive equilibrium if $h\geq1$ and $(\alpha,\kappa,\sigma)\in\mathbb{R}_+^3$ or $0<h<1$ and $(\alpha,\kappa,\sigma)\in\mathcal{P}_1\cup\mathcal{P}_2\cup\mathcal{P}_4$;
System (\ref{(3)}) has either a unique positive equilibrium $E_1(x_1,y_1)$ if $0<h<1$ and $(\alpha,\kappa,\sigma)\in\mathcal{P}_5\cup\mathcal{P}_6$
or a unique positive equilibrium $E_*(x_*,y_*)$ if $0<h<1$ and $(\alpha,\kappa,\sigma)\in\mathcal{P}_3$;
System (\ref{(3)}) has two positive equilibria $E_1(x_1,y_1)$ and $E_2(x_2,y_2)$ if $0<h<1$ and $(\alpha,\kappa,\sigma)\in\mathcal{P}_7$, where $y_i=\sigma x_i(1-\frac{x_i}{\kappa})$ $(i=1, 2, *)$.

In the following, we study the qualitative properties of the equilibria. By a direct computation, the Jacobian matrix of system (\ref{(3)}) at any equilibrium has the form
$$
J:=
\left(\begin{array}{cc} J_{11} &
 x^2 (\alpha h \kappa  \sigma x-\alpha h \sigma x^2-2 \alpha \kappa y-\kappa) \\
 2 \kappa (1-h) (\alpha y+1)x y &  \kappa \alpha (1-h)x^2y
\end{array}\right)
$$
with
\begin{eqnarray*}
\begin{array}{l}
J_{11}:=\sigma \{1+h (\alpha y+1) x^2\} (\kappa-2 x)+2 x (h \kappa \sigma x-h \sigma x^2-\kappa y) (\alpha y+1).
\end{array}
\end{eqnarray*}
We use the symbols $T$ and $D$ to denote the trace and determinant of Jacobian matrix $J$ respectively.
It is obvious that $E_0$ is a saddle since $D|_{E_0}=-\sigma\kappa^2<0$.
The computation yields that
\begin{eqnarray*}
\begin{array}{l}
D|_{E_{\kappa}}=\sigma\kappa^2(h\kappa^2+1)\{\kappa^2(h-1)+1\},~~
T|_{E_{\kappa}}=-\kappa\{\sigma(h\kappa^2+1)+\kappa^2(h-1)+1\},\\
T|_{E_{\kappa}}^2-4D|_{E_{\kappa}}=\kappa^2\{\sigma(h\kappa^2+1)+\kappa^2(1-h)-1\}^2,
\end{array}
\end{eqnarray*}
we therefore have $D|_{E_{\kappa}}<0$ if $0<h<1$ and $\kappa>\kappa_1$ implying that $E_{\kappa}$ is a saddle,
$D|_{E_{\kappa}}=0$ if $0<h<1$ and $\kappa=\kappa_1$ implying that $E_{\kappa}$ is degenerate,
$D|_{E_{\kappa}}>0$, $T|_{E_{\kappa}}<0$ and $T|_{E_{\kappa}}^2-4D|_{E_{\kappa}}\geq0$ if $h\geq1$ or $0<h<1$ and $0<\kappa<\kappa_1$ implying that $E_{\kappa}$ is a stable node.
At the positive equilibria $E_{i}$, $i=1, 2, *$,
we use the branch $1+h (\alpha y+1) x^2=(1+\alpha y)x^2$ in $J_{11}$ and the pseudo-division reduction to
obtain determinant $D|_{E_{i}}$ and trace $T|_{E_{i}}$ of the Jacobian matrix $J$ as follows
\begin{eqnarray}
\begin{array}{l}
D|_{E_{i}}=\frac{ \sigma (-\alpha \sigma x_i^2+\alpha \kappa \sigma x_i+\kappa) (x_i-\kappa)x_i^3}{\kappa}F'(x_i),\\
T|_{E_{i}}=\frac{\alpha^2 \sigma^2 \kappa \{-\kappa (h-1)^2 x_i^2-2 h \sigma x_i+\kappa (2 h \sigma-h-\sigma+1)\}}{\kappa \alpha^2 \sigma^2 (1-h)},
\end{array}
\label{(2-8)}
\end{eqnarray}
where $F'(x_i)$ is the derivative of $F(x)$ at $x_i$ given by (\ref{(2-6)}).
It is obvious that the sign of $D|_{E_{i}}$ is opposite to that of $F'(x_i)$ and the sign of $T|_{E_{i}}$ is same as that of $T(x_1)$, which is the factor of  numerator of  $T|_{E_{i}}$ given in Theorem \ref{thm1}.
Because $F'(x_1)<0$, $F'(x_2)>0$ and $F'(x_*)=0$, we obtain  $D|_{E_{1}}>0$, $D|_{E_{2}}<0$ and $D|_{E_{*}}=0$, which imply that $E_{2}$ is a saddle,
$E_{*}$ is degenerate and $E_{1}$ is a stable node or focus (or an unstable node or focus or a center or weak focus) if $T(x_1)<0$ (or $>0$ or $=0$).
\qquad$\Box$

\begin{rem}
System (\ref{(3)}) has more equilibria than the system without cooperative hunting among predators (i.e., $\alpha=0$ in system (\ref{(3)})).
More concretely, the two systems have the same boundary equilibria $E_{0}$ and $E_{\kappa}$ and they have the same qualitative properties as shown in Theorem \ref{thm1}.
Nevertheless, the two systems have different positive equilibrium. System (\ref{(3)}) has at most two positive equilibria as in Theorem \ref{thm1},
but the system without cooperative hunting has a unique positive equilibrium $(\kappa_1,\sigma\kappa_1(1-\frac{\kappa_1}{\kappa}))$ for $\kappa>\kappa_1$ and $0<h<1$, which is a stable node or focus (or an unstable node or focus or a center or weak focus) if $\kappa>\kappa_1$ and $0<h\leq\frac{1}{2}$ or $\kappa_1<\kappa<\frac{2h\kappa_1}{2h-1}$ and $\frac{1}{2}<h<1$ (or $\kappa>\frac{2h\kappa_1}{2h-1}$ and $\frac{1}{2}<h<1$ or $\kappa=\frac{2h\kappa_1}{2h-1}$ and $\frac{1}{2}<h<1$),
and has no positive equilibrium if  $\kappa>0$ and $h\geq1$ or $0<\kappa\leq\kappa_1$ and $0<h<1$.
Therefore, the main effect of the cooperative hunting among predators is on the coexistence of species.
\label{rem1}
\end{rem}

\section{Bifurcations }
From Theorem \ref{thm1}, we know that system (\ref{(3)}) has two degenerate equilibria $E_{\kappa}$ and $E_*$ and a center or weak focus $E_1$. In this section, we further identify the topological types of these nonhyperbolic equilibria and display all bifurcations at them. Concretely, system (\ref{(3)}) may exhibit transcritical and pitchfork bifurcations around equilibrium $E_{\kappa}$, saddle-node and Bogdanov-Takens bifurcations around  equilibrium $E_*$ and Hopf bifurcation around  equilibrium  $E_1$.

\subsection{Transcritical and Pitchfork Bifurcations}

In this subsection, we identify the topological type of the degenerate equilibrium $E_{\kappa}$  and
show that both transcritical and pitchfork bifurcations may occur at $E_{\kappa}$. In fact,
Theorem \ref{thm1} shows that $E_{\kappa}$ is degenerate if $\kappa=\kappa_1$ and $0<h<1$ with $D|_{E_{\kappa}}=0$ and $T|_{E_{\kappa}}<0$.
The following theorem indicates that $E_{\kappa}$ is a saddle-node and system (\ref{(3)}) may undergo transcritical and pitchfork bifurcations at $E_{\kappa}$.

\begin{thm}
For $\kappa=\kappa_1$ and $0<h<1$, $E_{\kappa}$ is a saddle-node of system (\ref{(3)}). Moreover,
(i) when $\alpha\neq\frac{2}{\sigma\kappa_1}$, a transcritical bifurcation happens at $E_{\kappa}$ if $\kappa$ varies from $\kappa>\kappa_1$ to $\kappa<\kappa_1$;
(ii) when  $\alpha=\frac{2}{\sigma\kappa_1}$, a pitchfork bifurcation
happens at $E_{\kappa}$ if $\kappa$ varies from $\kappa>\kappa_1$ to $\kappa<\kappa_1$.
\label{thm2}
\end{thm}
\textbf{Proof}.
Let $\epsilon=\kappa-\kappa_1$ and consider $|\epsilon|$ sufficiently small.
Suspending  system (\ref{(3)}) with the variable $\epsilon$ and
using the linear transformation $x=u+v+\kappa_1+\epsilon$ and $y=-\sigma v$  together with  the time-rescaling $\tau=\frac{\kappa_1\sigma}{h-1} t$ to translate $E_{\kappa}$ to $(0, 0)$  and diagonalize the linear part of the suspended system, we obtain
\begin{eqnarray}
\left\{
\begin{array}{l}
\frac{du}{dt}=-\frac{(\alpha \kappa_1 \sigma-2) (h-1) }{\sigma \kappa_1}u^2+\frac{2 (h-1)}{\sigma \kappa_1}u v+\frac{2 (h-1)}{\sigma \kappa_1} u\epsilon
-\frac{(2 \alpha \sigma+h \kappa_1-\kappa_1) (h-1) }{\sigma \kappa_1}u^3\\
\phantom{\frac{du}{dt}=}
-\frac{2 (\alpha \sigma+h \kappa_1-\kappa_1) (h-1) }{\sigma \kappa_1}u^2 v
-\frac{(3 \alpha \sigma+4 h \kappa_1-4 \kappa_1) (h-1) }{\sigma \kappa_1}u^2 \epsilon
-\frac{(h-1)^2 }{\sigma}u v^2
-\frac{4 (h-1)^2 }{\sigma} u v \epsilon\\
\phantom{\frac{du}{dt}=}
-\frac{3 (h-1)^2 }{\sigma}u \epsilon^2
+O(\parallel(u,v,\epsilon)\parallel^4),\\

\frac{dv}{dt}=v-\frac{\alpha \sigma^2 (h-1) \kappa_1^3+\alpha \sigma (h \sigma-h+1) \kappa_1-2 h \sigma+2 h+\sigma-2}{\kappa_1 \sigma}u^2+\frac{2 h+1}{\kappa_1}v^2-\frac{3 (\sigma-1) (h-1)^2 }{\sigma}u \epsilon^2\\
\phantom{\frac{du}{dt}=}
+\frac{2 (\sigma-1) (h-1) }{\sigma \kappa_1} u\epsilon
-\frac{(h-1) (3 h \sigma-h-\sigma+1) \kappa_1+\alpha \sigma (3 h \sigma-2 h-2 \sigma+2) }{\sigma \kappa_1}u^3\\
\phantom{\frac{du}{dt}=}
-\frac{\alpha h \kappa_1 \sigma^2-4 h \sigma+2 h-2}{\kappa_1 \sigma}u v
-\frac{(h-1) (9 h \sigma-2 h-2 \sigma+2) \kappa_1+2 \alpha \sigma (3 h \sigma-h-\sigma+1)}{\sigma \kappa_1}u^2 v\\
\phantom{\frac{du}{dt}=}
-\frac{(h-1) \{(6 h \sigma-4 h-4 \sigma+4) \kappa_1+3 \alpha \sigma (\sigma-1)\}}{\kappa_1 \sigma}u^2 \epsilon
-\frac{(h-1) (9 h \sigma-h-\sigma+1) \kappa_1+3 \alpha h \sigma^2}{\kappa_1 \sigma}u v^2\\
\phantom{\frac{du}{dt}=}
+\frac{2 h+1}{\kappa_1}v \epsilon
-\frac{4 (h-1) (3 h \sigma-h-\sigma+1) \kappa_1+3 \alpha h \sigma^2}{\kappa_1 \sigma} u v \epsilon
-3 h (h-1) v^3-6 h (h-1) v^2 \epsilon\\
\phantom{\frac{du}{dt}=}
-3 h (h-1) v \epsilon^2
+O(\parallel(u,v,\epsilon)\parallel^4),\\
\frac{d\epsilon}{dt}=0,
\end{array}
\right.
\label{(3-1)}
\end{eqnarray}
where $\tau$ is still denoted as $t$.
Theorem 1 of \cite{Carr} shows that
system (\ref{(3-1)}) has a $C^\infty$ center manifold $W^c: v=h_1(u,\epsilon)$ near $(u,v,\epsilon)=(0,0,0)$,
which is tangent to the plane $v=0$ at $(u,v,\epsilon)=(0,0,0)$. Let
\begin{equation}
v=h_1(u,\epsilon)=a_1u^2+b_1\epsilon^2+c_1u\epsilon+O(\parallel(u,\epsilon)\parallel^3).
\label{(3-2)}
\end{equation}
By the invariant property of  center manifold (\ref{(3-2)}) to the solutions of (\ref{(3-1)}),
we differentiate both sides of (\ref{(3-2)}) and obtain
$\dot{v}=h_{1u}\dot{u}+h_{1\epsilon}\dot{\epsilon}$.
Substituting (\ref{(3-1)}) in the equality
and comparing the coefficients
of $u^2$, $\epsilon^2$ and $u\epsilon$, we get $b_1=0$ and
\begin{eqnarray*}
\begin{array}{l}
a_1=\frac{\alpha \sigma^2 (h-1) \kappa_1^3+\alpha \sigma (h \sigma-h+1) \kappa_1-2 h \sigma+2 h+\sigma-2}{\kappa_1 \sigma},~~
c_1=-\frac{2 (\sigma-1) (h-1)}{\kappa_1 \sigma}.
\end{array}
\end{eqnarray*}
Thus, restricted to the center manifold (\ref{(3-2)}), system (\ref{(3-1)}) can be written as
\begin{eqnarray}
\begin{array}{l}
\frac{du}{dt}=\frac{2 (h-1) }{\kappa_1 \sigma} u\epsilon-\frac{(\alpha \kappa_1 \sigma-2) (h-1) }{\sigma \kappa_1}u^2
-\frac{(h-1) \{-4 (h-1) (h \sigma-h-2 \sigma+1) \kappa_1+3 \alpha \sigma^2\}}{\sigma^2 \kappa_1} u^2 \epsilon\\
\phantom{\frac{du}{dt}=}
-\frac{3 (h-1)^2}{\sigma}u \epsilon^2
+\frac{(h-1) \{(h-1) (4 h \sigma-4 h-3 \sigma+4) \kappa_1+2 \alpha \sigma (h \sigma-h-2 \sigma+1)\} }{\sigma^2 \kappa_1}u^3
+O(\parallel(u,\epsilon)\parallel^4).
\end{array}
\label{(3-3)}
\end{eqnarray}

If $\alpha\neq\frac{2}{\sigma\kappa_1}$, then the coefficient $\frac{(\alpha \kappa_1 \sigma-2) (h-1) }{\sigma \kappa_1}\neq0$ in system (\ref{(3-3)}), the origin is the unique equilibrium of system (\ref{(3-3)})
as $\epsilon=0$ and the other  equilibrium arises from the origin as $\epsilon\neq0$. Therefore, the boundary equilibrium $E_{\kappa}$ is a saddle-node as $\kappa=\kappa_1$ and system (\ref{(3)}) undergoes a transcritical bifurcation(\cite{Guckenheimer}) at $E_{\kappa}$. Furthermore,
if $\alpha<\frac{2}{\sigma\kappa_1}$ and $\kappa$ varies from $\kappa<\kappa_1$ to $\kappa>\kappa_1$, then a stable node $E_{\kappa}$ changes into a saddle $E_{\kappa}$ and a stable node $E_1$;
if $\alpha>\frac{2}{\sigma\kappa_1}$ and $\kappa$ varies from $\kappa>\kappa_1$ to $\kappa<\kappa_1$, then a saddle $E_{\kappa}$ changes into a stable node $E_{\kappa}$ and a saddle $E_2$.
If $\alpha=\frac{2}{\sigma\kappa_1}$, then $\frac{(\alpha \kappa_1 \sigma-2) (h-1) }{\sigma \kappa_1}=0$ and the coefficient of $u^3$ is $\frac{5 \sigma (1-h)}{ \sigma^2\kappa_1^2}\neq0$ in (\ref{(3-3)}) and the origin is the unique equilibrium
as $\epsilon=0$ and the other two equilibria arise from the origin as $\epsilon>0$.
Therefore, $E_{\kappa}$ is a saddle-node as $\kappa=\kappa_1$ and  system (\ref{(3)}) undergoes a  pitchfork bifurcation(\cite{Guckenheimer}) at $E_{\kappa}$ as $\kappa$ varies from $\kappa<\kappa_1$ to $\kappa>\kappa_1$ such that a stable node $E_{\kappa}$ changes into a saddle $E_{\kappa}$ and a stable node $E_1$.  This completes the proof of the theorem.
\qquad$\Box$

\subsection{Saddle-node Bifurcation}
In this subsection, we prove the topological type of the degenerate equilibrium $E_*$  and
show that the saddle-node bifurcation may occur at $E_*$ if $D|_{E_*}=0$ and $T|_{E_*}\neq0$.
From $F(x_*)=F'(x_*)=0$, we can express $\alpha$ and $\kappa$ by $x_*$, $\sigma$ and $h$ as follows
\begin{eqnarray}
\begin{array}{l}
\alpha=\alpha_3:=\frac{2\{(h-1)x_*^2+2\}}{(1-h)\sigma x_*^3},~~~ \kappa=\kappa_2:=\frac{2x_*\{(h-1)x_*^2+2\}}{(h-1)x_*^2+3}
\end{array}
\label{(3-4)}
\end{eqnarray}
with $0<x_*\leq 1$ and $0<h<1$ or $x_*>1$ and $\frac{x_*^2-1}{x_*^2}<h<1$.
Substituting (\ref{(3-4)}) into $T|_{E_*}$, we obtain that
$T|_{E_*}<0$ if $\sigma>\sigma_1$, $T|_{E_*}>0$ if $\sigma<\sigma_1$ and $T|_{E_*}=0$ if $\sigma=\sigma_1$  with
\begin{eqnarray*}
\begin{array}{l}
\sigma_1:=\frac{(h-1)\{(h-1)x_*^2+2\}\{(h-1)x_*^2+1\}}{(h-1)^2x_*^2+h-2}.
\end{array}
\end{eqnarray*}

\begin{thm}
For $\alpha=\alpha_3$, $\kappa=\kappa_2$ and $\sigma\neq\sigma_1$ with $0<x_*\leq 1$ and $0<h<1$ or $x_*>1$ and $\frac{x_*^2-1}{x_*^2}<h<1$,
$E_*$ is a saddle-node of system (\ref{(3)}) and a saddle-node bifurcation occurs at $E_*$ if $\alpha$ varies from  $\alpha<\alpha_3$ to $\alpha>\alpha_3$.
\label{thm3}
\end{thm}
\textbf{Proof}.
Let $\epsilon=\alpha-\alpha_3$ and consider $|\epsilon|$ sufficiently small.
Translating $E_{*}$ to the origin $(0,0)$  and suspending  system (\ref{(3)}) with the variable $\epsilon$ we obtain the following system
\begin{eqnarray*}
\left\{
\begin{array}{l}
\frac{dx}{dt}=a_{100} x+a_{010} y+a_{001} \epsilon+a_{200} x^2+a_{110} x y+a_{020} y^2+a_{101} x \epsilon+a_{011} y \epsilon\\
\phantom{\frac{dx}{dt}=}
+O(\parallel(x,y,\epsilon)\parallel^3),\\
\frac{dy}{dt}=b_{100} x+b_{010} y+b_{001} \epsilon+b_{200} x^2+b_{110} x y+b_{101} x \epsilon+b_{011} \epsilon y+b_{020} y^2\\
\phantom{\frac{dx}{dt}=}
+O(\parallel(x,y,\epsilon)\parallel^3),\\
\frac{d\epsilon}{dt}=0,
\end{array}
\right.
\end{eqnarray*}
where the coefficients $a_{ijk}$ and $b_{ijk}$ are listed as follows.
\begin{eqnarray*}
\left.
\begin{array}{l}
a_{100} := -\frac{2 x_* \sigma \{(h-1)^2 x_*^2+h-2\}}{(h x_*^2-x_*^2+3) (h-1)},~
a_{010} := -\frac{2 x_* (h x_*^2-x_*^2+2) \{(h-1)^2 x_*^2+h-2\}}{(h-1) (h x_*^2-x_*^2+3)},\\
a_{001} := \frac{(h-1) (h x_*^2-x_*^2+1)^2 \sigma^2 x_*^5}{2(h x_*^2-x_*^2+3) (h x_*^2-x_*^2+2)},~
a_{200} := -\frac{\sigma \{(h-1) (h-2) x_*^2-3 h-4\}}{(h-1) (h x_*^2-x_*^2+3)},\\
a_{110} := -\frac{4 (h x_*^2-x_*^2+2) \{-2+(h-1)^2 x_*^2\}}{(h-1) (h x_*^2-x_*^2+3)},~
a_{101} := \frac{(h x_*-x_*+1) (h x_*-x_*-1) (h x_*^2-x_*^2+1) \sigma^2 x_*^4}{(h x_*^2-x_*^2+3) (h x_*^2-x_*^2+2)},\\
a_{020} := \frac{4 (h x_*^2-x_*^2+2)^2}{(h-1) \sigma (h x_*^2-x_*^2+3)},~
a_{011} := \frac{(h-2) (h x_*^2-x_*^2+1) \sigma x_*^4}{h x_*^2-x_*^2+3},~
b_{100} := -\frac{2 (-h x_*^2+x_*^2-1) \sigma x_*}{h x_*^2-x_*^2+3},\\
b_{010} := \frac{2 x_* (h x_*^2-x_*^2+2) (h x_*^2-x_*^2+1)}{h x_*^2-x_*^2+3},~
b_{001} := -\frac{(h-1) (h x_*^2-x_*^2+1)^2 \sigma^2 x_*^5}{2(h x_*^2-x_*^2+3) (h x_*^2-x_*^2+2)},\\
b_{200} := -\frac{(-h x_*^2+x_*^2-1) \sigma}{h x_*^2-x_*^2+3},~
b_{110} := \frac{4 (h x_*^2-x_*^2+2)^2}{h x_*^2-x_*^2+3},~
b_{101} := -\frac{(h-1) (h x_*^2-x_*^2+1)^2 \sigma^2 x_*^4}{(h x_*^2-x_*^2+2) (h x_*^2-x_*^2+3)},\\
b_{020} :=\frac{ 4 (h x_*^2-x_*^2+2)^2}{\sigma (h x_*^2-x_*^2+3)},~
b_{011} := -\frac{2 (h-1) (h x_*^2-x_*^2+1) \sigma x_*^4}{h x_*^2-x_*^2+3}.
\end{array}
\right.
\end{eqnarray*}
Since $\sigma\neq\sigma_1$ as well as $0<x_*\leq 1$ and $0<h<1$ or $x_*>1$ and $\frac{x_*^2-1}{x_*^2}<h<1$, we get
\begin{eqnarray*}
\left.
\begin{array}{l}
a_{100}+b_{010}=\frac{2 x_* \{(-(h-1)^2 x_*^2-h+2) \sigma+(h-1) (h x_*^2-x_*^2+2) (h x_*^2-x_*^2+1)\}}{(h x_*^2-x_*^2+3) (h-1)}\neq0.
\end{array}
\right.
\end{eqnarray*}
Therefore,  by applying the linear transformation
$x=u+\frac{a_{100}}{b_{100}}v-\frac{b_{001}(b_{010}-b_{100})}{b_{100}(a_{100}+b_{010})}\epsilon$ and $y=-\frac{b_{100}}{b_{010}}u+v$
to diagonalize the linear part, we can change the suspended system into the following system
\begin{eqnarray}
\left\{
\begin{array}{l}
\frac{du}{dt}=p_{001}\epsilon+p_{200}u^2+p_{020}v^2+p_{002}\epsilon^2+p_{110}uv+p_{011} v\epsilon+p_{101} u\epsilon\\
\phantom{\frac{du}{dt}=}+O(\parallel(u,v,\epsilon)\parallel^3),\\
\frac{dv}{dt}=q_{010}v+q_{200}u^2+q_{020}v^2+q_{002}\epsilon^2+q_{110}uv+q_{101}u\epsilon+q_{011}v\epsilon\\
\phantom{\frac{dx}{dt}=}+O(\parallel(u,v,\epsilon)\parallel^3),\\
\displaystyle \frac{d\epsilon}{dt}=0,
\end{array}
\right.
\label{(3-5)}
\end{eqnarray}
where $p_{ijk}$ and $q_{ijk}$ are listed  in the Appendix.
System (\ref{(3-5)}) similarly has a two-dimensional center manifold $W^c: v=h_2(u,\epsilon)$ near $(u,v,\epsilon)=(0,0,0)$ as follows
\begin{eqnarray}
v=h_2(u,\epsilon)=a_2u^2+b_2\epsilon^2+c_2u\epsilon+O(\parallel(u,\epsilon)\parallel^3)
\label{(3-6)}
\end{eqnarray}
with coefficients $a_2$, $b_2$ and $c_2$  given in the Appendix.
Restricted to the center manifold (\ref{(3-6)}), system (\ref{(3-5)}) can be written as
\begin{eqnarray}
\begin{array}{l}
\frac{du}{dt}=d_0(\epsilon)+d_1(\epsilon)u+d_2(\epsilon)u^2+O(|u|^3),
\end{array}
\label{(3-7)}
\end{eqnarray}
where
\begin{eqnarray*}
\begin{array}{l}
d_0(\epsilon)=\frac{(h x_*^2-x_*^2+1)^2 (h-1) \sigma^2 x_*^5}{2\{(-(h-1)^2 x_*^2-h+2) \sigma+(h-1) (h x_*^2-x_*^2+2) (h x_*^2-x_*^2+1)\} (h x_*^2-x_*^2+3)} \epsilon+O(|\epsilon|^2),\\
d_1(\epsilon)=\frac{(h-1) (h x_*^2-x_*^2+1) \sigma^2 x_*^4 }{2\{(-(h-1)^2 x_*^2-h+2) \sigma+(h-1) (h x_*^2-x_*^2+2) (h x_*^2-x_*^2+1)\}^2 (h x_*^2-x_*^2+3)}
\{(-(h-1)^4 x_*^6\\
\phantom{d_1(\epsilon)=}
-8 (h-1)^3 x_*^4-(h-1) (11 h-13) x_*^2-4 h+2) \sigma+(h-1) (h x_*^2-x_*^2+2)\\
\phantom{d_1(\epsilon)=}
\cdot (h x_*^2-x_*^2+1) ((h-1)^2 x_*^4+7(h-1) x_*^2+4)\}
\epsilon+O(|\epsilon|^2),\\
d_2(\epsilon)=\frac{(h x_*^2-x_*^2+2) (h x_*^2-x_*^2+6) (h x_*^2-x_*^2+1)\sigma }{\{(-(h-1)^2 x_*^2-h+2) \sigma+(h-1) (h x_*^2-x_*^2+2) (h x_*^2-x_*^2+1)\} (h x_*^2-x_*^2+3)} +O(|\epsilon|).
\end{array}
\end{eqnarray*}
We can check that $d_0(0)=d_1(0)=0$ and $d_2(0)\neq0$ for $\sigma\neq\sigma_1$ as well as $0<x_*\leq 1$ and $0<h<1$ or $x_*>1$ and $\frac{x_*^2-1}{x_*^2}<h<1$.
Concretely, we have $d_2(0)>0$ for $\sigma>\sigma_1$ and  $d_2(0)<0$ for $\sigma<\sigma_1$.
Furthermore, making the translation $u=w-\frac{d_1(\epsilon)}{2d_2(\epsilon)}$ and the time-rescaling $\tau:=d_2(\epsilon)t$ to system (\ref{(3-7)}), we obtain
\begin{eqnarray}
\begin{array}{l}
\frac{dw}{d\tau}=\zeta(\epsilon)+w^2+O(|w|^3),
\end{array}
\label{(3-8)}
\end{eqnarray}
where
$\zeta(\epsilon):=\frac{4d_0(\epsilon)d_2(\epsilon)-d_1^2(\epsilon)}{4d_2^2(\epsilon)}$.
The computation shows
$\zeta(0)=0$ and
\begin{eqnarray*}
\begin{array}{l}
\zeta'(0)=\frac{d_0'(0)}{d_2(0)}
=\frac{(h-1) (h x_*^2-x_*^2+1) \sigma x_*^5}{2(h x_*^2-x_*^2+6) (h x_*^2-x_*^2+2)}<0
\end{array}
\end{eqnarray*}
for $\sigma\neq\sigma_1$ as well as $0<x_*\leq 1$ and $0<h<1$ or $x_*>1$ and $\frac{x_*^2-1}{x_*^2}<h<1$.
Hence,  the origin is the equilibrium of (\ref{(3-8)}) for $\epsilon=0$ and two nonzero  equilibria arise as $\epsilon>0$. Therefore, for $\alpha=\alpha_3$, $\kappa=\kappa_2$ and $\sigma\neq\sigma_1$ with $0<x_*\leq 1$ and $0<h<1$ or $x_*>1$ and $\frac{x_*^2-1}{x_*^2}<h<1$, equilibrium $E_*$ is a saddle-node and
the node $E_1$ and the saddle $E_2$ arise from the saddle-node bifurcation(\cite{Guckenheimer}) at $E_*$  as $\alpha$ changes from $\alpha<\alpha_3$ to $\alpha>\alpha_3$. The proof is completed.
\qquad$\Box$


\subsection{Bogdanov-Takens Bifurcation}
As shown in subsection 3.2, the equilibrium $E_{*}$ is degenerate with $D|_{E_{*}}=0$ and $T|_{E_{*}}=0$ if $\alpha=\alpha_*:=\alpha_3$, $\sigma=\sigma_*:=\sigma_1$ and $\kappa=\kappa_2$ with $0<x_*\leq 1$ and $0<h<1$ or $x_*>1$ and $\frac{x_*^2-1}{x_*^2}<h<1$.
The following theorem shows that the codimension of the degenerate equilibrium $E_{*}$ is two (i.e., a cusp of codimension two) and system (\ref{(3)}) may exhibit the Bogdanov-Takens bifurcation of codimension two
at $E_{*}$ under a small parameter perturbation around $(\alpha_*,\sigma_*)$.

\begin{thm}
If $\alpha=\alpha_*$, $\sigma=\sigma_*$ and $\kappa=\kappa_2$ with $0<x_*\leq 1$ and $0<h<1$ or $x_*>1$ and $\frac{x_*^2-1}{x_*^2}<h<1$,
then the equilibrium $E_{*}$ is a degenerate equilibrium of codimension two (i.e., a cusp of codimension two) and system (\ref{(3)}) undergoes a Bogdanov-Takens bifurcation of codimension two around $E_*$ as the pair of parameters $(\alpha, \sigma)$ varies near $(\alpha_*, \sigma_*)$. Concretely, there is a small neighborhood $U$ of $(\alpha_*, \sigma_*)$ in the $(\alpha, \sigma)$-parameter space and four curves
\begin{eqnarray*}
\begin{array}{l}
\mathcal{SN}^+:=\{(\alpha, \sigma)\in U:\alpha=
\alpha_*+f_{11}(\sigma-\sigma_*)+f_{12}(\sigma-\sigma_*)^2+O(|\sigma-\sigma_*|^3), \sigma<\sigma_*\},\\
\mathcal{SN}^-:=\{(\alpha, \sigma)\in U:\alpha=
\alpha_*+f_{11}(\sigma-\sigma_*)+f_{12}(\sigma-\sigma_*)^2+O(|\sigma-\sigma_*|^3), \sigma>\sigma_*\},\\
\mathcal{H}:=\{(\alpha, \sigma)\in U:
\alpha=
\alpha_*+f_{11}(\sigma-\sigma_*)+f_{22}(\sigma-\sigma_*)^2+O(|\sigma-\sigma_*|^3), \sigma>\sigma_*
\},\\
\mathcal{HL}:=\{(\alpha, \sigma)\in U:
\alpha=\alpha_*+f_{11}(\sigma-\sigma_*)+f_{32}(\sigma-\sigma_*)^2+O(|\sigma-\sigma_*|^3)
, \sigma>\sigma_*
\},
\end{array}
\end{eqnarray*}
where $f_{11}$,  $f_{12}$,  $f_{22}$ and $f_{32}$ are given in the Appendix, such that system (\ref{(3)})
undergoes the saddle-node, the Hopf and the homoclinic bifurcations around $E_*$ as $(\alpha, \sigma)$ passes across curves
$\mathcal{SN}^+\cup \mathcal{SN}^-$, $\mathcal{H}$ and $\mathcal{HL}$ respectively.
\label{thmBT}
\end{thm}
\textbf{Proof}.
Let $(\epsilon_1, \epsilon_2):=(\alpha-\alpha_*, \sigma-\sigma_*)$ be a parameter vector near $(0,0)$.
Translating $E_*$ of system (\ref{(3)}) to the origin, using the linear transformations
\begin{eqnarray*}
\left.
\begin{array}{l}
x=\frac{-h^2 x_*^2+(2 x_*^2-1) h-x_*^2+2}{(h-1) (h x_*^2-x_*^2+1)}u+v,~~~
y=u
\end{array}
\right.
\end{eqnarray*}
and time-rescaling
\begin{eqnarray*}
\left.
\begin{array}{l}
t=\frac{\{-h^2 x_*^2+(2 x_*^2-1) h-x_*^2+2\} (h x_*^2-x_*^2+3) }{2 x_* (-h x_*^2+x_*^2-2) (h x_*^2-x_*^2+1)^2 (h-1)} \tau
\end{array}
\right.
\end{eqnarray*}
and expanding the system at origin,  we obtain
\begin{eqnarray}
\left\{
\begin{array}{l}
\frac{du}{dt}=c_{00}(\epsilon_1, \epsilon_2)+c_{10}(\epsilon_1, \epsilon_2)u+c_{01}(\epsilon_1, \epsilon_2)v+c_{11}(\epsilon_1, \epsilon_2)uv+c_{02}v^2\\
\phantom{\frac{du}{dt}=}
+c_{20}(\epsilon_1, \epsilon_2)u^2+O(\parallel(u,v)\parallel^3),\\
\frac{dv}{dt}=d_{00}(\epsilon_1, \epsilon_2)+d_{10}(\epsilon_1, \epsilon_2)u+d_{01}(\epsilon_1, \epsilon_2)v+d_{11}(\epsilon_1, \epsilon_2)uv
+d_{02}(\epsilon_1, \epsilon_2)v^2\\
\phantom{\frac{du}{dt}=}
+d_{20}(\epsilon_1, \epsilon_2)u^2+O(\parallel(u,v)\parallel^3),
\end{array}
\right.
\label{(3-9)}
\end{eqnarray}
where $\tau$ is still denoted as $t$ and these coefficients $c_{ij}$ and $d_{ij}$ are given in the Appendix.
Additionally,  the following near-identity transformation
\begin{eqnarray*}
\left.
\begin{array}{l}
u_1:=u,\\
v_1:=c_{00}(\epsilon_1, \epsilon_2)+c_{10}(\epsilon_1, \epsilon_2)u+c_{01}(\epsilon_1, \epsilon_2)v+c_{11}(\epsilon_1, \epsilon_2)uv+c_{02}v^2\\
\phantom{\frac{du}{dt}=}
+c_{20}(\epsilon_1, \epsilon_2)u^2+O(\parallel(u,v)\parallel^3)
\end{array}
\right.
\end{eqnarray*}
brings system (\ref{(3-9)}) to the Kukles from
\begin{eqnarray}
\left\{
\begin{array}{l}
\frac{du_1}{dt}=v_1,\\
\frac{dv_1}{dt}=e_{00}(\epsilon_1, \epsilon_2)+ e_{01}(\epsilon_1, \epsilon_2)v_1+ e_{10}(\epsilon_1, \epsilon_2)u_1+e_{20}(\epsilon_1, \epsilon_2)u_1^2\\
\phantom{\frac{du_1}{dt}=}
+e_{02}(\epsilon_1, \epsilon_2) v_1^2+e_{11}(\epsilon_1, \epsilon_2) u_1 v_1+O(\parallel(u_1,v_1)\parallel^3),
\end{array}
\right.
\label{(3-10)}
\end{eqnarray}
where these coefficients $e_{ij}$ are given in the Appendix. Since the calculation yields
\begin{eqnarray*}
\left.
\begin{array}{l}
e_{11}(0,0)=\frac{\{3 h^3 x_*^4-2 x_*^2 (4 x_*^2-3) h^2+(x_*^2-1) (7 x_*^2-3) h-2 x_*^2 (x_*^2-2)\} \{-h^2 x_*^2+(2 x_*^2-1) h-x_*^2+2\}}{x_* (h x_*^2-x_*^2+1)^3 (h-1)^2}>0
\end{array}
\right.
\end{eqnarray*}
under the conditions of Theorem \ref{thmBT}, we use the transformation
\begin{eqnarray*}
\left.
\begin{array}{l}
u_2:=u_1+\frac{e_{01}(\epsilon_1, \epsilon_2)}{e_{11}(\epsilon_1, \epsilon_2)},~~~
v_2:=v_1
\end{array}
\right.
\end{eqnarray*}
to eliminate the term $e_{01}(\epsilon_1, \epsilon_2)v_1$ in system (\ref{(3-10)}) and obtain
\begin{eqnarray}
\left\{
\begin{array}{l}
\frac{du_2}{dt}=v_2,\\
\frac{dv_2}{dt}=f_{00}(\epsilon_1, \epsilon_2)+ f_{10}(\epsilon_1, \epsilon_2)u_2+e_{20}(\epsilon_1, \epsilon_2)u_2^2+e_{02}(\epsilon_1, \epsilon_2)v_2^2\\
\phantom{\frac{du_2}{dt}=}
+e_{11}(\epsilon_1, \epsilon_2) u_2 v_2+O(\parallel(u_2,v_2)\parallel^3),
\end{array}
\right.
\label{(3-11)}
\end{eqnarray}
where
\begin{eqnarray*}
\left.
\begin{array}{l}
f_{00}(\epsilon_1, \epsilon_2):=e_{00}(\epsilon_1, \epsilon_2)-\frac{e_{10}(\epsilon_1, \epsilon_2) e_{01}(\epsilon_1, \epsilon_2)}{e_{11}(\epsilon_1, \epsilon_2)}+\frac{e_{20}(\epsilon_1, \epsilon_2) e_{01}^2(\epsilon_1, \epsilon_2)}{e_{11}^2(\epsilon_1, \epsilon_2)},\\
f_{10}(\epsilon_1, \epsilon_2):=e_{10}(\epsilon_1, \epsilon_2)-\frac{2 e_{20}(\epsilon_1, \epsilon_2) e_{01}(\epsilon_1, \epsilon_2)}{e_{11}(\epsilon_1, \epsilon_2)}.
\end{array}
\right.
\end{eqnarray*}
Under the near-identity and time-rescaling transformations
\begin{eqnarray*}
\left.
\begin{array}{l}
u_3:=u_2,~~~v_3:=v_2-e_{02}(\epsilon_1, \epsilon_2)u_2v_2, ~~~\tau:=(1+e_{02}(\epsilon_1, \epsilon_2)u_3)t,
\end{array}
\right.
\end{eqnarray*}
 system (\ref{(3-11)}) becomes
\begin{eqnarray}
\left\{
\begin{array}{l}
\frac{du_3}{d\tau}=v_3,\\
\frac{dv_3}{d\tau}=\mu_1(\epsilon_1, \epsilon_2)+\mu_2(\epsilon_1, \epsilon_2)u_3+A(\epsilon_1, \epsilon_2)u_3^2+B(\epsilon_1, \epsilon_2)u_3v_3,
\end{array}
\right.
\label{(3-12)}
\end{eqnarray}
where
\begin{eqnarray*}
\begin{array}{l}
\mu_1(\epsilon_1, \epsilon_2):=f_{00}(\epsilon_1,\epsilon_2),\\
\mu_2(\epsilon_1, \epsilon_2):=f_{10}(\epsilon_1,\epsilon_2)-2e_{02}(\epsilon_1,\epsilon_2)f_{00}(\epsilon_1,\epsilon_2),\\
A(\epsilon_1, \epsilon_2):=e_{20}(\epsilon_1,\epsilon_2)+2e_{02}(\epsilon_1,\epsilon_2)(e_{02}(\epsilon_1,\epsilon_2)f_{00}(\epsilon_1,\epsilon_2)-f_{10}(\epsilon_1,\epsilon_2)),\\
B(\epsilon_1, \epsilon_2):=e_{11}(\epsilon_1,\epsilon_2).
\end{array}
\end{eqnarray*}
Since
\begin{eqnarray*}
\begin{array}{l}
A(0, 0)=\frac{\{h^2 x_*^2+(-2 x_*^2+1) h+x_*^2-2\}^2 (h x_*^2-x_*^2+6)}{2 x_* (h-1)^3 (h x_*^2-x_*^2+1)^3}<0
\end{array}
\end{eqnarray*}
and $B(0, 0)=e_{11}(0,0)>0$ under the conditions of Theorem \ref{thmBT},
with the rescaling
\begin{eqnarray*}
\begin{array}{l}
u_4:=\frac{B^2(\epsilon_1, \epsilon_2)}{A(\epsilon_1, \epsilon_2)}u_3, ~~ v_4:=-\frac{B^3(\epsilon_1, \epsilon_2)}{A^2(\epsilon_1, \epsilon_2)}v_3 ~\mbox{and}~ t:=-\frac{A(\epsilon_1, \epsilon_2)}{B(\epsilon_1, \epsilon_2)}\tau,
\end{array}
\end{eqnarray*}
system (\ref{(3-12)}) can be reduced to the following system
\begin{eqnarray}
\left\{
\begin{array}{l}
\frac{du_4}{dt}=v_4,\\
\frac{dv_4}{dt}=\beta_1(\epsilon_1, \epsilon_2)+\beta_2(\epsilon_1, \epsilon_2)u_4+u_4^2-u_4v_4,
\end{array}
\right.
\label{(3-13)}
\end{eqnarray}
where
\begin{eqnarray}
\begin{array}{l}
\beta_1(\epsilon_1, \epsilon_2):=\frac{B^4(\epsilon_1, \epsilon_2)}{A^3(\epsilon_1, \epsilon_2)}\mu_1(\epsilon_1, \epsilon_2),~~
\beta_2(\epsilon_1, \epsilon_2):=\frac{B^2(\epsilon_1, \epsilon_2)}{A^2(\epsilon_1, \epsilon_2)}\mu_2(\epsilon_1, \epsilon_2).
\end{array}
\label{(3-14)}
\end{eqnarray}
Since $\mu_1(0, 0)=0$ and $\mu_2(0, 0)=0$, we can get that $\beta_1(0, 0)=0$ and $\beta_2(0, 0)=0$. Therefore, the equilibrium $E_*$ is a degenerate equilibrium of codimension two (i.e., a cusp of codimension two).
Moreover, since the Jacobian  determinant of (\ref{(3-14)}) at $(0,0)$ is given by
$$
\begin{vmatrix}
\frac{\partial\beta_1(\epsilon_1, \epsilon_2)}{\partial\epsilon_1}&\frac{\partial\beta_1(\epsilon_1, \epsilon_2)}{\partial\epsilon_2}\\
\frac{\partial\beta_2(\epsilon_1, \epsilon_2)}{\partial\epsilon_1}&\frac{\partial\beta_2(\epsilon_1, \epsilon_2)}{\partial\epsilon_2}
\end{vmatrix}_{(\epsilon_1, \epsilon_2)=(0, 0)}
=\frac{B^6(0, 0)J_0}{A^5(0, 0)}
\neq0
$$
with
\begin{eqnarray*}
\begin{array}{l}
J_0:=\frac{-x_*^4 (h x_*^2-x_*^2+6) \{h^2 x_*^2+(-2 x_*^2+1) h+x_*^2-2\}^2}{4(h-1)^2 (h x_*^2-x_*^2+1) (h x_*^2-x_*^2+2) \{3 h^3 x_*^4-2 x_*^2 (4 x_*^2-3) h^2+(x_*^2-1) (7 x_*^2-3) h-2 x_*^2 (x_*^2-2)\} }
\end{array}
\end{eqnarray*}
under the conditions of Theorem \ref{thmBT},
the parameter transformation (\ref{(3-14)}) is a homeomorphism in a small neighborhood of the origin and parameters $\beta_1(\epsilon_1, \epsilon_2)$ and $\beta_2(\epsilon_1, \epsilon_2)$ are independent.
The results in section 8.4 of \cite{Kuznetsov} indicate that  system (\ref{(3-13)}) undergoes the
Bogdanov-Takens bifurcation of codimension two. More specifically, system (\ref{(3-13)}) undergoes
a saddle-node bifurcation as $(\beta_1(\epsilon_1, \epsilon_2),\beta_2(\epsilon_1, \epsilon_2))$ crossing $\mathcal{SN}^+\cup\mathcal{SN}^-$, where
\begin{eqnarray*}
\begin{array}{l}
\mathcal{SN}^+=\{(\beta_1(\epsilon_1, \epsilon_2),\beta_2(\epsilon_1, \epsilon_2))\in U:\beta_1(\epsilon_1, \epsilon_2)=\frac{1}{4}\beta_2^2(\epsilon_1, \epsilon_2), \beta_2(\epsilon_1, \epsilon_2)>0\},\\
\mathcal{SN}^-=\{(\beta_1(\epsilon_1, \epsilon_2),\beta_2(\epsilon_1, \epsilon_2))\in U:\beta_1(\epsilon_1, \epsilon_2)=\frac{1}{4}\beta_2^2(\epsilon_1, \epsilon_2), \beta_2(\epsilon_1, \epsilon_2)<0\},
\end{array}
\end{eqnarray*}
a Hopf bifurcation as $(\beta_1(\epsilon_1, \epsilon_2),\beta_2(\epsilon_1, \epsilon_2))$ crossing $\mathcal{H}$, where
\begin{eqnarray*}
\begin{array}{l}
\mathcal{H}=\{(\beta_1(\epsilon_1, \epsilon_2),\beta_2(\epsilon_1, \epsilon_2))\in U:\beta_1(\epsilon_1, \epsilon_2)=0, \beta_2(\epsilon_1, \epsilon_2)<0\},
\end{array}
\end{eqnarray*}
and a homoclinic bifurcation as $(\beta_1(\epsilon_1, \epsilon_2),\beta_2(\epsilon_1, \epsilon_2))$ crossing $\mathcal{HL}$, where
\begin{eqnarray*}
\begin{array}{l}
\mathcal{HL}=\{(\beta_1(\epsilon_1, \epsilon_2),\beta_2(\epsilon_1, \epsilon_2))\in U:\beta_1(\epsilon_1, \epsilon_2)=-\frac{6}{25}\beta_2^2(\epsilon_1, \epsilon_2)+O(|\beta_2(\epsilon_1, \epsilon_2)|^3), \\
\phantom{\mathcal{HL}:=}
\beta_2(\epsilon_1, \epsilon_2)<0\}.
\end{array}
\end{eqnarray*}

In the following, we express the four bifurcation curves $\mathcal{SN}$, $\mathcal{H}$ and $\mathcal{HL}$ in terms of $\epsilon_1$ and $\epsilon_2$.
From (\ref{(3-14)}), we solve $\epsilon_1$ and $\epsilon_2$ as follows
\begin{eqnarray*}
\begin{array}{l}
\epsilon_1=g_{10}\beta_1+g_{01}\beta_2+O(\parallel(\beta_1,\beta_2)\parallel^2),~~
\epsilon_2=h_{10}\beta_1+h_{01}\beta_2+O(\parallel(\beta_1,\beta_2)\parallel^2),
\end{array}
\end{eqnarray*}
where the coefficients $g_{ij}$ and $h_{ij}$ are given in the Appendix.
For the saddle-node bifurcation curves $\mathcal{SN}$, we consider $\Gamma:=\beta_1(\epsilon_1, \epsilon_2)-\frac{\beta_2^2(\epsilon_1, \epsilon_2)}{4}=0$.
Since
\begin{eqnarray*}
\begin{array}{l}
\frac{\partial\Gamma}{\partial\epsilon_1}=\frac{2 x_*^3 (h-1)^2 \{3 h^3 x_*^4-2 x_*^2 (4 x_*^2-3) h^2+(x_*^2-1) (7 x_*^2-3) h-2 x_*^2 (x_*^2-2)\}^4}{(h x_*^2-x_*^2+6)^3 (h x_*^2-x_*^2+1)^2 \{h^2 x_*^2+(-2 x_*^2+1) h+x_*^2-2\}^3}\neq0
\end{array}
\end{eqnarray*}
at the origin under the conditions of Theorem \ref{thmBT},
the implicit function theorem indicates that there exists a unique function $\epsilon_1(\epsilon_2)$ such that $\epsilon_1(0)=0$ and $\Gamma(\epsilon_1(\epsilon_2),\epsilon_2)=0$, which can be obtained as an expansion
\begin{eqnarray*}
\begin{array}{l}
\epsilon_1(\epsilon_2)=f_{11}\epsilon_2+f_{12}\epsilon_2^2+O(|\epsilon_2|^3),
\end{array}
\end{eqnarray*}
where the coefficients $f_{11}$ and $f_{12}$ are given in the Appendix.
Furthermore, we get
$\epsilon_2=h_{01}\beta_2+O(|\beta_2|^2)$
restricted on the curve $\Gamma=0$ and the coefficient of $\beta_2$ is negative, which implies that
 $\epsilon_2$ and $\beta_2$ have the opposite sign. We therefore obtain the bifurcation curves $\mathcal{SN}$ as
\begin{eqnarray*}
\begin{array}{l}
\mathcal{SN}^+=\{(\epsilon_1, \epsilon_2)\in U:\epsilon_1=f_{11}\epsilon_2+f_{12}\epsilon_2^2+O(|\epsilon_2|^3), \epsilon_2<0\},\\
\mathcal{SN}^-=\{(\epsilon_1, \epsilon_2)\in U:\epsilon_1=f_{11}\epsilon_2+f_{12}\epsilon_2^2+O(|\epsilon_2|^3), \epsilon_2>0\}.
\end{array}
\end{eqnarray*}
For the Hopf bifurcation curve $\mathcal{H}$, we need to consider the curve $\beta_1(\epsilon_1, \epsilon_2)=0$. Since $\frac{\partial\beta_1}{\partial\epsilon_1}\neq0$ at the origin, we similarly obtain that there exists a unique function
\begin{eqnarray*}
\begin{array}{l}
\epsilon_1(\epsilon_2)=f_{21}\epsilon_2+f_{22}\epsilon_2^2+O(|\epsilon_2|^3),
\end{array}
\end{eqnarray*}
where $f_{21}:=f_{11}$ and the coefficient  $f_{22}$ is given in the Appendix.
Furthermore, we get $\epsilon_2=h_{01}\beta_2+O(|\beta_2|^2)$
restricted on the curve $\beta_1(\epsilon_1, \epsilon_2)=0$.
It follows that $\epsilon_2>0$ if $\beta_2<0$. Thus, we get the bifurcation curve $\mathcal{H}$ as
\begin{eqnarray*}
\begin{array}{l}
\mathcal{H}=\{(\epsilon_1, \epsilon_2)\in U:\epsilon_1=f_{21}\epsilon_2+f_{22}\epsilon_2^2+O(|\epsilon_2|^3), \epsilon_2>0\}.
\end{array}
\end{eqnarray*}
For the  homoclinic bifurcation curve $\mathcal{HL}$, we consider the curve $\Xi:=\beta_1+\frac{6}{25}\beta_2^2+O(|\beta_2|^3)$.
Since $\frac{\partial\Xi}{\partial\epsilon_1}\neq0$ at the origin,
we also obtain that there exists a unique function
\begin{eqnarray*}
\begin{array}{l}
\epsilon_1(\epsilon_2)=f_{31}\epsilon_2+f_{32}\epsilon_2^2+O(|\epsilon_2|^3),
\end{array}
\end{eqnarray*}
where  $f_{31}:=f_{11}$ and the coefficient $f_{32}$ is given in the Appendix.
Furthermore, we have
$\epsilon_2=h_{01}\beta_2+O(|\beta_2|^2)$
restricted on the curve $\Xi=0$ and the coefficient of $\beta_2$ is negative, which implies that
 $\epsilon_2>0$ if $\beta_2<0$ . Therefore, we obtain the  bifurcation curve $\mathcal{HL}$ as
\begin{eqnarray*}
\begin{array}{l}
\mathcal{HL}=\{(\epsilon_1, \epsilon_2)\in U:\epsilon_1=f_{31}\epsilon_2+f_{32}\epsilon_2^2+O(|\epsilon_2|^3), \epsilon_2>0\}.
\end{array}
\end{eqnarray*}
With the linear transformation $\epsilon_1=\alpha-\alpha_*$ and $\epsilon_2=\sigma-\sigma_*$, we can rewrite the above four bifurcation curves as in Theorem \ref{thmBT}.
The proof is completed.
\qquad$\Box$


\subsection{Hopf Bifurcation}
Theorem \ref{thm1} shows that $E_1$ is a center or weak focus as $T(x_1)=0$ under the conditions $0<h<1$ and $(\alpha,\kappa,\sigma)\in\mathcal{P}_7$.
In this subsection, we are devoted to the center-focus determination and obtain that the multiplicity of weak focus $E_1$ is at most three and the Hopf bifurcation occurs.

In order to have a Hopf bifurcation at $E_1$, we look for some parameter conditions such that $E_1$ is a nonhyperbolic equilibrium satisfying
$F(x_1)=0$, $T(x_1)=0$ and $D|_{E_{1}}>0$, where
\begin{eqnarray*}
\begin{array}{l}
D|_{E_{1}}=\frac{\{(h-1)x_1^2+3\}\kappa-2x_1\{(h-1)x_1^2+2\}}{(\kappa-x_1)(1-h)x_1^2}.
\end{array}
\end{eqnarray*}
Thus, we can express $\alpha$ and $\sigma$ by $\kappa$, $h$ and $x_1$ as follows
\begin{eqnarray}
\begin{array}{l}
\alpha=\frac{(2h-1)\kappa-2hx_1}{x_1^3(h-1)^2(x_1-\kappa)},~~~
\sigma=\frac{\kappa(h-1)\{(h-1)x_1^2+1\}}{(2h-1)\kappa-2hx_1}
\end{array}
\label{(4-1)}
\end{eqnarray}
and $\kappa$, $h$ and $x_1$ satisfy $(\kappa,h,x_1)\in\mathcal{P}$ with
\begin{eqnarray}
\begin{array}{l}
\mathcal{P}:=\{(\kappa,h,x_1)\in\mathbb{R}_+^3:
\kappa>\kappa_2, h\leq\frac{1}{2}, x_1\leq1 ~~\mbox{or}~~
\kappa>\kappa_2, \frac{x_1^2-1}{x_1^2}<h\leq\frac{1}{2},\\
\phantom{\mathcal{P}:=\{}
 1<x_1<\sqrt{2} ~~\mbox{or}~~\kappa_2<\kappa<\frac{2hx_1}{2h-1}, \frac{1}{2}<h<1, x_1\leq\sqrt{2} ~~\mbox{or}~~\\
\phantom{\mathcal{P}:=\{}
\kappa_2<\kappa<\frac{2hx_1}{2h-1}, \frac{x_1^2-1}{x_1^2}<h<1, x_1>\sqrt{2}
\},
\end{array}
\label{(4-4)}
\end{eqnarray}
where $\kappa_2$ is given in subsection 3.2.
Then under the critical conditions (\ref{(4-1)}) and (\ref{(4-4)}),  the linearized system of (\ref{(3)}) at equilibrium $E_1$
has a pair of purely imaginary  eigenvalues $\pm i\omega$ with $\omega:=\sqrt{D|_{E_{1}}}$.
Hence we consider the Hopf bifurcation at $E_1$ under conditions (\ref{(4-1)}) and (\ref{(4-4)})
and obtain the following result.
\begin{thm}
Under conditions (\ref{(4-1)}) and $(\kappa,h,x_1)\in\mathcal{P}$,
equilibrium $E_1$ is a weak focus of multiplicity at most three.
More concretely,
$E_1$ is a weak focus of multiplicity one, two and three as
$(\kappa,h,x_1)\in\mathcal{P}_1$, $(\kappa,h,x_1)\in\mathcal{P}_2$ and $(\kappa,h,x_1)\in\mathcal{P}_3$ respectively, where  $\mathcal{P}_1:=\mathcal{P}\backslash(\mathcal{P}_2\cup\mathcal{P}_3)$,
$\mathcal{P}_2:=\{(\kappa,h,x_1)\in\mathcal{P}: f_1=0, f_2\neq0\}$,
$\mathcal{P}_3:=\{(\kappa,h,x_1)\in\mathcal{P}: f_1=f_2=0\}$ with
$f_1$ and $f_2$ given in the Appendix.
Thus, at most three limit cycles arise from the Hopf bifurcation at equilibrium $E_1$.
\label{thm4}
\end{thm}
\textbf{Proof}.
Translating $E_1$ to $(0,0)$ and using the linear transformation
\begin{eqnarray*}
\begin{array}{l}
x=\frac{(2 h-1) \kappa-2 h x_1}{2\kappa (\kappa-x_1) \{(h-1) x_1^2+1\} (h-1)}u-\frac{ (2 h-1) \kappa-2 h x_1}{2\omega (\kappa-x_1) (h-1)}v  ~~\mbox{and}~~  y=\frac{1}{\omega}v
\end{array}
\end{eqnarray*}
and the time-rescaling $\tau:=\omega t$ with $\omega$ given above to normalize the linear part of system (\ref{(3)}),
we can change system (\ref{(3)}) into the following system
\begin{eqnarray}
\left\{
\begin{array}{l}
\frac{du}{dt}=-v+a_{20}u^2+a_{02}v^2+a_{11}uv+a_{30}u^3+a_{21}u^2v+a_{12}uv^2+a_{03}v^3\\
\phantom{\frac{du}{dt}=}
+a_{40}u^4+a_{31}u^3v+a_{22}u^2v^2+a_{13}uv^3+a_{04}v^4+a_{41}u^4v+a_{32}u^3v^2\\
\phantom{\frac{du}{dt}=}
+a_{23}u^2v^3+a_{14}uv^4+a_{05}v^5,\\
\frac{dv}{dt}=
u+b_{20}u^2+b_{11}uv+b_{02}v^2+b_{21}u^2v+b_{12}uv^2+b_{03}v^3+b_{22}u^2v^2\\
\phantom{\frac{du}{dt}=}
+b_{13}uv^3+b_{04}v^4,\\
\end{array}
\right.
\label{(4-5)}
\end{eqnarray}
where these coefficients are given in the Appendix.
In the following, we are devoted to identifying a weak focus from a center and determining the multiplicity of the weak focus by using the successive function method(\cite{Zhang}).
We obtain the first three focal values $l_i$$(i=1,2,3)$, where
\begin{eqnarray}
\begin{array}{l}
l_1:= \frac{\kappa f_1}{16(h-1)^2 (\kappa-x_1)^2 \omega^3 \{(h x_1^2-x_1^2+3) \kappa-2 x_1 (h x_1^2-x_1^2+2)\} x_1},\\
l_2:=\frac{\kappa^5 f_2}{3072(h-1)^4 (\kappa-x_1)^4 \omega^9 \{(h x_1^2-x_1^2+3) \kappa-2 x_1 (h x_1^2-x_1^2+2)\} (2 h \kappa-2 h x_1-\kappa)^2 x_1^3},\\
l_3:=\frac{\kappa^9f_3}{70778880(h-1)^6 (\kappa-x_1)^6 \omega^{15} \{(h x_1^2-x_1^2+3) \kappa-2 x_1 (h x_1^2-x_1^2+2)\} (2 h \kappa-2 h x_1-\kappa)^4 x_1^5},\\
\end{array}
\label{(4-6)}
\end{eqnarray}
where $f_1$ and $f_2$ are given in the Appendix and the polynomial $f_3$ of 1755 terms is too long to display.
Because the denominators of the three focal values are positive, the real zeros and signs of $l_i$$(i = 1, 2, 3)$ are
determined by the factors $f_i$$(i = 1, 2, 3)$ respectively.
To determine the multiplicity of weak focus $E_1$, we need to compute the two algebraic varieties $V(f_1, f_2)$ and $V(f_1, f_2, f_3)$ in $\mathcal{P}$, which are the sets of common zeros of two polynomials $f_i$$(i = 1, 2)$ and three polynomials $f_i$$(i = 1, 2, 3)$ respectively(\cite{Chen,Winkler}). The tedious symbolic computation of $V(f_1, f_2)$ and $V(f_1, f_2, f_3)$ is displayed in the Appendix.
It is shown that $V(f_1, f_2, f_3)=\emptyset$ but $V(f_1, f_2)\neq\emptyset$ in $\mathcal{P}$, which implies that
the multiplicity of weak focus $E_1$ is at most three and at most three limit cycles arise from the Hopf bifurcation at $E_1$.
The proof is completed.
\qquad$\Box$

Theorem \ref{thm4} shows that the degenerate Hopf bifurcation at the weak focus of multiplicity three produces at most three hyperbolic limit cycles.
Conversely, these limit cycles may coalesce into a multiple limit cycle through the multiple limit cycle bifurcation,
or destroy through either the Hopf bifurcation or the homoclinic bifurcation or even the multiple limit cycle bifurcation as perturbing the parameters.

\begin{rem}
System (\ref{(3)}) has more complicated bifurcation phenomena
than the system without cooperative hunting among predators (i.e., $\alpha=0$ in system (\ref{(3)})),
because the system without cooperative hunting only undergoes the transcritical bifurcation at the saddle-node
$E_{\kappa}$ as $\kappa$ varies across the bifurcation value $\kappa=\kappa_1$ and the Hopf bifurcation around
the weak focus of multiplicity one $(\kappa_1,\sigma\kappa_1(1-\frac{\kappa_1}{\kappa}))$ if  $\kappa=\frac{2h\kappa_1}{2h-1}$ and $\frac{1}{2}<h<1$ and at most one limit cycle bifurcates from  the positive equilibrium(\cite{Bazykin,ChenZhang1986}).
\label{rem2}
\end{rem}

\section{Simulations and Discussions}

The cooperative hunting among predators is a ubiquitous behavior in ecological system and
plays an important role in determining the dynamics of predator-prey system.
In this paper, we consider the dynamics of the predator-prey system with cooperative hunting among predators and Holling III functional response.
We first give the parametric conditions for the existence and the qualitative properties of the equilibria(Theorem \ref{thm1}), and then discuss the various
bifurcations of the system around the nonhyperbolic  equilibria such as the transcritical and the pitchfork bifurcations at the degenerate boundary equilibrium(Theorem \ref{thm2}), the saddle-node and the Bogdanov-Takens bifurcations at the degenerate positive equilibrium(Theorems \ref{thm3} and \ref{thmBT}) and the Hopf bifurcation around the weak focus(Theorem \ref{thm4}).
Remarks \ref{rem1} and \ref{rem2} indicate that system (\ref{(3)})
has richer and more complicated  dynamics than the system without cooperative hunting,
because the later at most has three equilibria and only undergoes the transcritical and the Hopf  bifurcations(\cite{Bazykin,ChenZhang1986}).
Our theoretical analysis shows that system (\ref{(3)}) can display more bifurcation phenomena than the results of
Alves and Hilker(\cite{Alves}) showed through numerical simulations for some special parameter values.
The result in Theorem \ref{thm1} indicates that the predators will face extinction if their handing time to prey $h$ is too long ( i.e., $h\geq 1$), which is an understandable phenomenon,
because the predators spend too much time handing their prey so that they have much less effort to encounter prey.
Therefore, the predators may survive only if the handing time $h$ can not be too long (i.e., $h<1$).

\begin{figure}[H]
\centering
\subcaptionbox{\label{fig:subfig:a}}
{\includegraphics[width=7.5cm,height=5cm]{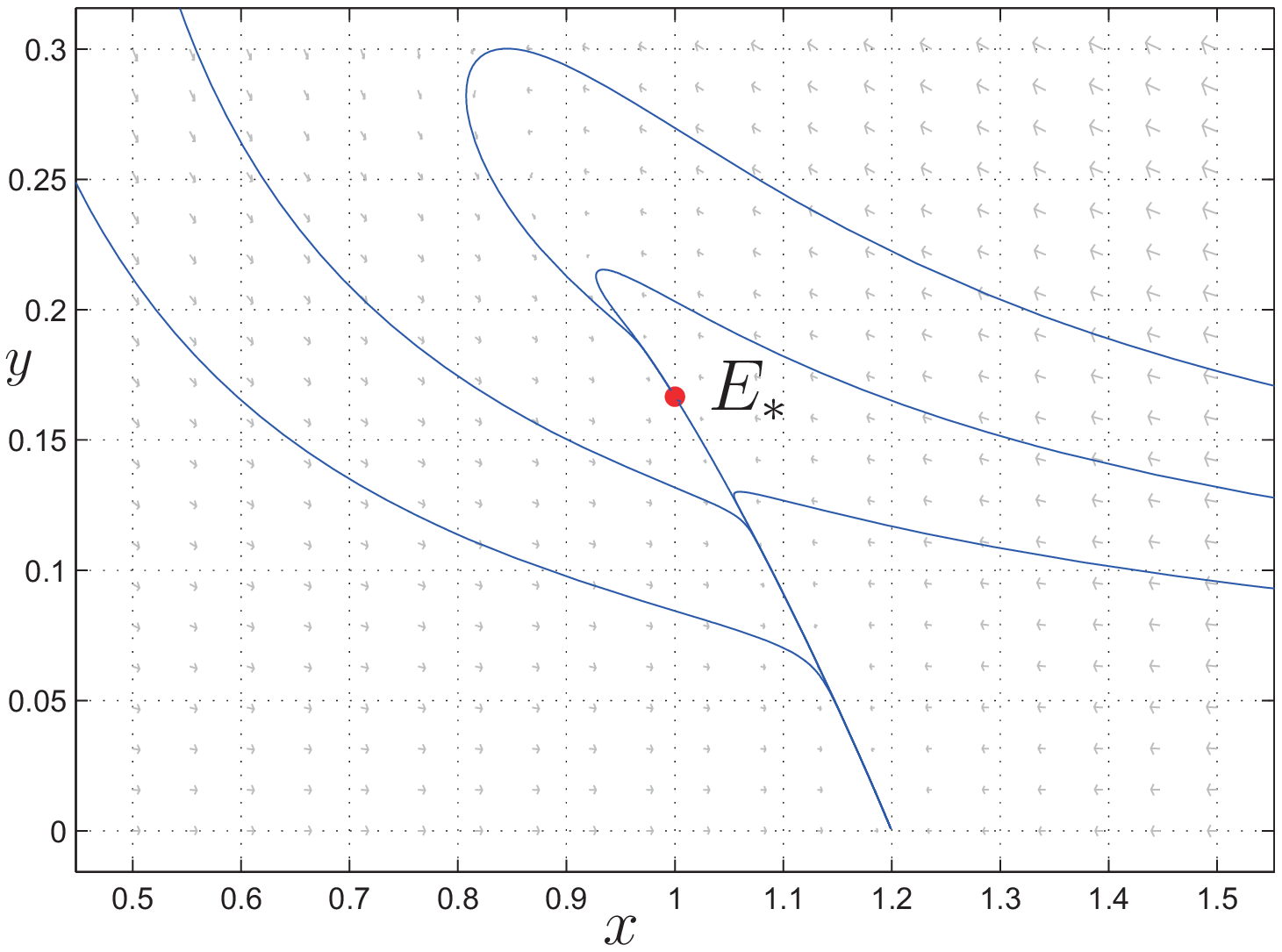}}\ \hspace{0.1cm}
\subcaptionbox{\label{fig:subfig:b}}
{\includegraphics[width=7.5cm,height=5cm]{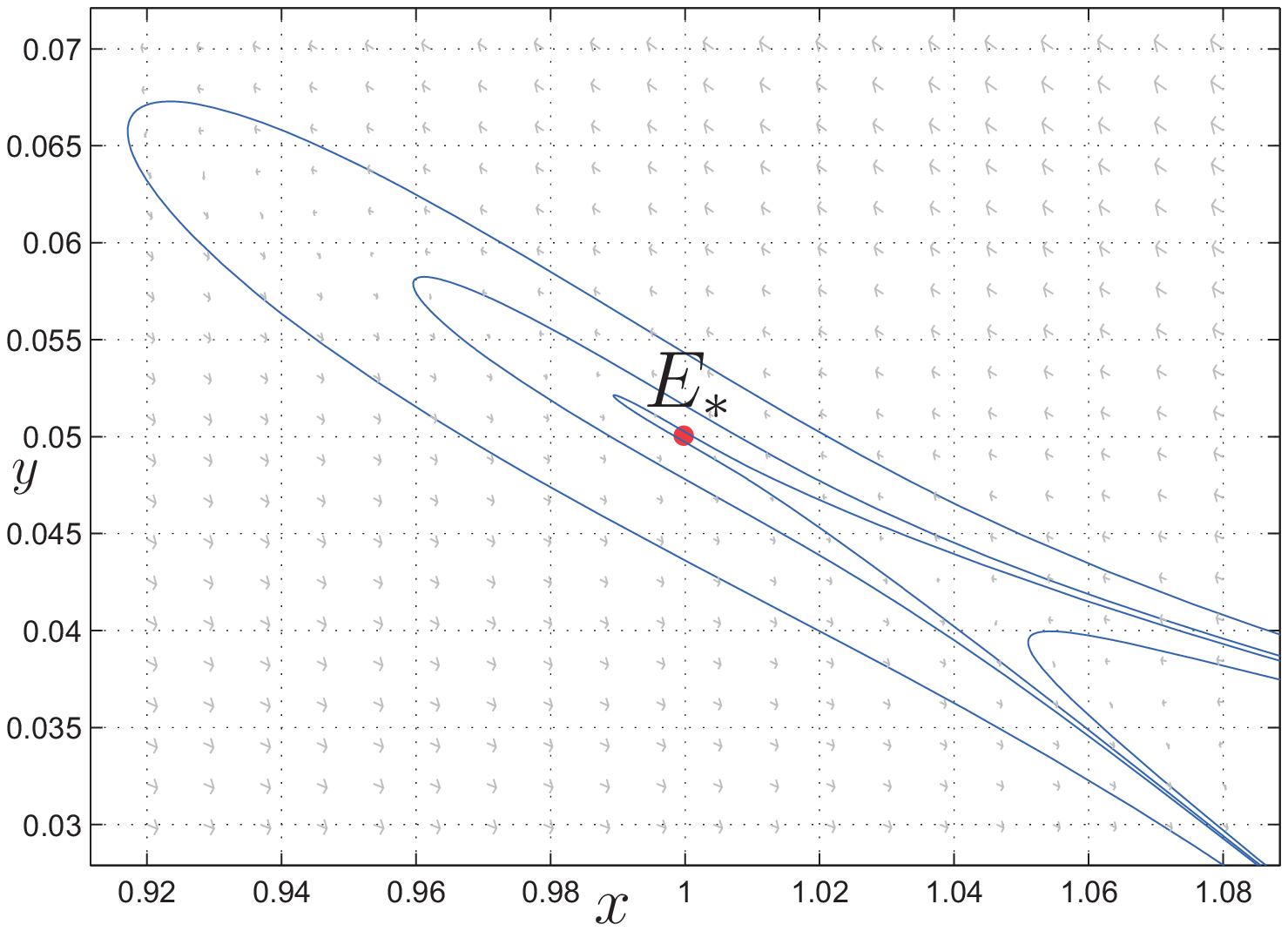}}\
\caption {The phase portraits of the degenerate equilibrium $E_*$ when $h=0.5$ and $\kappa=1.2$.
(a) Saddle-node $E_*$ when $\alpha=6$ and $\sigma=1$.
(b) Degenerate equilibrium (cusp) $E_*$ when $\alpha=20$ and $\sigma=0.3$.
}\label{Figure DegeEqui}
\end{figure}

\begin{figure}
\centering
\subcaptionbox{\label{fig:subfig:a}}
{\includegraphics[width=7.5cm,height=5cm]{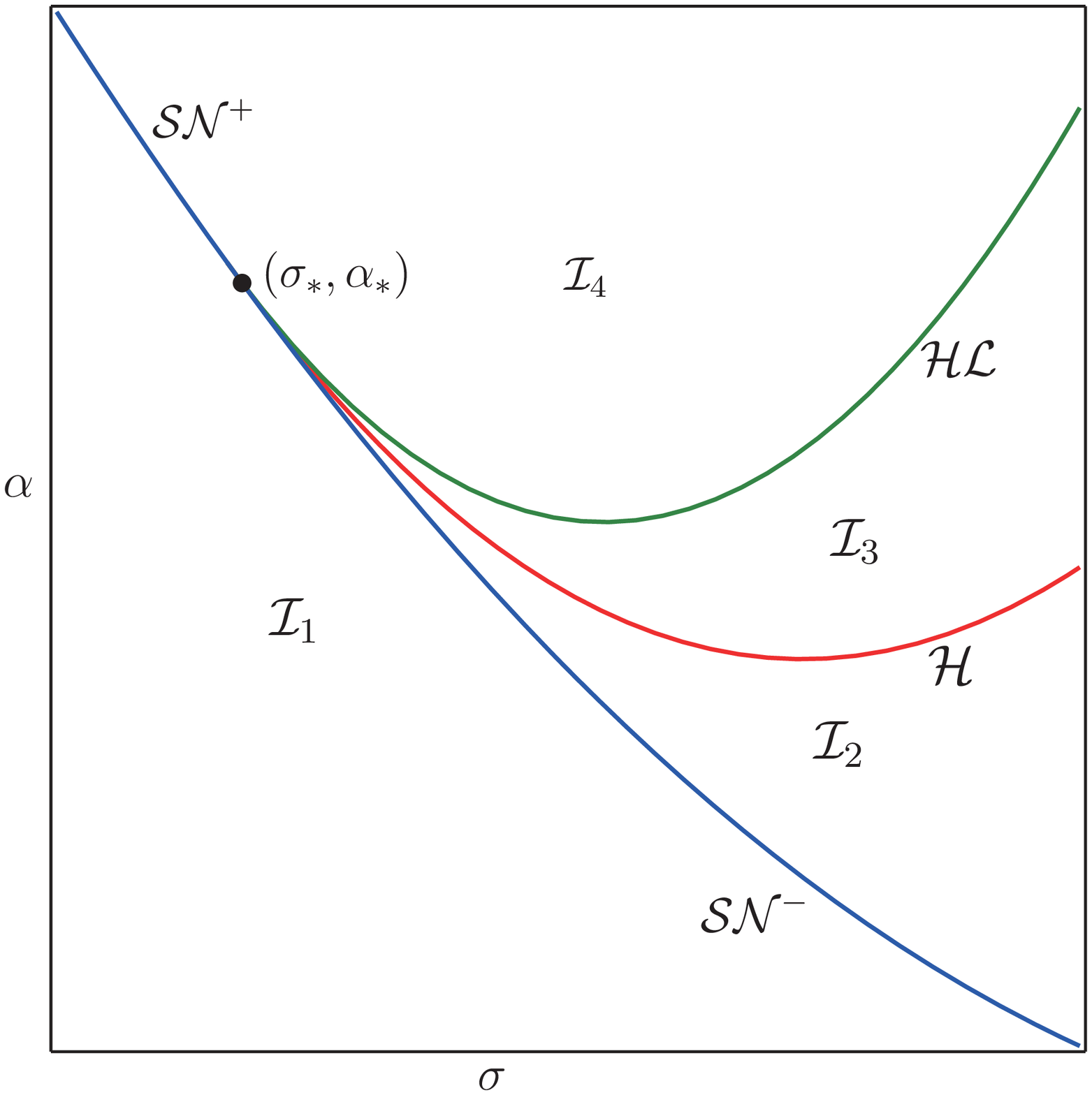}}\ \hspace{0.1cm}
\subcaptionbox{\label{fig:subfig:b}}
{\includegraphics[width=7.5cm,height=5cm]{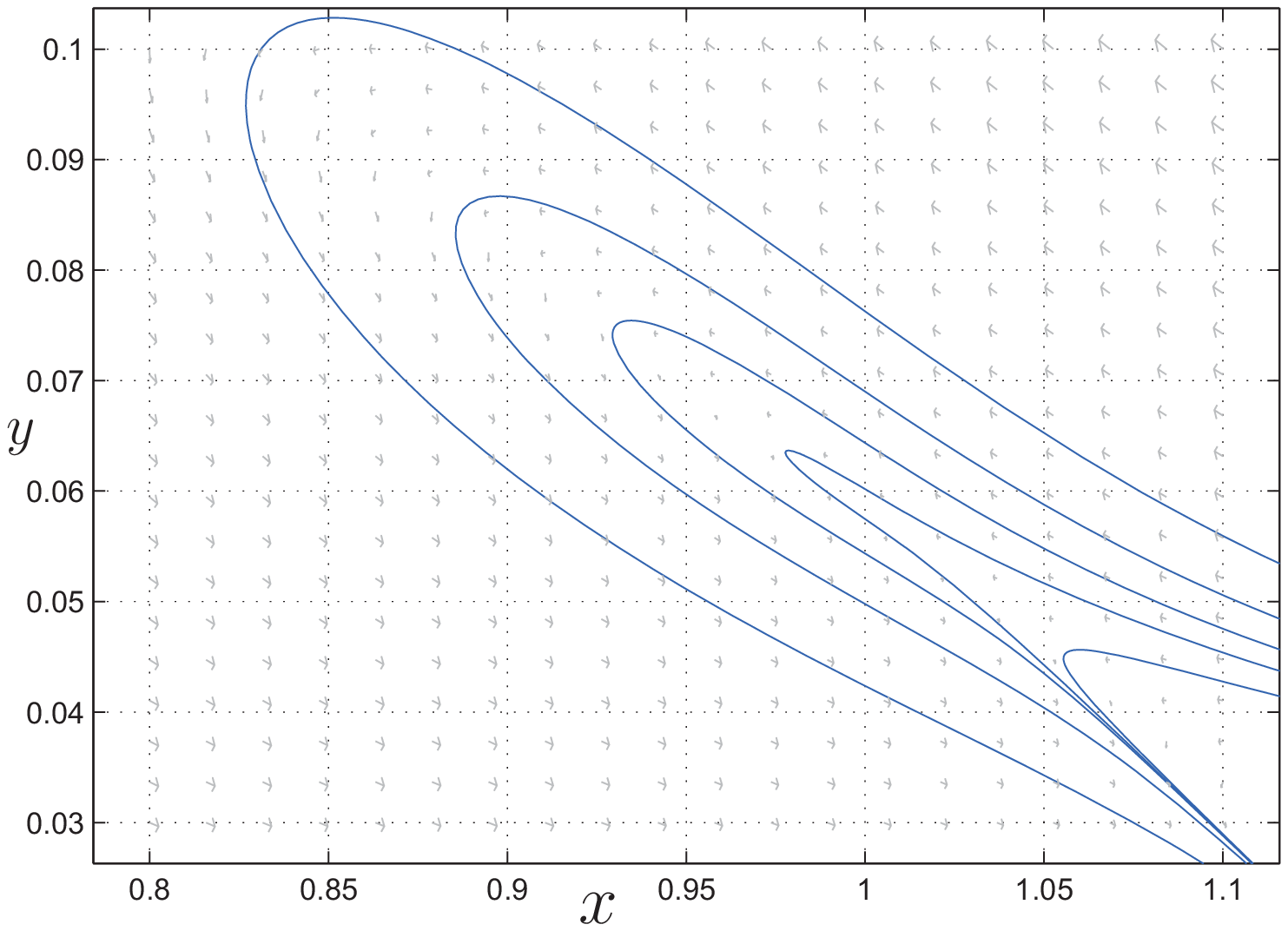}}\ \\
\subcaptionbox{\label{fig:subfig:c}}
{\includegraphics[width=7.5cm,height=5cm]{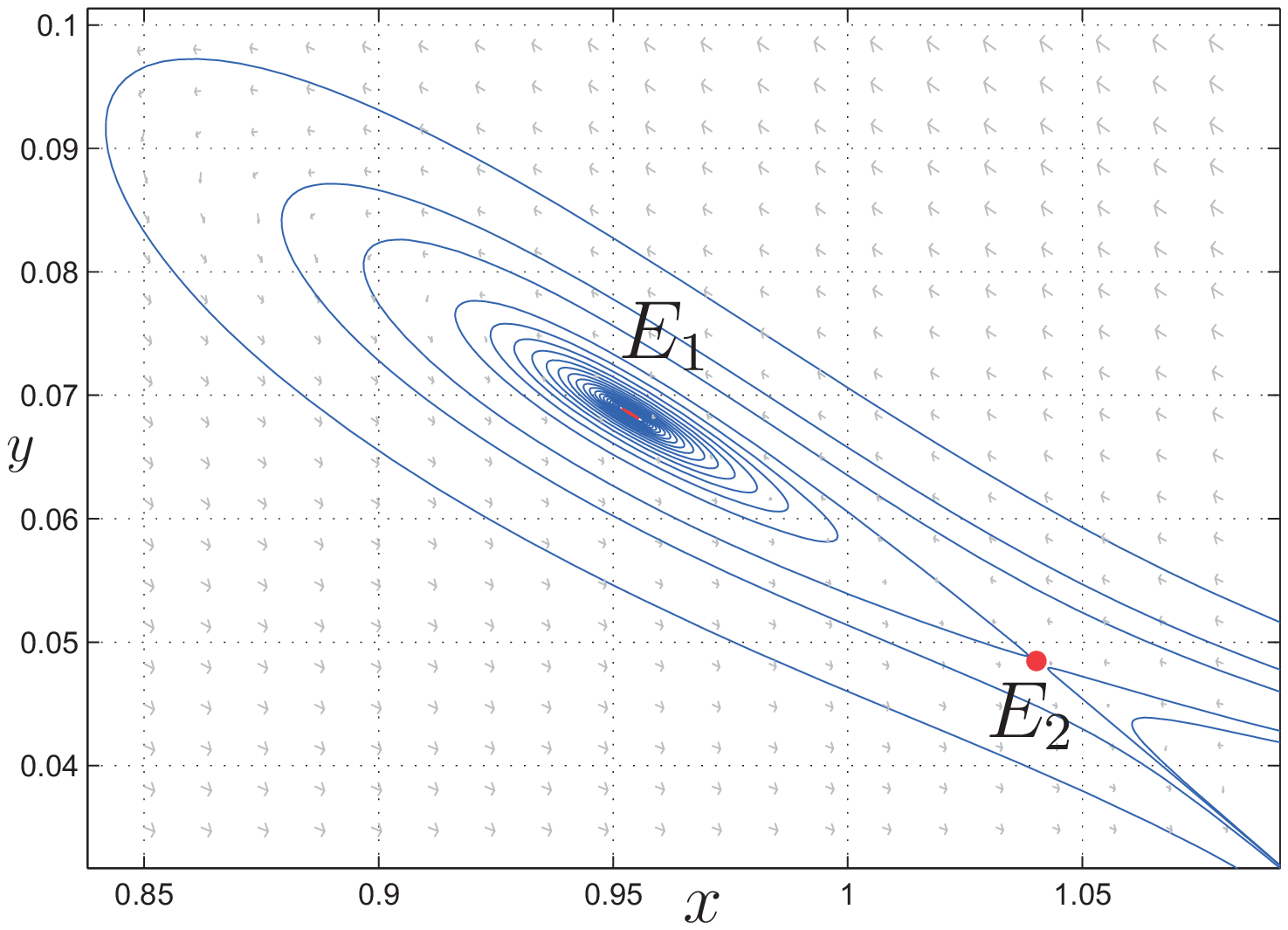}}\ \hspace{0.1cm}
\subcaptionbox{\label{fig:subfig:d}}
{\includegraphics[width=7.5cm,height=5cm]{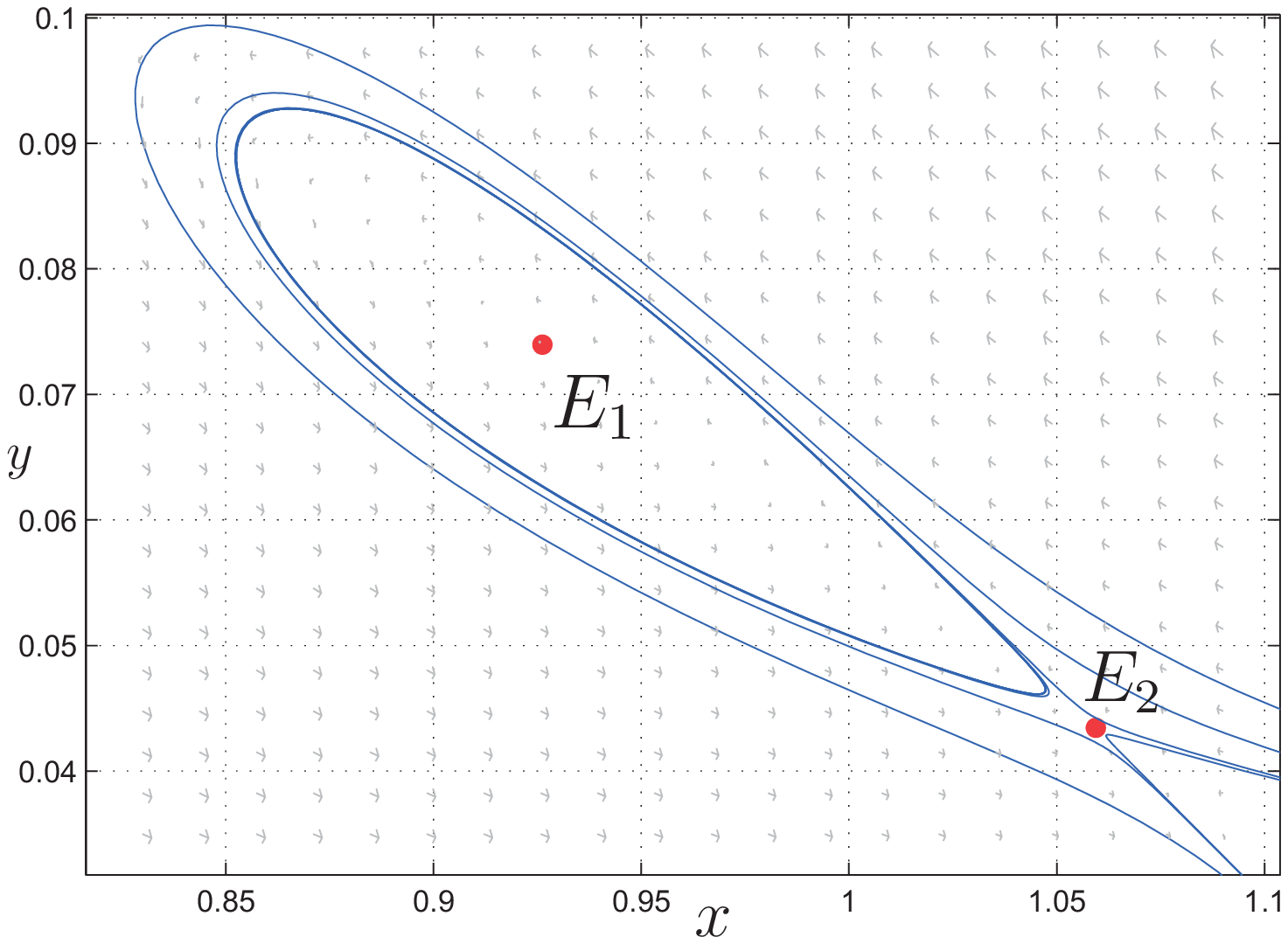}}\ \\
\subcaptionbox{\label{fig:subfig:e}}
{\includegraphics[width=7.5cm,height=5cm]{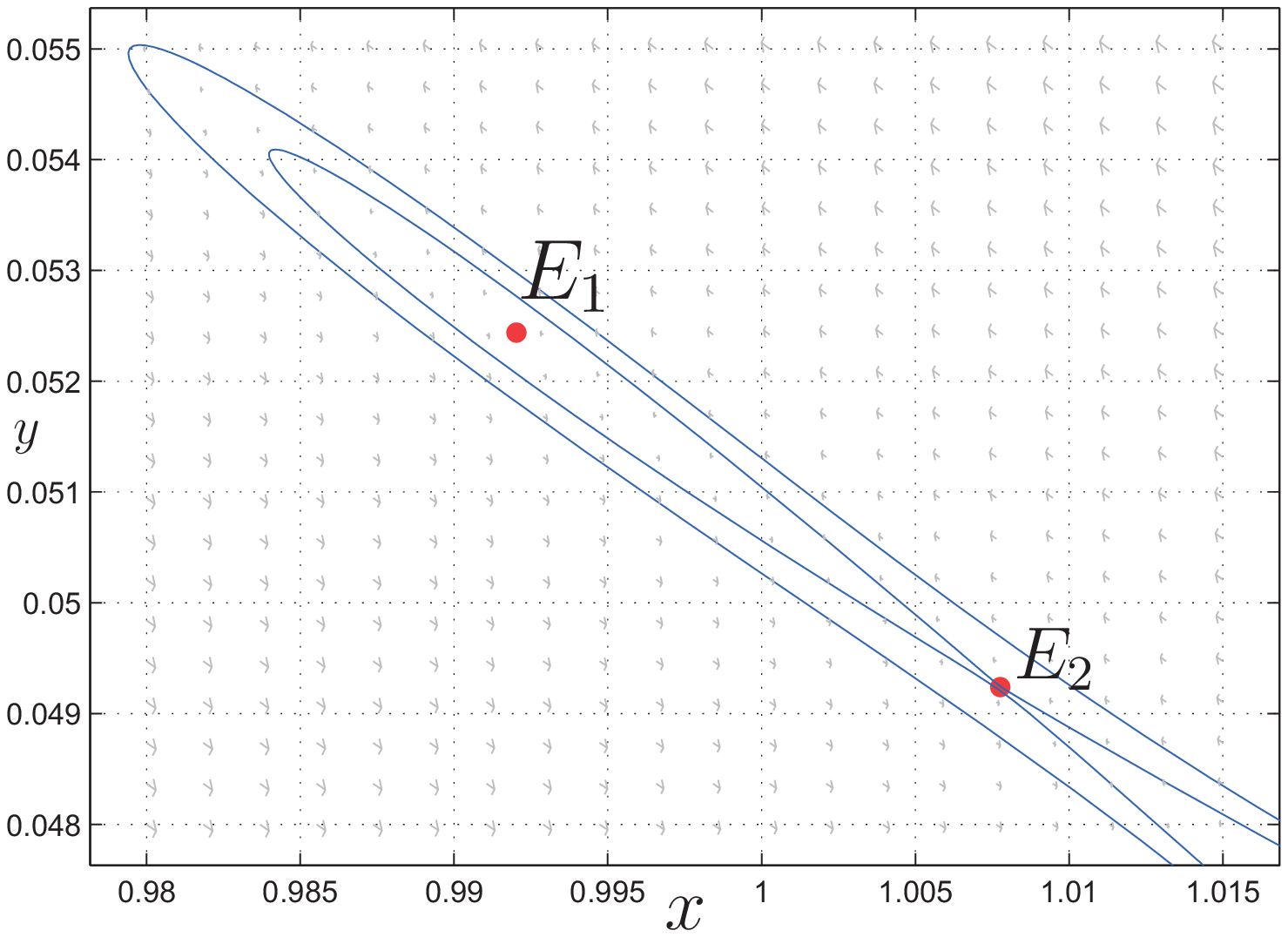}}\ \hspace{0.1cm}
\subcaptionbox{\label{fig:subfig:f}}
{\includegraphics[width=7.5cm,height=5cm]{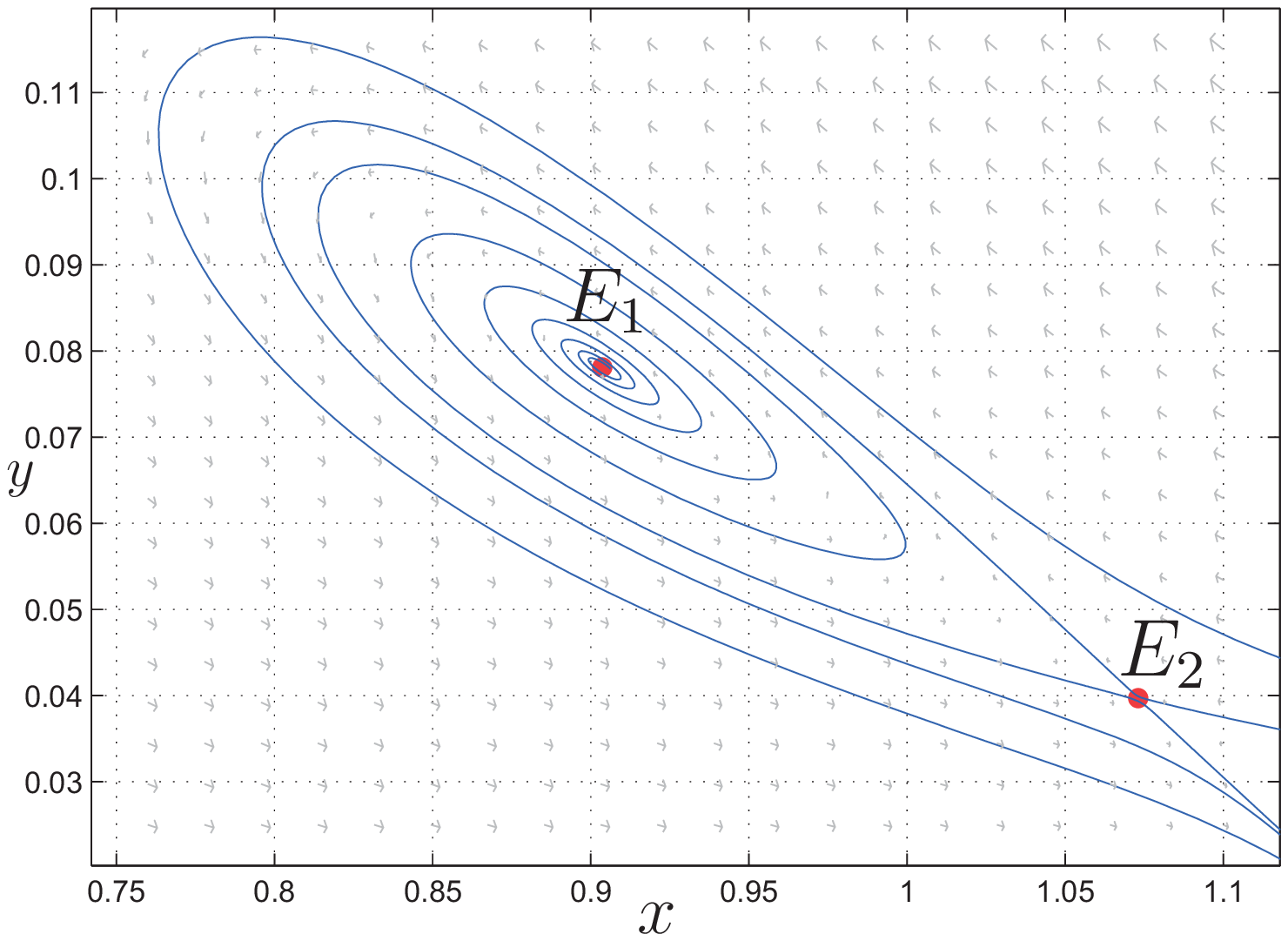}}\
\caption {The bifurcation diagram and corresponding phase portraits of system (\ref{(3)}) when $h=0.5$ and $\kappa=1.2$.
(a) Bifurcation diagram with $(\sigma_*,\alpha_*)=(0.3,20)$.
(b) No positive equilibrium when $(\sigma,\alpha)=(0.35,17.1)\in\mathcal{I}_{1}$.
(c) Stable focus $E_1$ and saddle $E_2$ when $(\sigma,\alpha)=(0.35,17.5)\in\mathcal{I}_{2}$.
(d) A stable limit cycle around unstable focus $E_1$ and saddle $E_2$ when $(\sigma,\alpha)=(0.35,18)\in\mathcal{I}_{3}$.
(e) A homoclinic orbit around unstable focus $E_1$ and saddle $E_2$ when $(\sigma,\alpha)=(0.305,19.6855)\in\mathcal{HL}$.
(f) Unstable focus $E_1$ and saddle $E_2$ when $(\sigma,\alpha)=(0.35,18.6)\in\mathcal{I}_{4}$.
}\label{Figure BT}
\end{figure}

\begin{figure}
\centering
\subcaptionbox{\label{fig:subfig:a}}
{\includegraphics[width=7.5cm,height=5.5cm]{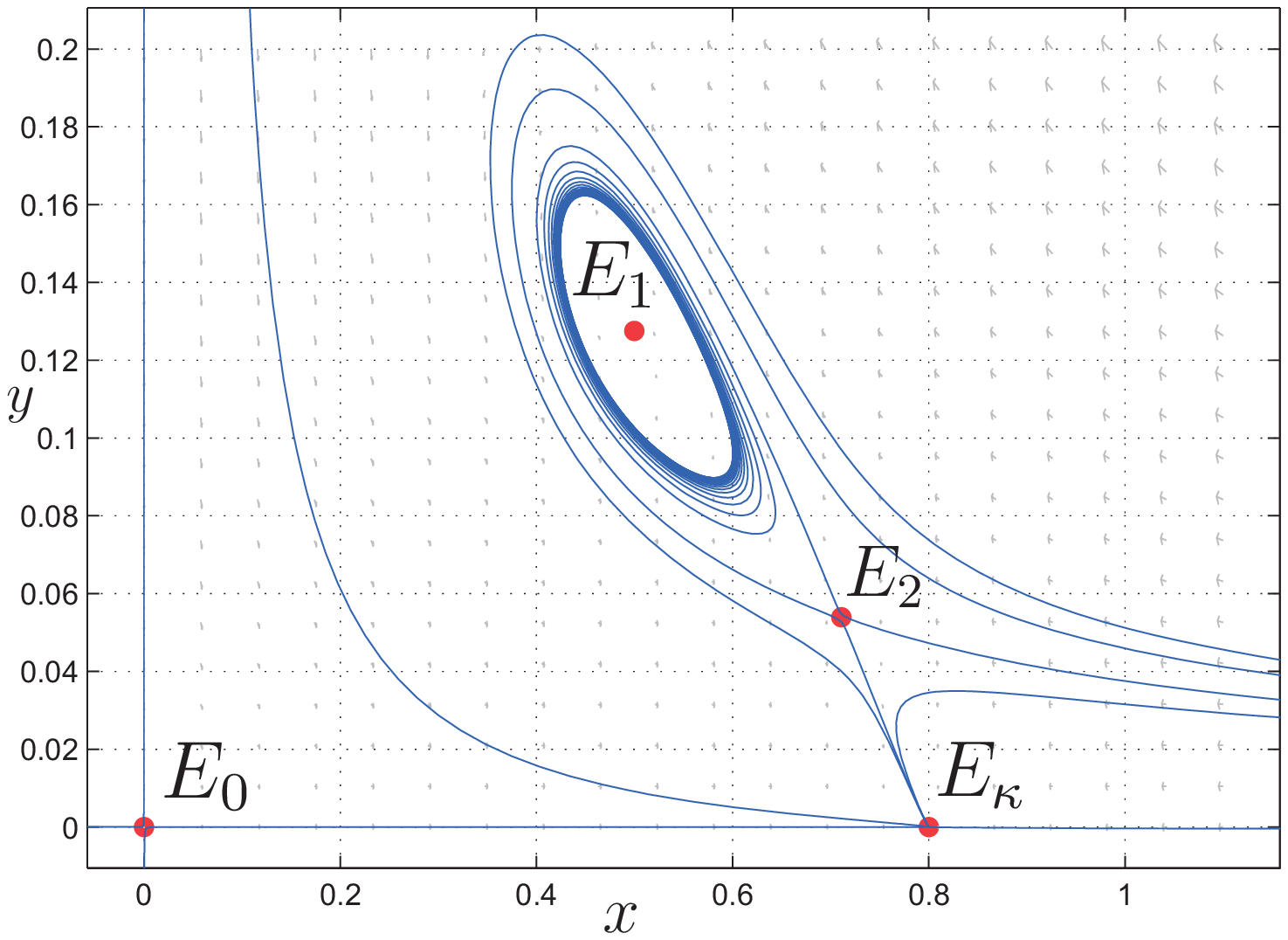}}\ \hspace{0.1cm}
\subcaptionbox{\label{fig:subfig:b}}
{\includegraphics[width=7.5cm,height=5.5cm]{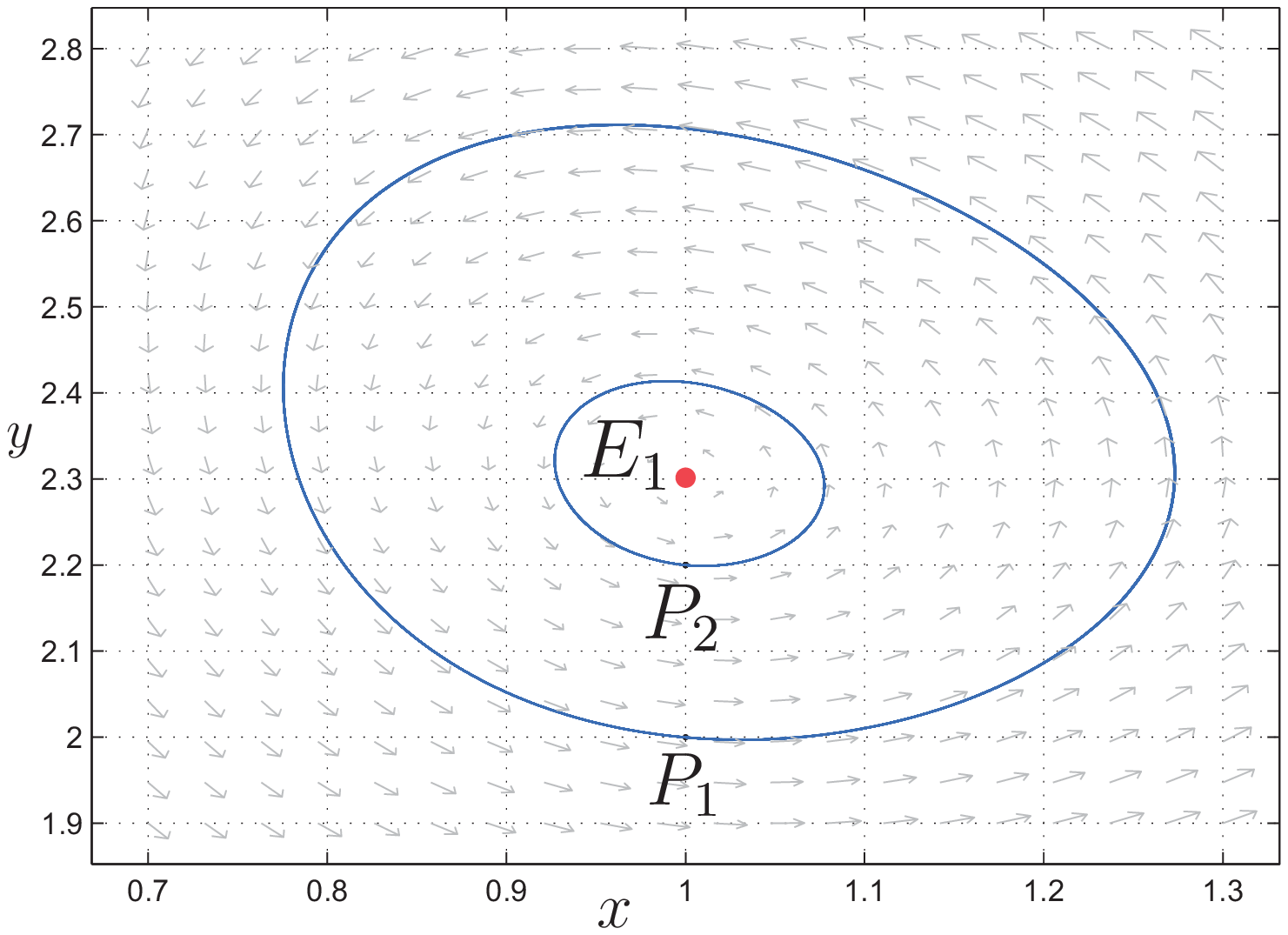}}\ \\
\subcaptionbox{\label{fig:subfig:c}}
{\includegraphics[width=7.5cm,height=5.5cm]{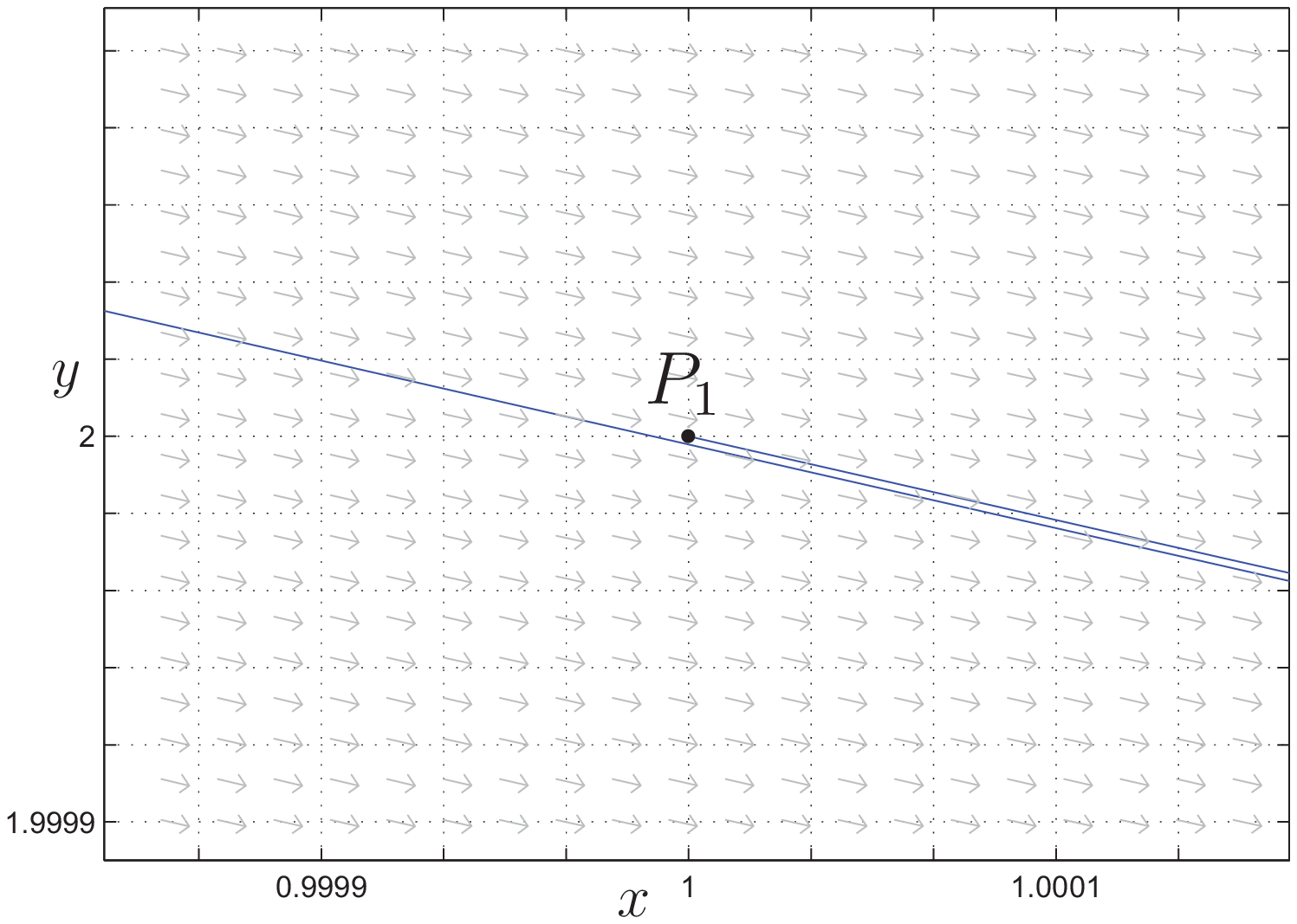}}\ \hspace{0.1cm}
\subcaptionbox{\label{fig:subfig:d}}
{\includegraphics[width=7.5cm,height=5.5cm]{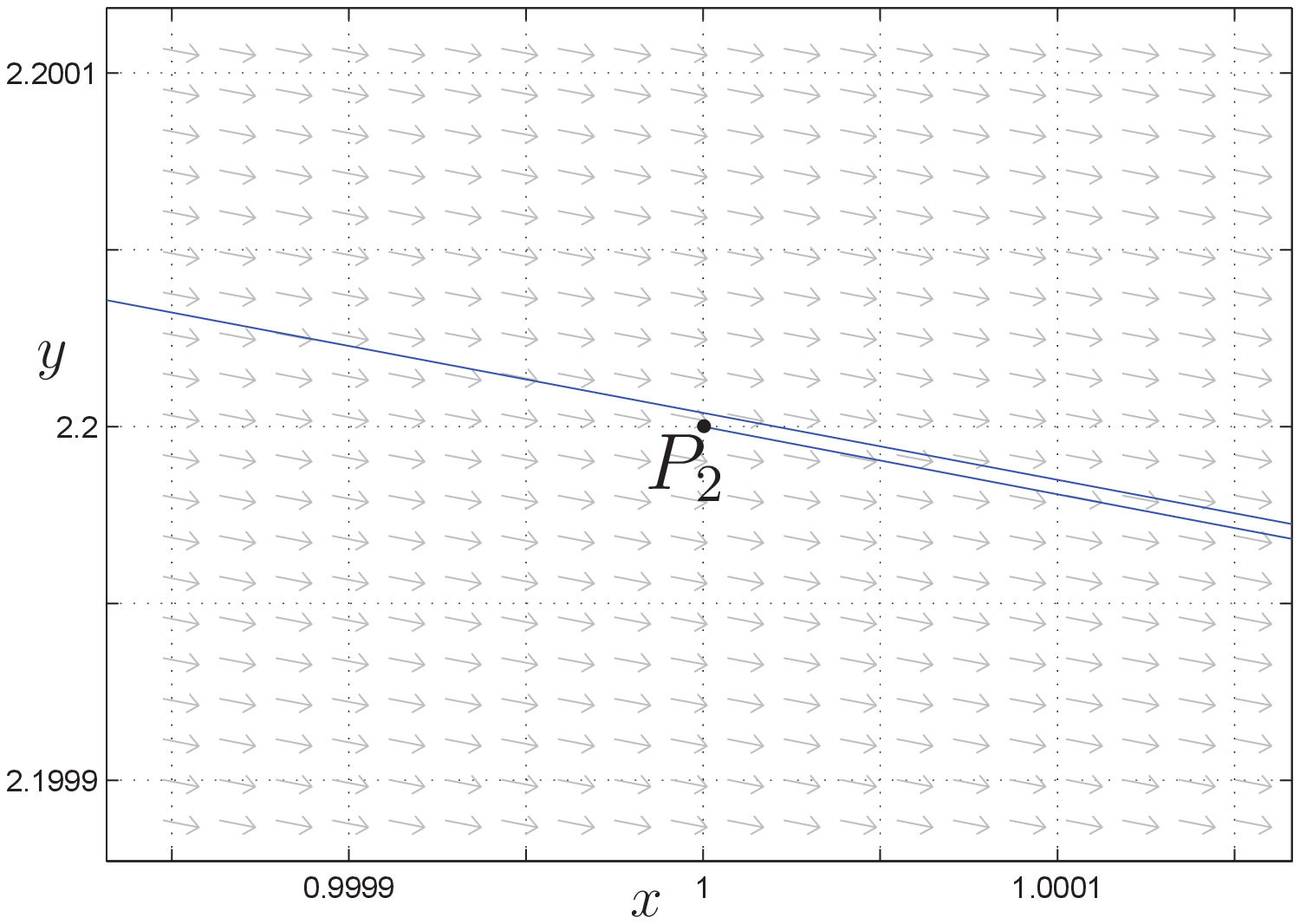}}\
\caption {Limit cycles induced by Hopf bifurcation at $E_1$ in system (\ref{(3)}).
(a) A stable limit cycle with $h=0.5$, $\kappa=0.8$, $\alpha=54.902$ and $\sigma=0.68$.
(b) Two limit cycles with $h=0.45$, $\kappa=133.7629$, $\alpha=0.3555$ and $\sigma=2.319$.
(c) The orbit from $P_1(1,2)$  spirals outward.
(d) The orbit from $P_2(1,2.2)$  spirals  inward.
}\label{Figure Hopf}
\end{figure}

In the section, we provide some numerical simulations to demonstrate our theoretical results and give some biological interpretations for the case $h<1$.
If taking $h=0.5$ and $\kappa=1.2$, then $E_*$ is a saddle-node when $\alpha=6$ and $\sigma=1$ (see Figure \ref{Figure DegeEqui}(a)), $E_*$ is a degenerate equilibrium of codimension two (i.e., a cusp) when  $\sigma=\sigma_*=0.3$ and $\alpha=\alpha_*=20$ (see Figure \ref{Figure DegeEqui}(b)). Moreover, the predator extinction equilibrium $E_{\kappa}$ is a stable node in both cases.
Thus, the numerical simulations indicate that
if both the environmental capacity of prey and the intensity of cooperative hunting among predators are
relatively small (i.e., $\kappa<\kappa_1$ and $\alpha\leq\alpha_2$),
then the intensity of cooperative hunting is so weak that it can not sustain the survival of predators,
which thereby causes the survival of prey and the extinction of predators.
Theorem \ref{thmBT} indicates that the neighborhood $U$ of point $(\sigma_*,\alpha_*)$  can be divided into four regions, i.e., $U=\mathcal{SN}^+\cup\mathcal{SN}^-\cup\mathcal{H}\cup\mathcal{HL}\cup\mathcal{I}_1\cup\mathcal{I}_2\cup\mathcal{I}_3\cup\mathcal{I}_4$  with
\begin{eqnarray*}
\begin{array}{l}
\mathcal{I}_1:=\{(\sigma, \alpha)\in U:\alpha<
\alpha_*+f_{11}(\sigma-\sigma_*)+f_{12}(\sigma-\sigma_*)^2+O(|\sigma-\sigma_*|^3)\},\\
\mathcal{I}_2:=\{(\sigma, \alpha)\in U:
\alpha_*+f_{11}(\sigma-\sigma_*)+f_{12}(\sigma-\sigma_*)^2+O(|\sigma-\sigma_*|^3)
<\alpha
\\
\phantom{\mathcal{I}_1:=}
<
\alpha_*+f_{21}(\sigma-\sigma_*)+f_{22}(\sigma-\sigma_*)^2+O(|\sigma-\sigma_*|^3), \sigma>\sigma_*
\},\\
\mathcal{I}_3:=\{(\sigma, \alpha)\in U:
\alpha_*+f_{21}(\sigma-\sigma_*)+f_{22}(\sigma-\sigma_*)^2+O(|\sigma-\sigma_*|^3)<\alpha\\
\phantom{\mathcal{I}_1:=}
<\alpha_*+f_{31}(\sigma-\sigma_*)+f_{32}(\sigma-\sigma_*)^2+O(|\sigma-\sigma_*|^3)
, \sigma>\sigma_*
\},\\
\mathcal{I}_4:=\{(\sigma, \alpha)\in U:
\alpha>\alpha_*+f_{31}(\sigma-\sigma_*)+f_{32}(\sigma-\sigma_*)^2+O(|\sigma-\sigma_*|^3), \sigma\geq\sigma_*\}\\
\phantom{\mathcal{I}_1:=}
\cup\{(\sigma, \alpha)\in U:
\alpha>\alpha_*+f_{11}(\sigma-\sigma_*)+f_{12}(\sigma-\sigma_*)^2+O(|\sigma-\sigma_*|^3), \sigma<\sigma_*\}.\\
\end{array}
\end{eqnarray*}
The corresponding dynamical properties of system (\ref{(3)}) near $E_*$ for the parameters $\sigma$ and $\alpha$ in the neighborhood $U$ of point
$(\sigma_*,\alpha_*)$ are displayed in Table 1.
\begin{table}[htb]
\begin{tabular}{lll}
\multicolumn{3}{c}{Table 1. Dynamical behaviors near $E_*$.}\\
\hline
$(\sigma, \alpha)\in$  &  Positive equilibria and properties  & Closed or homoclinic orbits\\
\hline
$\mathcal{I}_{1}$  & No equilibria & No\\
$\mathcal{SN}^-$  &  $E_*$(saddle node)  &  No\\
$\mathcal{I}_2$ &   $E_1$(stable focus or node) $E_2$(saddle)  &No\\
$\mathcal{H}$ &   $E_1$(stable weak focus) $E_2$(saddle)  &  No\\
$\mathcal{I}_{3}$ &  $E_1$(unstable focus) $E_2$(saddle) & A stable limit cycle\\
$\mathcal{HL}$ &  $E_1$(unstable focus) $E_2$(saddle)  & A homoclinic orbit\\
$\mathcal{I}_{4}$ &  $E_1$(unstable focus or node) $E_2$(saddle)  & No\\
$\mathcal{SN}^+$  &  $E_*$(saddle node)  &  No\\
$(\sigma_*,\alpha_*)$ &  $E_*$(degenerate equilibrium, i.e., cusp)  &  No\\
\hline
\end{tabular}
\end{table}
The bifurcation diagram for Bogdanov-Takens bifurcation of codimension two and the corresponding phase portraits are shown in Figure \ref{Figure BT} for $h=0.5$ and $\kappa=1.2$.
The bifurcation diagram  is displayed in Figure \ref{Figure BT}(a).
The following explains the various phase portraits of Bogdanov-Takens bifurcation.
If $(\sigma,\alpha)=(0.35,17.1)\in\mathcal{I}_{1}$, then system  (\ref{(3)}) has no any positive equilibrium(see Figure \ref{Figure BT}(b)), which also shows that the prey can survive but the predators go extinct because the intensity of cooperative hunting is too weak to sustain the survival of predators.
If $(\sigma,\alpha)=(0.35,17.5)\in\mathcal{I}_{2}$, then system  (\ref{(3)}) has two positive equilibria arising from the saddle-node bifurcation, and one is the stable focus $E_1$ and the other one is the saddle $E_2$(see Figure \ref{Figure BT}(c)),
which means that the predators and their prey can coexist at stable coexistence equilibrium  or predator extinction equilibrium for different initial values.
If $(\sigma,\alpha)=(0.35,18)\in\mathcal{I}_{3}$, then system  (\ref{(3)}) has a stable limit cycle around the unstable focus $E_1$ induced by the Hopf bifurcation and the saddle $E_2$(see Figure \ref{Figure BT}(d)).
If $(\sigma,\alpha)=(0.305,19.6855)\in\mathcal{HL}$, then system  (\ref{(3)}) has the saddle $E_2$ and a homoclinic orbit around the unstable focus $E_1$(see Figure \ref{Figure BT}(e)).  The disappearance of the limit cycle and the appearance of the homoclinic orbit are due to the homoclinic bifurcation.
If $(\sigma,\alpha)=(0.35,18.6)\in\mathcal{I}_{4}$, then system  (\ref{(3)}) still has two positive equilibria, and more specifically one is the unstable focus $E_1$ and the other one is the saddle $E_2$(see Figure \ref{Figure BT}(f)).
Figures \ref{Figure BT}(d) and (f) show that
if the intensity of cooperative hunting is too large, then the predators also go extinct, because the system
undergoes Hopf bifurcation leading that the stable coexistence equilibrium loses stability and the limit
cycle emerges, or the system undergoes homoclinic bifurcation leading to the disappearance of the limit cycle.
The cooperative hunting therefore can turn detrimental to predators if the prey density is relatively small and the
intensity of cooperative hunting is too large because of the increased  predation pressure of predators.
Because the Hopf bifurcation at the weak focus $E_1$ of multiplicity three depends on more bifurcation parameters and even is very sensitive to the perturbation of parameters,
we display one limit cycle and two limit cycles around $E_1$ arising from the Hopf bifurcation through numerical simulations in Figure \ref{Figure Hopf}.

(i) If $h=0.5$, $\kappa=0.8$, $\alpha=54.902$ and $\sigma=0.68$, then system (\ref{(3)}) has four equilibria, i.e., the saddles $E_0$ and $E_2$, the stable node $E_{\kappa}$ and the unstable focus $E_1$ with the trace $T|_{E_1}=0.02$ and the first order Lyapunov quantity $l_1=-21.2827$, which indicates that system (\ref{(3)}) has a stable limit cycle around the unstable focus $E_1$ (see Figure \ref{Figure Hopf}(a));

(ii) If $h=0.45$, $\kappa=133.7629$, $\alpha=0.3555$ and $\sigma=2.319$, then system (\ref{(3)}) has two limit cycles around the unstable focus $E_1(1,2.3016)$
with the trace $T|_{E_1}=5.1927\times 10^{-25}$ and the first two Lyapunov quantities $l_1=-2.6823\times 10^{-19}$ and $l_2=1.2546\times 10^{-15}$.
Two orbits spiraling inward and outward respectively and the unstable focus $E_1$ form two annular regions.  From Poincar\'e-Bendixson Theorem(\cite{Zhang}), there are two closed orbits in the two annular regions.
We plot the two orbits starting from $P_1(1, 2)$ and $P_2(1,2.2)$ respectively in Figure \ref{Figure Hopf}(b).
However, it is very blurry to distinguish whether the two orbits starting from $P_1$ and $P_2$ spiral outward or inward as
$t\rightarrow +\infty$ from Figure \ref{Figure Hopf}(b).
Therefore, by zooming in the orbits near $P_1$ and $P_2$ in Figure \ref{Figure Hopf}(c) and \ref{Figure Hopf}(d) respectively, we obtain that the orbit from
$P_1$ spirals outward  while the orbit from $P_2$ spirals inward as $t\rightarrow +\infty$.
Therefore, there is an unstable limit cycle lying in the annular region bounded by the two orbits starting from $P_1$ and $P_2$,
and a stable limit cycle lying in the annular region bounded by the orbit starting from $P_2$ and the unstable focus $E_1$.
These numerical simulations indicate that
if only the environmental capacity of prey is relatively small but the intensity of cooperative hunting
is relatively large (i.e.,  $\kappa<\kappa_1$ and $\alpha>\alpha_2$),
then the cooperative hunting can  promote the probability of survival of predators, which implies that
both of predators and prey can coexist at the positive equilibrium or limit cycles for different initial values.
Therefore, the cooperative hunting can be beneficial to the predators and the system exhibits two types of bistability phenomenon,
i.e., either a stable coexistence equilibrium together with a stable predator extinction equilibrium or a stable limit cycle together with a stable predator extinction equilibrium.
This phenomenon actually is the Allee effect in the predators(\cite{Stephens}).
Our theoretical results indicate that the cooperative hunting among predators can lead to the changes of existence and stabilities for invariant sets, such as the equilibria, the limit cycles and the homoclinic orbit, which are induced by various bifurcations.
The cooperative hunting among predators therefore also has a destabilizing effect on the dynamics in the system with Holling III functional response
just as in the systems with Lotka-Volterra or Holling II or Holling IV functional responses respectively(\cite{Alves,Yao2020,ZhangZhang}).
By the bifurcation analysis, we provide some thresholds to control the dynamical behaviours of the predator-prey system, which are the
critical values to promote the coexistence of predators and their prey.
Hence, the analysis results reveal that appropriate intensity of cooperative hunting among predators is beneficial for the persistence of predators and the diversity of ecosystem.

\section*{ Declarations}

\section*{\large Funding}
This work is supported by the National Natural Science
Foundation of China (Grant Nos. 12101470, 11961023, 11701163) and
the Science Foundation of Wuhan Institute of Technology (Grant No. K2021077).

\section*{\large Conflict of interest}
The authors declare that they have no competing interests.

\section*{\large Ethics approval and consent to participate}
The research in this paper does not involve any ethical research.

\section*{\large Consent for publication}
All authors agree to publish this paper.


\bibliographystyle{plain}
\footnotesize\bibliography{references}
\newpage

\newpage
\section*{Appendix}
Coefficients in system \eqref{(3-5)} are defined as follows:
\begin{eqnarray*}
\left.
\begin{array}{l}
p_{001}:=-\frac{(a_{100}+b_{100}) b_{001} b_{010}}{b_{100} (a_{100}+b_{010})},~q_{010}:=a_{100}+b_{010},\\
p_{200}:=\frac{a_{020} b_{100}^3-a_{100} b_{010}^2 b_{200}+a_{100} b_{010} b_{100} b_{110}-a_{100} b_{020} b_{100}^2-a_{110} b_{010} b_{100}^2+a_{200} b_{010}^2 b_{100}}{(a_{100}+b_{010}) b_{010}b_{100}},\\
p_{020}:=\frac{(a_{020} b_{100}^3-a_{100}^3 b_{200}+a_{100}^2 a_{200} b_{100}-a_{100}^2 b_{100} b_{110}+a_{100} a_{110} b_{100}^2-a_{100} b_{020} b_{100}^2) b_{010}}{ (a_{100}+b_{010}) b_{100}^3},\\
p_{002}:=\frac{1}{ b_{100}^3 (a_{100}+b_{010})^3}\{b_{010} b_{001} (b_{010}-b_{100}) (a_{100}^2 b_{100} b_{101}-a_{100} a_{101} b_{100}^2-a_{100} b_{001} b_{010} b_{200}\\
\phantom{p_{002}:=}
+a_{100} b_{001} b_{100} b_{200}+a_{100} b_{010} b_{100} b_{101}-a_{101} b_{010} b_{100}^2+a_{200} b_{001} b_{010} b_{100}-a_{200} b_{001} b_{100}^2)\},\\
p_{110}:=\frac{-1}{b_{100}^2 (a_{100}+b_{010})}\{2 a_{020} b_{100}^3+2 a_{100}^2 b_{010} b_{200}-a_{100}^2 b_{100} b_{110}+a_{100} a_{110} b_{100}^2\\
\phantom{p_{002}:=}
-2 a_{100} a_{200} b_{010} b_{100}+a_{100} b_{010} b_{100} b_{110}-2 a_{100} b_{020} b_{100}^2-a_{110} b_{010} b_{100}^2\},\\
p_{011}:=\frac{ b_{010}}{(a_{100}+b_{010})^2 b_{100}^3}\{a_{011} a_{100} b_{100}^3+a_{011} b_{010} b_{100}^3-a_{100}^3 b_{100} b_{101}+a_{100}^2 a_{101} b_{100}^2\\
\end{array}
\right.
\end{eqnarray*}
Coefficients of center manifold $W^c$ in \eqref{(3-6)} are
\begin{eqnarray*}
\left.
\begin{array}{l}
a_2:=-\frac{a_{020} b_{100}^3-a_{110} b_{010} b_{100}^2+a_{200} b_{010}^2 b_{100}+b_{010}^3 b_{200}-b_{010}^2 b_{100} b_{110}+b_{010} b_{020} b_{100}^2}{(a_{100}+b_{010})^2 b_{010}^2},\\
b_2:=\frac{1}{(a_{100}+b_{010})^4 b_{010} b_{100}}\{a_{011} a_{100}^2 b_{100}^3+2 a_{011} a_{100} b_{010} b_{100}^3+a_{011} b_{010}^2 b_{100}^3+2 a_{020} a_{100} b_{001} b_{100}^3\\
\phantom{b_2:=}
+2 a_{020} b_{001} b_{100}^4-a_{100}^2 a_{101} b_{010} b_{100}^2-a_{100}^2 b_{010}^2 b_{100} b_{101}+a_{100}^2 b_{010} b_{011} b_{100}^2\\
\phantom{b_2:=}
-2 a_{100} a_{101} b_{010}^2 b_{100}^2-3 a_{100} a_{110} b_{001} b_{010} b_{100}^2+a_{100} a_{110} b_{001} b_{100}^3+4 a_{100} a_{200} b_{001} b_{010}^2 b_{100}\\
\phantom{b_2:=}
-2 a_{100} a_{200} b_{001} b_{010} b_{100}^2+4 a_{100} b_{001} b_{010}^3 b_{200}-3 a_{100} b_{001} b_{010}^2 b_{100} b_{110}\\
\phantom{b_2:=}
-2 a_{100} b_{001} b_{010}^2 b_{100} b_{200}+2 a_{100} b_{001} b_{010} b_{020} b_{100}^2+a_{100} b_{001} b_{010} b_{100}^2 b_{110}-2 a_{100} b_{010}^3 b_{100} b_{101}\\
\phantom{b_2:=}
+2 a_{100} b_{010}^2 b_{011} b_{100}^2-a_{101} b_{010}^3 b_{100}^2-a_{110} b_{001} b_{010}^2 b_{100}^2-a_{110} b_{001} b_{010} b_{100}^3+2 a_{200} b_{001} b_{010}^3 b_{100}\\
\phantom{b_2:=}
+2 b_{001} b_{010}^4 b_{200}-b_{001} b_{010}^3 b_{100} b_{110}-b_{001} b_{010}^2 b_{100}^2 b_{110}+2 b_{001} b_{010} b_{020} b_{100}^3-b_{010}^4 b_{100} b_{101}\\
\phantom{b_2:=}
+b_{010}^3 b_{011} b_{100}^2\},\\
\end{array}
\right.
\end{eqnarray*}
\begin{eqnarray*}
\left.
\begin{array}{l}
c_2:=\frac{-b_{001} }{(a_{100}+b_{010})^6 b_{100}^2}\{a_{011} a_{100}^3 b_{100}^3+2 a_{011} a_{100}^2 b_{010} b_{100}^3+a_{011} a_{100}^2 b_{100}^4+a_{011} a_{100} b_{010}^2 b_{100}^3\\
\phantom{b_2:=}
+2 a_{011} a_{100} b_{010} b_{100}^4+a_{011} b_{010}^2 b_{100}^4+2 a_{020} a_{100}^2 b_{001} b_{100}^3+4 a_{020} a_{100} b_{001} b_{100}^4+2 a_{020} b_{001} b_{100}^5\\
\phantom{b_2:=}
-2 a_{100}^3 a_{101} b_{010} b_{100}^2+a_{100}^3 a_{101} b_{100}^3-2 a_{100}^3 b_{010}^2 b_{100} b_{101}+a_{100}^3 b_{010} b_{011} b_{100}^2+a_{100}^3 b_{010} b_{100}^2 b_{101}\\
\phantom{b_2:=}
-5 a_{100}^2 a_{101} b_{010}^2 b_{100}^2+2 a_{100}^2 a_{101} b_{010} b_{100}^3-3 a_{100}^2 a_{110} b_{001} b_{010} b_{100}^2+a_{100}^2 a_{110} b_{001} b_{100}^3\\
\phantom{b_2:=}
+5 a_{100}^2 a_{200} b_{001} b_{010}^2 b_{100}-4 a_{100}^2 a_{200} b_{001} b_{010} b_{100}^2+a_{100}^2 a_{200} b_{001} b_{100}^3+5 a_{100}^2 b_{001} b_{010}^3 b_{200}\\
\phantom{b_2:=}
-3 a_{100}^2 b_{001} b_{010}^2 b_{100} b_{110}-4 a_{100}^2 b_{001} b_{010}^2 b_{100} b_{200}+2 a_{100}^2 b_{001} b_{010} b_{020} b_{100}^2\\
\phantom{b_2:=}
+a_{100}^2 b_{001} b_{010} b_{100}^2 b_{110}+a_{100}^2 b_{001} b_{010} b_{100}^2 b_{200}-5 a_{100}^2 b_{010}^3 b_{100} b_{101}+2 a_{100}^2 b_{010}^2 b_{011} b_{100}^2\\
\phantom{b_2:=}
+2 a_{100}^2 b_{010}^2 b_{100}^2 b_{101}+a_{100}^2 b_{010} b_{011} b_{100}^3-4 a_{100} a_{101} b_{010}^3 b_{100}^2+a_{100} a_{101} b_{010}^2 b_{100}^3\\
\end{array}
\right.
\end{eqnarray*}
Coefficients  in system \eqref{(3-9)} are defined as follows:
\begin{eqnarray*}
\left.
\begin{array}{l}
c_{00}(\epsilon_1, \epsilon_2):=\frac{-(h-1)^2 (h x_*^2-x_*^2+1)^2 x_*^4 \epsilon_1}{4\{h^2 x_*^2+(-2 x_*^2+1) h+x_*^2-2\}},~~
c_{10}(\epsilon_1, \epsilon_2):=\frac{1}{2} (h x_*^2-x_*^2-1) \epsilon_1 (h-1) x_*^3,\\
c_{01}(\epsilon_1, \epsilon_2):=-\frac{x_*^3 (h-1)^2 (h x_*^2-x_*^2+1)^2 \epsilon_1-2 h^2 x_*^2+(4 x_*^2-2) h-2 x_*^2+4}{2\{h^2 x_*^2+(-2 x_*^2+1) h+x_*^2-2\}},\\
c_{02}(\epsilon_1, \epsilon_2):=-\frac{x_*^3 (h-1)^2 (h x_*^2-x_*^2+1)^2 \epsilon_1-2 h^2 x_*^2+(4 x_*^2-2) h-2 x_*^2+4}{4 x_* (h^2 x_*^2-2 h x_*^2+x_*^2+h-2)},\\
c_{11}(\epsilon_1, \epsilon_2):=\frac{x_*^3 (h-1)^2 (h x_*^2-x_*^2-3) (h x_*^2-x_*^2+1)^2 \epsilon_1+2 (h x_*^2-x_*^2+3) (h^2 x_*^2-2 h x_*^2+x_*^2+h-2)}{2 x_* (h-1) (h x_*^2-x_*^2+1)^2},\\
c_{20}(\epsilon_1, \epsilon_2):=\frac{-(h^2 x_*^2-2 h x_*^2+x_*^2+h-2) \{x_*^3 (-1+h)^2 (x_*^4 h^2-2 x_*^2 (x_*^2+3) h+x_*^4+6 x_*^2-3) \epsilon_1+6 h^2 x_*^2+(-12 x_*^2+6) h+6 x_*^2-12\}}{4\{x_* (-1+h)^2 (h x_*^2-x_*^2+1)^2\}},\\
d_{00}(\epsilon_1, \epsilon_2):=\{(h x_*^3 (h-1)^2 (h^2 x_*^2-2 h x_*^2+x_*^2+h-2) (h x_*^2-x_*^2+1)^2 \epsilon_1 \epsilon_2+x_*^3 (h-1)^3 (h x_*^2-x_*^2+2)\\
\phantom{d_{10}(\epsilon_1, \epsilon_2):=}
\cdot (h x_*^2-x_*^2+1)^2 \epsilon_1-2 (h^2 x_*^2-2 h x_*^2+x_*^2+h-2)^2 \epsilon_2) x_*\}/\{4(h x_*^2-x_*^2+2) (h^2 x_*^2-2 h x_*^2\\
\phantom{d_{10}(\epsilon_1, \epsilon_2):=}
+x_*^2+h-2) (h-1)^2 (h x_*^2-x_*^2+1)\},\\
\end{array}
\right.
\end{eqnarray*}
Coefficients  in system \eqref{(3-10)} are defined as follows:
\begin{eqnarray*}
\left.
\begin{array}{l}
e_{00}(\epsilon_1, \epsilon_2) := -\{(c_{00}^4(\epsilon_1, \epsilon_2) c_{02}^2(\epsilon_1, \epsilon_2) d_{02}(\epsilon_1, \epsilon_2)+2 c_{00}^3(\epsilon_1, \epsilon_2) c_{01}^2(\epsilon_1, \epsilon_2) c_{02}(\epsilon_1, \epsilon_2) d_{02}(\epsilon_1, \epsilon_2)\\
\phantom{e_{00}(\epsilon_1, \epsilon_2):=}
+c_{00}^2(\epsilon_1, \epsilon_2)c_{01}^4(\epsilon_1, \epsilon_2) d_{02}(\epsilon_1, \epsilon_2)-c_{00}^2(\epsilon_1, \epsilon_2) c_{01}^3(\epsilon_1, \epsilon_2) c_{02}(\epsilon_1, \epsilon_2) d_{01}(\epsilon_1, \epsilon_2)-c_{00}(\epsilon_1, \epsilon_2)\\
\phantom{e_{00}(\epsilon_1, \epsilon_2):=}
\cdot c_{01}^5(\epsilon_1, \epsilon_2) d_{01}(\epsilon_1, \epsilon_2)+c_{01}^6(\epsilon_1, \epsilon_2) d_{00}(\epsilon_1, \epsilon_2)) (2 c_{00}^2(\epsilon_1, \epsilon_2) c_{02}^2(\epsilon_1, \epsilon_2)+2 c_{00}(\epsilon_1, \epsilon_2) c_{01}^2(\epsilon_1, \epsilon_2)\\
\phantom{e_{00}(\epsilon_1, \epsilon_2):=}
\cdot c_{02}(\epsilon_1, \epsilon_2)-c_{01}^4(\epsilon_1, \epsilon_2))\}/c_{01}^9(\epsilon_1, \epsilon_2),\\
e_{10}(\epsilon_1, \epsilon_2) := \frac{1}{c_{01}^9(\epsilon_1, \epsilon_2)}\{6 c_{00}^5(\epsilon_1, \epsilon_2) c_{01}(\epsilon_1, \epsilon_2) c_{02}^3(\epsilon_1, \epsilon_2) c_{11}(\epsilon_1, \epsilon_2) d_{02}(\epsilon_1, \epsilon_2)-12 c_{00}^5(\epsilon_1, \epsilon_2) c_{02}^4(\epsilon_1, \epsilon_2)\\
\phantom{e_{00}(\epsilon_1, \epsilon_2):=}
\cdot c_{10}(\epsilon_1, \epsilon_2) d_{02}(\epsilon_1, \epsilon_2)+2 c_{00}^4(\epsilon_1, \epsilon_2) c_{01}^3(\epsilon_1, \epsilon_2) c_{02}^3(\epsilon_1, \epsilon_2) d_{11}(\epsilon_1, \epsilon_2)+13 c_{00}^4(\epsilon_1, \epsilon_2) c_{01}^3(\epsilon_1, \epsilon_2) \\
\phantom{e_{00}(\epsilon_1, \epsilon_2):=}
\cdot c_{02}^2(\epsilon_1, \epsilon_2)c_{11}(\epsilon_1, \epsilon_2) d_{02}(\epsilon_1, \epsilon_2)-30 c_{00}^4(\epsilon_1, \epsilon_2) c_{01}^2(\epsilon_1, \epsilon_2) c_{02}^3(\epsilon_1, \epsilon_2) c_{10}(\epsilon_1, \epsilon_2) d_{02}(\epsilon_1, \epsilon_2)\\
\end{array}
\right.
\end{eqnarray*}
\begin{eqnarray*}
\left.
\begin{array}{l}
\phantom{e_{00}(\epsilon_1, \epsilon_2):=}
 \cdot c_{02}^2(\epsilon_1, \epsilon_2) c_{10}^2(\epsilon_1, \epsilon_2) d_{01}(\epsilon_1, \epsilon_2)-c_{01}^{10}(\epsilon_1, \epsilon_2) d_{20}(\epsilon_1, \epsilon_2)+c_{01}^9(\epsilon_1, \epsilon_2) c_{10}(\epsilon_1, \epsilon_2) d_{11}(\epsilon_1, \epsilon_2)\\
\phantom{d_{10}(\epsilon_1, \epsilon_2):=}
-c_{01}^9(\epsilon_1, \epsilon_2) c_{11}(\epsilon_1, \epsilon_2) d_{10}(\epsilon_1, \epsilon_2)+c_{01}^9(\epsilon_1, \epsilon_2) c_{20}(\epsilon_1, \epsilon_2) d_{01}(\epsilon_1, \epsilon_2)+2 c_{01}^8(\epsilon_1, \epsilon_2) c_{02}(\epsilon_1, \epsilon_2)\\
\phantom{e_{00}(\epsilon_1, \epsilon_2):=}
\cdot c_{10}(\epsilon_1, \epsilon_2)d_{10}(\epsilon_1, \epsilon_2)+2 c_{01}^8(\epsilon_1, \epsilon_2) c_{02}(\epsilon_1, \epsilon_2) c_{20}(\epsilon_1, \epsilon_2) d_{00}(\epsilon_1, \epsilon_2)-c_{01}^8(\epsilon_1, \epsilon_2) c_{10}^2(\epsilon_1, \epsilon_2)\\
\phantom{e_{00}(\epsilon_1, \epsilon_2):=}
\cdot d_{02}(\epsilon_1, \epsilon_2)-c_{01}^7(\epsilon_1, \epsilon_2)c_{02}(\epsilon_1, \epsilon_2) c_{10}^2(\epsilon_1, \epsilon_2) d_{01}(\epsilon_1, \epsilon_2)-2 c_{01}^7(\epsilon_1, \epsilon_2) c_{02}(\epsilon_1, \epsilon_2) c_{10}(\epsilon_1, \epsilon_2) \\
\phantom{e_{00}(\epsilon_1, \epsilon_2):=}
\cdot c_{11}(\epsilon_1, \epsilon_2) d_{00}(\epsilon_1, \epsilon_2)+2 c_{01}^6(\epsilon_1, \epsilon_2) c_{02}^2(\epsilon_1, \epsilon_2) c_{10}^2(\epsilon_1, \epsilon_2) d_{00}(\epsilon_1, \epsilon_2)\},\\
e_{02}(\epsilon_1, \epsilon_2) := -\frac{1}{c_{01}^9(\epsilon_1, \epsilon_2)}\{30 c_{00}^4(\epsilon_1, \epsilon_2) c_{02}^4(\epsilon_1, \epsilon_2) d_{02}(\epsilon_1, \epsilon_2)+60 c_{00}^3(\epsilon_1, \epsilon_2) c_{01}^2(\epsilon_1, \epsilon_2) c_{02}^3(\epsilon_1, \epsilon_2) d_{02}(\epsilon_1, \epsilon_2)\\
\end{array}
\right.
\end{eqnarray*}
Coefficients $g_{01}$, $g_{10}$, $h_{01}$, $h_{10}$, $f_{11}$, $f_{12}$, $f_{22}$ and $f_{32}$ in the proof of theorem \ref{thmBT}  are defined as follows:
\begin{eqnarray*}
\left.
\begin{array}{l}
g_{01} := -\frac{h^3 x_*^4-x_*^2 (3 x_*^2-7) h^2+(3 x_*^4-15 x_*^2+6) h-(x_*^2-6) (x_*^2-2)}{2 x_*^3 \{3 h^5 x_*^4-2 x_*^2 (7 x_*^2-3) h^4+(26 x_*^4-22 x_*^2+3) h^3-6 (x_*^2-1) (4 x_*^2-1) h^2+(11 x_*^4-18 x_*^2+3) h-2 x_*^2 (x_*^2-2)\}},\\
g_{10} := -\{h^{11} x_*^{16}-10 x_*^{14} (x_*^2-2) h^{10}+x_*^{12} (46 x_*^4-178 x_*^2+157) h^9-x_*^{10} (129 x_*^6-726 x_*^4+1202 x_*^2\\
\phantom{g_{10} := }
-628) h^8+x_*^8 (246 x_*^8-1801 x_*^6+4199 x_*^4-3856 x_*^2+1415) h^7-x_*^6 (336 x_*^{10}-3031 x_*^8+8949 x_*^6\\
\phantom{g_{10} := }
-10596 x_*^4+5950 x_*^2-1876) h^6+x_*^4 (336 x_*^{12}-3633 x_*^{10}+13045 x_*^8-18001 x_*^6+9172 x_*^4-3410 x_*^2\\
\phantom{g_{10} := }
+1455) h^5-x_*^2 (246 x_*^{14}-3143 x_*^{12}+13597 x_*^{10}-22065 x_*^8+7298 x_*^6+6986 x_*^4-1482 x_*^2-612) h^4\\
\phantom{g_{10} := }
+(129 x_*^{16}-1931 x_*^{14}+10077 x_*^{12}-20378 x_*^{10}+6947 x_*^8+21402 x_*^6-18657 x_*^4+2628 x_*^2+108) h^3\\
\phantom{g_{10} := }
-(x_*^2-6) (46 x_*^{14}-525 x_*^{12}+1889 x_*^{10}-1872 x_*^8-1994 x_*^6+4366 x_*^4-2004 x_*^2+144) h^2+(x_*^2-1) \\
\phantom{g_{10} := }
\cdot (x_*^2-2) (10 x_*^8-51 x_*^6+9 x_*^4+88 x_*^2-36) (x_*^2-6)^2 h-(x_*^2+1) (x_*-1)^2 (x_*+1)^2(x_*^2-2)^2\\
\end{array}
\right.
\end{eqnarray*}
\begin{eqnarray*}
\begin{array}{l}
h_{10} := - \{3 h^9 x_*^{12}-x_*^{10} (26 x_*^2-33) h^8+x_*^8 (101 x_*^4-256 x_*^2+132) h^7-x_*^6 (231 x_*^6-882 x_*^4+896 x_*^2\\
\phantom{g_{10} := }
-258) h^6+x_*^4 (343 x_*^8-1764 x_*^6+2683 x_*^4-1486 x_*^2+267) h^5-x_*^2 (343 x_*^{10}-2240 x_*^8+4595 x_*^6\\
\phantom{g_{10} := }
-3781 x_*^4+1244 x_*^2-141) h^4+(231 x_*^{12}-1848 x_*^{10}+4850 x_*^8-5424 x_*^6+2628 x_*^4-490 x_*^2+30) h^3\\
\phantom{g_{10} := }
-(x_*^2-1) (101 x_*^{10}-865 x_*^8+2277 x_*^6-2299 x_*^4+791 x_*^2-66) h^2+(x_*^2-2) (26 x_*^6-188 x_*^4+269 x_*^2\\
\phantom{g_{10} := }
-54) (x_*^2-1)^2 h-(3 (x_*^2-6)) (x_*^2-2)^2 (x_*^2-1)^3\}\{h^4 x_*^8-2 x_*^6 (2 x_*^2-7) h^3+x_*^4 (6 x_*^4-42 x_*^2+61) h^2\\
\phantom{g_{10} := }
-2 x_*^2 (x_*^2-1) (x_*^2-6) (2 x_*^2-7) h+(x_*^2-1)^2 (x_*^2-6)^2\}/\{8 (h^2 x_*^2-2 h x_*^2+x_*^2+h-2) (81 h^{12} x_*^{16}\\
\phantom{g_{10} := }
-216 x_*^{14} (4 x_*^2-3) h^{11}+108 x_*^{12} (39 x_*^4-58 x_*^2+21) h^{10}-24 x_*^{10} (517 x_*^6-1143 x_*^4+810 x_*^2-189) h^9\\
\phantom{g_{10} := }
+2 x_*^8 (12299 x_*^8-35916 x_*^6+37314 x_*^4-16740 x_*^2+2835) h^8-8 x_*^6 (4321 x_*^{10}-15619 x_*^8+21117 x_*^6\\
\end{array}
\end{eqnarray*}
\begin{eqnarray*}
\begin{array}{l}
f_{32}:=-2 \{274 h^{12} x_*^{14}-3 x_*^{12} (1046 x_*^2-721) h^{11}+x_*^{10} (16459 x_*^4-23663 x_*^2+7224) h^{10}-5 x_*^8 (10456 x_*^6\\
\phantom{f_{11} := }
-23478 x_*^4+15107 x_*^2-2653) h^9+5 x_*^6 (22401 x_*^8-69729 x_*^6+70765 x_*^4-26512 x_*^2+2898) h^8\\
\phantom{f_{11} := }
-x_*^4 (170508 x_*^{10}-688890 x_*^8+978220 x_*^6-584940 x_*^4+138320 x_*^2-9429) h^7+x_*^2 (189126 x_*^{12}\\
\phantom{f_{11} := }
-951006 x_*^{10}+1768415 x_*^8-1497050 x_*^6+573345 x_*^4-85903 x_*^2+3388) h^6+(-154008 x_*^{14}\\
\phantom{f_{11} := }
+936096 x_*^{12}-2185108 x_*^{10}+2450755 x_*^8-1349400 x_*^6+333434 x_*^4-29423 x_*^2+519) h^5\\
\phantom{f_{11} := }
+5 (x_*^2-2) (18276 x_*^{12}-94866 x_*^{10}+184195 x_*^8-164238 x_*^6+66408 x_*^4-10445 x_*^2+429) h^4\\
\phantom{f_{11} := }
-5 (7706 x_*^{10}-33655 x_*^8+53396 x_*^6-36238 x_*^4+9648 x_*^2-723) (x_*^2-2)^2 h^3+(10959 x_*^8\\
\phantom{f_{11} := }
-39561 x_*^6+50471 x_*^4-23960 x_*^2+3165) (x_*^2-2)^3 h^2-(1888 x_*^6-5514 x_*^4+5641 x_*^2-1470)\\
\end{array}
\end{eqnarray*}
Coefficients in system \eqref{(4-5)} are defined as follows:
\begin{eqnarray*}
\begin{array}{l}
a_{20} := -\frac{4 h \kappa-8 h x_1-\kappa}{4x_1 (h-1) (\kappa-x_1) \omega},~~
a_{11} := \frac{\{(2 h+1) (h x_1^2-x_1^2+3) \kappa-4 x_1 (h+1) (h x_1^2-x_1^2+2)\} \kappa}{2x_1 (h-1) (\kappa-x_1) \omega^2},\\
a_{02} := -\frac{(h x_1^2-x_1^2+1)^2 (3 \kappa-8 x_1) \kappa^2}{4x_1 (h-1) (\kappa-x_1) \omega^3},~~
a_{30} := -\frac{(2 h \kappa-2 h x_1-\kappa) (\kappa-4 x_1) h}{4 x_1^2 \kappa \omega  (\kappa-x_1)^2  (h x_1^2-x_1^2+1) (h-1)^2 },\\
a_{21} := -\frac{\{(h^2 x_1^2-2 h x_1^2+x_1^2-h-2) \kappa+2 x_1 (h+1) (h x_1^2-x_1^2+2)\} (2 h \kappa-2 h x_1-\kappa)}{4 x_1^2 \omega^2 (h x_1^2-x_1^2+1) (\kappa-x_1)^2 (h-1)^2} ,\\
a_{12} := \frac{\kappa (2 h \kappa-2 h x_1-\kappa) \{(5 h^2 x_1^2-7 h x_1^2+2 x_1^2+9 h) \kappa-4 x_1 (2 h^2 x_1^2-3 h x_1^2+x_1^2+3 h)\}}{4 x_1^2 \omega^3 (\kappa-x_1)^2 (h-1)^2} ,\\
a_{03} := -\frac{\kappa^2 (2 h \kappa-2 h x_1-\kappa) (h x_1^2-x_1^2+1) \{(3 h^2 x_1^2-4 h x_1^2+x_1^2+9 h+2) \kappa-2 x_1 (3 h^2 x_1^2-4 h x_1^2+x_1^2+6 h+2)\}}{4 x_1^2 \omega^4 (\kappa-x_1)^2 (h-1)^2},\\
a_{40} := \frac{(2 h \kappa-2 h x_1-\kappa)^2 h}{8x_1^2 \omega (\kappa-x_1)^3 (h-1)^3 \kappa^2 (h x_1^2-x_1^2+1)^2},~~
a_{31} := -\frac{(2 h \kappa-2 h x_1-\kappa)^2 (-2 x_1+\kappa) h}{4x_1^3 \kappa (\kappa-x_1)^3 (h-1)^3 \omega^2 (h x_1^2-x_1^2+1)},\\
a_{22} := \frac{(2 h \kappa-2 h x_1-\kappa)^2 \{(3 h^2 x_1^2-3 h x_1^2+5 h+1) \kappa-x_1 (9 h^2 x_1^2-9 h x_1^2+11 h+2)\}}{4x_1^3 (\kappa-x_1)^3 (h-1)^3 \omega^3 (h x_1^2-x_1^2+1)},\\
a_{13} := -\frac{(2 h \kappa-2 h x_1-\kappa)^2 \{(3 h^2 x_1^2-3 h x_1^2+7 h+2) \kappa-2 x_1 (5 h^2 x_1^2-5 h x_1^2+7 h+2)\} \kappa}{4x_1^3 (\kappa-x_1)^3 (h-1)^3 \omega^4},\\
\end{array}
\end{eqnarray*}
Functions $f_1$ and $f_2$ in \eqref{(4-6)} are defined as follows:
\begin{eqnarray*}
\begin{array}{l}
f_1:=
x_1 (2 h-1) (h-1) (3 h^2 x_1^2-4 h x_1^2+x_1^2+15 h-6) \kappa^3+\{-22 h^4 x_1^4+5 x_1^2 (11 x_1^2-18) h^3+(-45 x_1^4\\
\phantom{f_1:=}
+137 x_1^2+12) h^2+(13 x_1^4-53 x_1^2-18) h-x_1^2 (x_1^2-6)\} \kappa^2+2 h x_1 \{14 h^3 x_1^4-x_1^2 (30 x_1^2-43) h^2\\
\phantom{f_1:=}
+(18 x_1^4-44 x_1^2-15) h-2 x_1^4+x_1^2+20\} \kappa-2 x_1^2 h \{5 h^3 x_1^4-x_1^2 (7 x_1^2-11) h^2+(-x_1^4+x_1^2-10) h\\
\phantom{f_1:=}
+3 (x_1^2-2)^2\},\\
f_2:=
-x_1 (2 h-1) (h-1) \{2301 x_1^{10} h^9-3 x_1^8 (5377 x_1^2-8637) h^8+x_1^6 (49247 x_1^4-157848 x_1^2+102738) h^7\\
\phantom{f_1:=}
-x_1^4 (85593 x_1^6-409636 x_1^4+532904 x_1^2-181134) h^6+x_1^2 (92833 x_1^8-587900 x_1^6+1138738 x_1^4\\
\phantom{f_1:=}
-778814 x_1^2+146961) h^5+(-64671 x_1^{10}+506010 x_1^8-1285336 x_1^6+1300160 x_1^4-490965 x_1^2\\
\phantom{f_1:=}
+44955) h^4+(28621 x_1^{10}-264096 x_1^8+813975 x_1^6-1052906 x_1^4+557967 x_1^2-94122) h^3+(-7627 x_1^{10}\\
\end{array}
\end{eqnarray*}
\begin{eqnarray*}
\begin{array}{l}
\phantom{f_1:=}
-1135266) h^5+(-1682050 x_1^{12}+13652892 x_1^{10}-39479708 x_1^8+49536546 x_1^6-19890133 x_1^4\\
\phantom{f_1:=}
-5819287 x_1^2+3120294) h^4+(198292 x_1^{12}-1807982 x_1^{10}+5875322 x_1^8-8298632 x_1^6+2668196 x_1^4\\
\phantom{f_1:=}
+3520094 x_1^2-1700980) h^3+(3781 x_1^{12}-19200 x_1^{10}-24646 x_1^8+401202 x_1^6-872746 x_1^4+878328 x_1^2\\
\phantom{f_1:=}
-492776) h^2+(-2184 x_1^{12}+17849 x_1^{10}-31474 x_1^8-97638 x_1^6+421964 x_1^4-549880 x_1^2+253200) h\\
\phantom{f_1:=}
+3 x_1^2 (x_1-1) (x_1+1) (45 x_1^2-206) (x_1^2-2)^3\} \kappa^3+2 h x_1^4 \{138307 h^{10} x_1^{12}-x_1^{10} (969619 x_1^2-1165484) h^9\\
\phantom{f_1:=}
+x_1^8 (2938428 x_1^4-7071976 x_1^2+3542569) h^8-x_1^6 (5012024 x_1^6-18080215 x_1^4+18272302 x_1^2\\
\phantom{f_1:=}
-4690088) h^7+x_1^4 (5238142 x_1^8-25093853 x_1^6+38189061 x_1^4-20053146 x_1^2+2244877) h^6\\
\phantom{f_1:=}
-x_1^2 (3409686 x_1^{10}-20212255 x_1^8+40816259 x_1^6-32246049 x_1^4+7093245 x_1^2+427780) h^5\\
\end{array}
\end{eqnarray*}
The computation of algebraic varieties $V(f_1, f_2)$ and $V(f_1, f_2, f_3)$ in the proof of Theorem \ref{thm4} is displayed as follows:

In order to discuss the algebraic varieties $V(f_1, f_2)$ and $V(f_1, f_2, f_3)$ in $\mathcal{P}$, we first compute the following two resultants by eliminating $\kappa$.
\begin{eqnarray*}
\begin{array}{l}
\mathcal{R}_1:=\mbox{res}(f_1, f_2, \kappa)
=-262144 h^4x_1^{13}(2h-1)^2(h-1)^2\{(h-1)x_1^2+6\}^2\\
\phantom{\mathcal{R}_1:=\mbox{res}(f_1, f_2, \kappa)=}
\cdot
\{(h-1)x_1^2+1\}^6\mathcal{R}_{11}^4\mathcal{R}_{12}^2\mathcal{R}_{13}\mathcal{R}_{14},\\
\mathcal{R}_2:=\mbox{res}(f_1, f_3, \kappa)
=67108864000 h^5x_1^{21}(2h-1)^2(h-1)^2\{(h-1)x_1^2+6\}^3\\
\phantom{\mathcal{R}_1:=\mbox{res}(f_1, f_2, \kappa)=}
\cdot
\{(h-1)x_1^2+1\}^8\mathcal{R}_{11}^6\mathcal{R}_{12}^2\mathcal{R}_{13}\mathcal{R}_{21},
\end{array}
\label{(4-7)}
\end{eqnarray*}
where
\begin{eqnarray*}
\begin{array}{l}
\mathcal{R}_{11}:=x_1^2h^2+(1-2 x_1^2) h+x_1^2-2,\\
\mathcal{R}_{12}:=3 x_1^4 h^3-2 x_1^2 (4 x_1^2-3) h^2+(x_1^2-1) (7 x_1^2-3) h-2 x_1^2 (x_1^2-2),\\
\mathcal{R}_{13}:=4 h^3 x_1^2+2 (2-5 x_1^2) h^2+3 (3 x_1^2-2) h-3 x_1^2+6,\\
\mathcal{R}_{14}:=
138648 h^{25} x_1^{28}-8 x_1^{26} (347308 x_1^2-323343) h^{24}+2 x_1^{24} (13317143 x_1^4-24533740 x_1^2+10858952) h^{23}\\
\phantom{\mathcal{R}_{14}:=}
-2 x_1^{22} (81237707 x_1^6-222104059 x_1^4+195029732 x_1^2-55059488) h^{22}2 x_1^{20} (353906169 x_1^8\\
\phantom{\mathcal{R}_{14}:=}
-1276176758 x_1^6+1667704576 x_1^4-943015444 x_1^2+186493860) h^{21}-2 x_1^{18} (1171201098 x_1^{10}\\
\phantom{\mathcal{R}_{14}:=}
-5220522111 x_1^8+9025299990 x_1^6-7678241166 x_1^4+3086958356 x_1^2-444040892) h^{20}\\
\end{array}
\end{eqnarray*}
\begin{eqnarray*}
\begin{array}{l}
\phantom{\mathcal{R}_{14}:=}
+x_1^{16} (6115145446 x_1^{12}-32329653117 x_1^{10}+69311955226 x_1^8-78988428536 x_1^6+48612471694 x_1^4\\
\phantom{\mathcal{R}_{14}:=}
-14411511192 x_1^2+1534377208) h^{19}-x_1^{14} (12908775390 x_1^{14}-78637602251 x_1^{12}+200633635357 x_1^{10}\\
\phantom{\mathcal{R}_{14}:=}
-287614352650 x_1^8+241912359062 x_1^6-111534423774 x_1^4+24703132992 x_1^2-1962921024) h^{18}\\
\phantom{\mathcal{R}_{14}:=}
+x_1^{12} (22405428498 x_1^{16}-153896999093 x_1^{14}+453952038813 x_1^{12}-787142692781 x_1^{10}\\
\phantom{\mathcal{R}_{14}:=}
+852263742554 x_1^8-546616392080 x_1^6+190197279598 x_1^4-31606296800 x_1^2+1879158536) h^{17}\\
\phantom{\mathcal{R}_{14}:=}
-x_1^{10} (32338485680 x_1^{18}-246174990571 x_1^{16}+821327523092 x_1^{14}-1677336359863 x_1^{12}\\
\phantom{\mathcal{R}_{14}:=}
+2255413901801 x_1^{10}-1898044671530 x_1^8+927895706398 x_1^6-243469724922 x_1^4+30369710208 x_1^2\\
\phantom{\mathcal{R}_{14}:=}
-1347380840) h^{16}
+x_1^4 (34062749604 x_1^{24}-325457107196 x_1^{22}+1421716056334 x_1^{20}\\
\end{array}
\end{eqnarray*}
\begin{eqnarray*}
\begin{array}{l}
\phantom{\mathcal{R}_{14}:=}
+23771039392 x_1^6-6435894848 x_1^4+962645248 x_1^2-63446784) h^4+(34782 x_1^{28}+50411 x_1^{26}\\
\phantom{\mathcal{R}_{14}:=}
-7350288 x_1^{24}+71097240 x_1^{22}-343921815 x_1^{20}+1063754991 x_1^{18}-2296652137 x_1^{16}\\
\phantom{\mathcal{R}_{14}:=}
+3580087278 x_1^{14}-3987726020 x_1^{12}+3000013528 x_1^{10}-1311065328 x_1^8+176227680 x_1^6\\
\phantom{\mathcal{R}_{14}:=}
+66846912 x_1^4-28947072 x_1^2+4068864) h^3+x_1^2 (234 x_1^{26}+19883 x_1^{24}-610007 x_1^{22}\\
\phantom{\mathcal{R}_{14}:=}
+6943918 x_1^{20}-41230349 x_1^{18}+141746221 x_1^{16}-286341716 x_1^{14}+291617852 x_1^{12}+16077984 x_1^{10}\\
\phantom{\mathcal{R}_{14}:=}
-419860560 x_1^8+489026304 x_1^6-246756672 x_1^4+63037440 x_1^2-7368192) h^2-3 x_1^4 (x_1^2-2)\\
\phantom{\mathcal{R}_{14}:=}
\cdot(18 x_1^{22}+1075 x_1^{20}-28307 x_1^{18}+284045 x_1^{16}-1583807 x_1^{14}+5537788 x_1^{12}-12618820 x_1^{10}\\
\phantom{\mathcal{R}_{14}:=}
+18427440 x_1^8-16420176 x_1^6+8135424 x_1^4-2461248 x_1^2+387072) h-9 x_1^6 (x_1^2-6) (13 x_1^4-60 x_1^2\\
\phantom{\mathcal{R}_{14}:=}
+24) (x_1^2-3)^2 (x_1^2-2)^5,\\
\end{array}
\end{eqnarray*}
the polynomial $\mathcal{R}_{21}$ of 1193 terms is too long to display.
It is obvious that the two factors $(h-1)x_1^2+6$ and $(h-1)x_1^2+1$ in $\mathcal{R}_1$ and $\mathcal{R}_2$ are positive for $(\kappa,h,x_1)\in\mathcal{P}$.
If $h=\frac{1}{2}$, then $\mathcal{R}_1=\mathcal{R}_2=0$ and the parameter condition $\mathcal{P}$ becomes
\begin{eqnarray*}
\begin{array}{l}
\mathcal{P}=\{(\kappa,h,x_1)\in\mathbb{R}_+^3:
\kappa>\frac{2 x_1 (x_1^2-4)}{x_1^2-6}, h=\frac{1}{2},x_1<\sqrt{2}
\}.
\end{array}
\end{eqnarray*}
Substituting $h=\frac{1}{2}$ into $f_1$ and $f_2$  and computing the resultant of $f_1$ and $f_2$ by eliminating variable $\kappa$ we obtain that $f_1$ and $f_2$ have no common zero for $(\kappa,h,x_1)\in\mathcal{P}$, because the essential condition of $f_1$ and $f_2$ having common zeros for $(\kappa,h,x_1)\in\mathcal{P}$ is that their resultant has zero $x_1\doteq 0.8168$ or $x_1\doteq 0.8966$ for $0<x_1<\sqrt{2}$, but if $x_1\doteq 0.8168$, then $f_1$ has no zero for $\kappa>\frac{2 x_1 (x_1^2-4)}{x_1^2-6}\doteq 1.021$ and if
$x_1\doteq 0.8966$, then $f_1$ has no zero for $\kappa>\frac{2 x_1 (x_1^2-4)}{x_1^2-6}\doteq 1.103$. Therefore, we conclude that $V(f_1, f_2)=\emptyset$ if $h=\frac{1}{2}$.
Through the discussion of the monotonic interval partition of $\mathcal{R}_{11}$, $\mathcal{R}_{12}$ and $\mathcal{R}_{13}$ and their sign of function value
at the endpoint of interval, we can have that $\mathcal{R}_{11}<0$, $\mathcal{R}_{12}>0$ and $\mathcal{R}_{13}>0$  for $(\kappa,h,x_1)\in\mathcal{P}$.
The above analysis indicates that to determine whether the three polynomials $f_1$, $f_2$, $f_3$ have common zeros
is reduced to consider whether $\mathcal{R}_{14}$ and $\mathcal{R}_{21}$ have common zeros for $(\kappa,h,x_1)\in\mathcal{P}$.
We compute the resultant of $\mathcal{R}_{14}$ and $\mathcal{R}_{21}$ by eliminating $h$ as follows
\begin{eqnarray*}
\begin{array}{l}
\mathcal{R}_3:=\mbox{res}(\mathcal{R}_{14}, \mathcal{R}_{21}, h)
=C_1 x_1^{1040}(x_1^2-1)(x_1^2-12)(x_1^2-2)^{16}(x_1^2-3)^2\\
\phantom{\mathcal{R}_3:=\mbox{res}(\mathcal{R}_{14}, \mathcal{R}_{21}, h)=}
\cdot
(x_1^2-6)^{12} \mathcal{R}_{31} \mathcal{R}_{32} \mathcal{R}_{33} \mathcal{R}_{34} \mathcal{R}_{35}\mathcal{R}_{36}\mathcal{R}_{37},
\end{array}
\end{eqnarray*}
where
\begin{eqnarray*}
\begin{array}{l}
\mathcal{R}_{31}:=x_1^4-32 x_1^2+6,\\
\mathcal{R}_{32}:=11 x_1^4-5340 x_1^2+121500,\\
\mathcal{R}_{33}:=35 x_1^6-1718 x_1^4+33512 x_1^2-363048,\\
\mathcal{R}_{34}:=275 x_1^8-9465 x_1^6+78603 x_1^4-32723 x_1^2+2550,\\
\mathcal{R}_{35}:=47933 x_1^{10}-278074 x_1^8+616501 x_1^6-663380 x_1^4+364512 x_1^2-87552,
\end{array}
\end{eqnarray*}
both the polynomial $\mathcal{R}_{36}$ of 92 terms and
the polynomial $\mathcal{R}_{37}$ of 448 terms with complicated coefficients
are too long to display and $C_1$ is a nonzero constant.
Solving the exact solutions $x_1$ of the third to the ninth factors of $\mathcal{R}_3$, and substituting these solutions into
$\mathcal{R}_{14}$ and $\mathcal{R}_{21}$ and then computing the resultants of $\mathcal{R}_{14}$ and $\mathcal{R}_{21}$ respectively by eliminating variable $h$, we obtain that all of these resultants are nonzero constants, which implies  that $\mathcal{R}_{14}$ and $\mathcal{R}_{21}$ have no common zero. Therefore,
$V(f_1, f_2, f_3)=\emptyset$ for $(\kappa,h,x_1)\in\mathcal{P}$ if the seven factors vanish.
For example, when the third factor of $\mathcal{R}_3$ vanishes, i.e., $x_1^2-1=0$, then $\mathcal{R}_{14}$ and $\mathcal{R}_{21}$ are reduced to
\begin{eqnarray*}
\begin{array}{l}
\mathcal{R}_{14}=(2 h-1)\mathcal{\tilde{R}}_{14},~~
\mathcal{R}_{21}=h(2 h-1)\mathcal{\tilde{R}}_{21}
\end{array}
\end{eqnarray*}
by substituting $x_1=1$ into $\mathcal{R}_{14}$ and $\mathcal{R}_{21}$,
where the unary polynomials $\mathcal{\tilde{R}}_{14}$ and $\mathcal{\tilde{R}}_{21}$ are given as follows
\begin{eqnarray*}
\begin{array}{l}
\mathcal{\tilde{R}}_{14}:=69324 h^{24}-61198 h^{23}-388244 h^{22}+701986 h^{21}-10736604 h^{20}+53976423 h^{19}-50833424 h^{18}\\
\phantom{\mathcal{R}_{14}:=}
-229718263 h^{17}+602774491 h^{16}-184712481 h^{15}-933829705 h^{14}+903723219 h^{13}\\
\phantom{\mathcal{R}_{14}:=}
+586788847 h^{12}-1112956059 h^{11}+35349657 h^{10}+534362021 h^9-51678155 h^8-248638382 h^7\\
\phantom{\mathcal{R}_{14}:=}
+115458013 h^6-8638340 h^5-2637476 h^4+4917126 h^3-4281924 h^2+1010208 h-4140,\\
\mathcal{\tilde{R}}_{21}:=91163350713139920 h^{52}-805109298179470056 h^{51}+2909271328946802624 h^{50}\\
\phantom{\mathcal{R}_{14}:=}
-3966082520625814664 h^{49}-23821713664311131384 h^{48}+242718540353932841692 h^{47}\\
\phantom{\mathcal{R}_{14}:=}
-1175023543086646733512 h^{46}+3424517787203314151826 h^{45}-4836923647710710501490 h^{44}\\
\phantom{\mathcal{R}_{14}:=}
-5637384268261131317470 h^{43}+43949791872619135643068 h^{42}\\
\phantom{\mathcal{R}_{14}:=}
-87400906068290274236886 h^{41}+26891242507868782630535 h^{40}\\
\phantom{\mathcal{R}_{14}:=}
+205621916819307542043553 h^{39}-335443721511666216324658 h^{38}\\
\phantom{\mathcal{R}_{14}:=}
-92330162240822034185269 h^{37}+775164905117747736475569 h^{36}\\
\phantom{\mathcal{R}_{14}:=}
-544447151051718791913610 h^{35}-712992722105920461334142 h^{34}\\
\phantom{\mathcal{R}_{14}:=}
+1231016018962378963979452 h^{33}-278932955966676127369675 h^{32}\\
\phantom{\mathcal{R}_{14}:=}
-355724800016614007199171 h^{31}+186704316968315124054102 h^{30}\\
\phantom{\mathcal{R}_{14}:=}
-1171409641360001265730321 h^{29}+2348814880571718796979420 h^{28}\\
\phantom{\mathcal{R}_{14}:=}
-79355226027089714816367 h^{27}-3724011186540134019966188 h^{26}\\
\phantom{\mathcal{R}_{14}:=}
+2774553653856813463258615 h^{25}+2161055945839090313433308 h^{24}\\
\phantom{\mathcal{R}_{14}:=}
-3604529673079775821639913 h^{23}+25235012986739699452236 h^{22}\\
\phantom{\mathcal{R}_{14}:=}
+2499372014516246587276039 h^{21}-1085304067185267330282786 h^{20}\\
\phantom{\mathcal{R}_{14}:=}
-894133262608692418042529 h^{19}+852939119210961601455628 h^{18}\\
\phantom{\mathcal{R}_{14}:=}
+30787919904210915045177 h^{17}-287547641986943296692551 h^{16}\\
\phantom{\mathcal{R}_{14}:=}
+84283427911435646408458 h^{15}+35723602437214754722238 h^{14}\\
\phantom{\mathcal{R}_{14}:=}
-23401145020831097780578 h^{13}+2141335582526414031207 h^{12}+1085571912987773743919 h^{11}\\
\phantom{\mathcal{R}_{14}:=}
-716794394571895703714 h^{10}+438590757874883392947 h^9-40295446973462531301 h^8\\
\phantom{\mathcal{R}_{14}:=}
-69975459688916292390 h^7+17346114210008893446 h^6+3484538717411476278 h^5\\
\phantom{\mathcal{R}_{14}:=}
-958350103974544176 h^4-149637655598055228 h^3+20639692772757252 h^2\\
\phantom{\mathcal{R}_{14}:=}
+3035472663530160 h-12485116569600.
\end{array}
\end{eqnarray*}
The above discussion indicates that $f_1$ and $f_2$ have no common zero for $(\kappa,h,x_1)\in\mathcal{P}$ if $h=\frac{1}{2}$. Thus,
we further compute the resultant of $\mathcal{\tilde{R}}_{14}$ and $\mathcal{\tilde{R}}_{21}$ by eliminating variable $h$, which is a positive constant.
Therefore,  the polynomials $f_1$, $f_2$ and $f_3$ have no common zero for $(\kappa,h,x_1)\in\mathcal{P}$ if the third factor of $\mathcal{R}_3$ vanishes.
For the rest of  $\mathcal{R}_{ij}$ in $\mathcal{R}_3$,  the numerical calculation indicates that they have positive common zeros $x_1$, but either $\mathcal{R}_{14}$ and $\mathcal{R}_{21}$ or $f_1$, $f_2$ and $f_3$ have no common zero located in $\mathcal{P}$.
Therefore, we conclude that the algebraic variety $V(f_1, f_2, f_3)$ is empty in $\mathcal{P}$ if the remaining five factors $\mathcal{R}_{ij}$ vanish.
For example, when the factor $\mathcal{R}_{36}$ vanishes, then there are fifteen positive zeros $x_1$ of $\mathcal{R}_{36}=0$ such that either $\mathcal{\tilde{R}}_{14}$ and $\mathcal{\tilde{R}}_{21}$ or $f_1$, $f_2$ and $f_3$ have no common zero in $\mathcal{P}$, which implies that
the polynomials $f_1$, $f_2$ and $f_3$ have no common zero for $(\kappa,h,x_1)\in\mathcal{P}$ if the factor $\mathcal{R}_{36}$ vanishes.
Hence, we obtain that  the algebraic variety $V(f_1, f_2, f_3)=\emptyset$ for $(\kappa,h,x_1)\in\mathcal{P}$.
The conclusion shows that it is impossible to have four limit cycles around equilibrium $E_1$ induced by the Hopf bifurcation.
Therefore, the best possibility is that there are  three limit cycles around equilibrium $E_1$ induced by the Hopf bifurcation,
which requires the conditions $f_1=f_2=0$ and $f_3\neq0$.
By taking the special case $x_1=\frac{1}{2}$,  we can prove that there are three limit cycles around $E_1$ arising from the Hopf bifurcation.
For $x_1=\frac{1}{2}$, $f_1$, $f_2$ and $f_3$ are reduced to
\begin{eqnarray*}
\begin{array}{l}
\tilde{f}_1:=\frac{1}{8} (h-1) (2 h-1) (3 h^2+56 h-23) \kappa^3+(-\frac{11}{8} h^4-\frac{305}{16} h^3+\frac{695}{16} h^2-\frac{487}{16} h+\frac{23}{16}) \kappa^2\\
\phantom{\tilde{f}_1:=}
+\frac{1}{8} h (7 h^3+71 h^2-199 h+161) \kappa-\frac{1}{32} h (5 h^3+37 h^2-157 h+147),\\
\tilde{f}_2:=-\frac{1}{2048} (h-1) (2 h-1) (2301 h^9+87513 h^8+1061663 h^7+4619063 h^6+3738961 h^5-15048887 h^4\\
\phantom{\tilde{f}_2:=}
-8931523 h^3-1084075 h^2+12144918 h-3824494) \kappa^7+(\frac{57936947}{1024} h^5+\frac{174041}{512} h^{10}-\frac{1639627}{512} h^3\\
\phantom{\tilde{f}_2:=}
-\frac{74349967}{1024} h^6-\frac{10847359}{1024} h^7+\frac{56596193}{1024} h^4-\frac{80025471}{1024} h^2+\frac{42682305}{1024} h+\frac{2925805}{256} h^8-\frac{4638393}{1024}+\frac{457603}{128} h^9\\
\phantom{\tilde{f}_2:=}
+\frac{10121}{1024} h^{11}) \kappa^6+(\frac{256640367}{2048} h^6-\frac{247755109}{4096} h+\frac{29067783}{8192}-\frac{18989}{1024} h^{11}+\frac{20230725}{1024} h^7+\frac{35055051}{512} h^3\\
\phantom{\tilde{f}_2:=}
-\frac{24491797}{4096} h^9-\frac{150089393}{8192} h^8+\frac{234734227}{2048} h^2-\frac{19163}{32} h^{10}-\frac{199332063}{2048} h^5-\frac{610456751}{4096} h^4) \kappa^5+(\frac{183847465}{4096} h\\
\phantom{\tilde{f}_2:=}
-\frac{43934037}{2048} h^7-\frac{1914119}{16} h^6-\frac{977813437}{8192} h^3+\frac{130383063}{8192} h^8-\frac{8371739}{8192}+\frac{405153609}{4096} h^5+\frac{385131499}{2048} h^4\\
\phantom{\tilde{f}_2:=}
+\frac{4754413}{8192} h^{10}+\frac{11262195}{2048} h^9+\frac{157025}{8192} h^{11}-\frac{733590357}{8192} h^2) \kappa^4+(-\frac{96051}{8192} h^{11}-\frac{24498235}{8192} h^9-\frac{143842415}{8192} h\\
\phantom{\tilde{f}_2:=}
-\frac{2160746711}{16384} h^4+\frac{1317216027}{32768} h^2+\frac{780930145}{8192} h^3+\frac{58938633}{4096} h^7-\frac{262458887}{32768} h^8-\frac{258362059}{4096} h^5+\frac{1119250459}{16384} h^6\\
\phantom{\tilde{f}_2:=}
+\frac{2404773}{32768}-\frac{10913089}{32768} h^{10}) \kappa^3+\frac{1}{32768} h (138307 h^{10}+3692317 h^9+31331628 h^8+75117636 h^7\\
\phantom{\tilde{f}_2:=}
-192825126 h^6-765791114 h^5+813588316 h^4+1750226516 h^3-1330168517 h^2-349407275 h\\
\phantom{\tilde{f}_2:=}
+105233184) \kappa^2-\frac{1}{65536} h (53965 h^{10}+1356916 h^9+10816453 h^8+21759360 h^7-88466518 h^6\\
\phantom{\tilde{f}_2:=}
-291471896 h^5+362162434 h^4+774666976 h^3-590003991 h^2-113792700 h+10739673) \kappa\\
\phantom{\tilde{f}_2:=}
+\frac{1}{65536} h^2 (4375 h^9+103643 h^8+765608 h^7+1090696 h^6-8785666 h^5-23893626 h^4+35088560 h^3\\
\phantom{\tilde{f}_2:=}
+73348336 h^2-54415725 h-10739673),\\
\tilde{f}_3:=\frac{1}{524288} (h-1)(2 h-1) (140805675 h^{16}+7244615014 h^{15}+139158913678 h^{14}+1265156179722 h^{13}\\
\phantom{\tilde{f}_2:=}
+6427647061938 h^{12}+37300793498734 h^{11}+319959385783550 h^{10}+1699041589290610 h^9\\
\phantom{\tilde{f}_2:=}
+3768001178644240 h^8-563210238257358 h^7-11130338920636646 h^6+430259282513022 h^5\\
\phantom{\tilde{f}_2:=}
+13828775175467166 h^4-7823676892075654 h^3-929819855238438 h^2+1585805085637510 h\\
\phantom{\tilde{f}_2:=}
-300086925057803) \kappa^{11}+(\frac{37836610685814319}{65536} h^5-\frac{28933067785322487}{1048576} h-\frac{1467713483283}{65536} h^{15}\\
\phantom{\tilde{f}_2:=}
-\frac{3017331869259}{1048576} h^{16}-\frac{8051201309996249}{131072} h^3-\frac{10811316295905353}{32768} h^7+\frac{48940599416021407}{524288} h^2-\frac{774184385182879}{131072} h^{12}\\
\phantom{\tilde{f}_2:=}
-\frac{1823549625}{524288} h^{18}-\frac{3622360116941469}{131072} h^{11}-\frac{8783744674294019}{262144} h^{10}-\frac{78567797467371}{131072} h^{13}-\frac{6000551217195}{65536} h^{14}\\
\phantom{\tilde{f}_2:=}
-\frac{18392217328617501}{65536} h^4+\frac{2607920428387489}{1048576}-\frac{173889793043}{1048576} h^{17}+\frac{56428265198672697}{524288} h^9+\frac{98318295117756985}{524288} h^8\\
\phantom{\tilde{f}_2:=}
-\frac{13301070903895537}{65536} h^6) \kappa^{10}+(\frac{3916196660238163}{262144} h^{12}+\frac{352673051554887}{262144} h^{13}+\frac{161939442055308807}{262144} h^6\\
\phantom{\tilde{f}_2:=}
+\frac{594374342567591}{8192} h^{11}-\frac{8571683172835871}{2097152}+\frac{198122527844521537}{262144} h^7+\frac{959767669343}{2097152} h^{17}+\frac{10584251805}{1048576} h^{18}\\
\phantom{\tilde{f}_2:=}
-\frac{519146845018170287}{1048576} h^8-\frac{262129247245476701}{1048576} h^9+\frac{43149482239296513}{131072} h^3+\frac{15703258539077}{2097152} h^{16}\\
\phantom{\tilde{f}_2:=}
-\frac{180662883941849255}{131072} h^5+\frac{49391534744377}{262144} h^{14}+\frac{133399420063046291}{2097152} h+\frac{14009420388893}{262144} h^{15}\\
\phantom{\tilde{f}_2:=}
+\frac{30138868187655269}{65536} h^4+\frac{50919335889055481}{524288} h^{10}-\frac{284481558973120911}{1048576} h^2) \kappa^9+(\frac{1631592546187268869}{2097152} h^8\\
\phantom{\tilde{f}_2:=}
-\frac{44164644626293731}{262144} h^{10}-\frac{373828424377745}{262144} h^{13}+\frac{701633845403018371}{2097152} h^9+\frac{13629452279521457}{4194304}\\
\phantom{\tilde{f}_2:=}
-\frac{160209984064619469}{524288} h^4-\frac{18136958145}{1048576} h^{18}-\frac{265150707542975635}{262144} h^7-\frac{327624048583336369}{4194304} h-\frac{42116258021715}{262144} h^{14}\\
\phantom{\tilde{f}_2:=}
-\frac{572338303457}{8192} h^{15}+\frac{438342605782982473}{1048576} h^2-\frac{2806571550242451}{131072} h^{12}-\frac{58856803281920273}{524288} h^{11}\\
\phantom{\tilde{f}_2:=}
-\frac{71910528092009819}{65536} h^6-\frac{3110902671901}{4194304} h^{17}-\frac{369808982663851785}{524288} h^3+\frac{1009882244567901435}{524288} h^5\\
\phantom{\tilde{f}_2:=}
-\frac{47313028912627}{4194304} h^{16}) \kappa^8+(\frac{81351619125}{4194304} h^{18}-\frac{904661121754792377}{524288} h^5-\frac{10764990449060821}{8388608}\\
\phantom{\tilde{f}_2:=}
+\frac{237123049165732323}{2097152} h^{11}-\frac{204101442145165293}{2097152} h^4-\frac{1177799048768276503}{4194304} h^9+\frac{105807374435039}{2097152} h^{15}\\
\phantom{\tilde{f}_2:=}
+\frac{472878761924713581}{8388608} h+\frac{1786840876165227489}{2097152} h^3-\frac{12117673683849}{262144} h^{14}+\frac{90195913848163}{8388608} h^{16}\\
\phantom{\tilde{f}_2:=}
+\frac{39208636221048247}{2097152} h^{12}+\frac{6544841008441}{8388608} h^{17}+\frac{353851473863187}{1048576} h^{13}-\frac{1621405769881797807}{4194304} h^2\\
\phantom{\tilde{f}_2:=}
+\frac{406356447512864569}{2097152} h^{10}+\frac{1816588892031505853}{2097152} h^7-\frac{3392140058111290011}{4194304} h^8+\frac{668088573865345243}{524288} h^6) \kappa^7\\
\phantom{\tilde{f}_2:=}
+(\frac{4906942566037650943}{8388608} h^8-\frac{649728266860559659}{4194304} h^{10}+\frac{3668430461372279}{16777216}-\frac{4279401476414670115}{4194304} h^6\\
\phantom{\tilde{f}_2:=}
+\frac{469847566540045935}{2097152} h^2-\frac{2702727036469346781}{4194304} h^3-\frac{2027060077374925129}{4194304} h^7-\frac{54909832140995}{4194304} h^{15}\\
\phantom{\tilde{f}_2:=}
+\frac{747313276356437947}{2097152} h^4+\frac{1029904104097003}{1048576} h^{13}-\frac{9309843187639}{16777216} h^{17}-\frac{5042439163436737}{524288} h^{12}\\
\phantom{\tilde{f}_2:=}
+\frac{1075315786862599}{4194304} h^{14}+\frac{2121211204107716279}{2097152} h^5-\frac{323972986568092983}{4194304} h^{11}+\frac{1246553832491357653}{8388608} h^9\\
\end{array}
\end{eqnarray*}
\begin{eqnarray*}
\begin{array}{l}
\phantom{\tilde{f}_2:=}
-\frac{110614193055149}{16777216} h^{16}-\frac{408234751472010203}{16777216} h-\frac{62482414455}{4194304} h^{18}) \kappa^6+(\frac{17663658730781923}{8388608} h^{12}\\
\phantom{\tilde{f}_2:=}
-\frac{748161364176247753}{16777216} h^9+\frac{9021391408501}{33554432} h^{17}-\frac{3106230583202516513}{8388608} h^5-\frac{2667804785651032537}{8388608} h^4\\
\phantom{\tilde{f}_2:=}
-\frac{5058813514299363335}{16777216} h^8+\frac{691174193024603405}{4194304} h^7+\frac{151604802151654929}{4194304} h^{11}+\frac{1320025786764121983}{4194304} h^3\\
\phantom{\tilde{f}_2:=}
+\frac{737423451445195389}{8388608} h^{10}-\frac{1355534311385745033}{16777216} h^2-\frac{287052507131985}{33554432}+\frac{82319816153271}{33554432} h^{16}\\
\phantom{\tilde{f}_2:=}
-\frac{11634249381412761}{8388608} h^{13}-\frac{1219928679573601}{4194304} h^{14}-\frac{42741108884397}{4194304} h^{15}+\frac{2411657695473752887}{4194304} h^6\\
\phantom{\tilde{f}_2:=}
+\frac{200000621140535685}{33554432} h+\frac{133653863415}{16777216} h^{18}) \kappa^5+(-\frac{819746264216279837}{8388608} h^3-\frac{1904307438061698571}{8388608} h^6\\
\phantom{\tilde{f}_2:=}
+\frac{3740725840940176885}{33554432} h^8-\frac{95911325678712155}{8388608} h^{11}+\frac{2637407836167327937}{16777216} h^4+\frac{98785446239441163}{33554432} h^9\\
\phantom{\tilde{f}_2:=}
+\frac{1544388019883439}{8388608} h^{14}+\frac{15728503582372479}{16777216} h^{13}-\frac{193201212377419843}{8388608} h^7-\frac{49533585015}{16777216} h^{18}\\
\phantom{\tilde{f}_2:=}
+\frac{10298833250390657}{16777216} h^{12}+\frac{1114769960820619167}{16777216} h^5+\frac{285342104847206251}{16777216} h^2+\frac{18009884805255}{67108864}\\
\phantom{\tilde{f}_2:=}
-\frac{149202897831001767}{4194304} h^{10}-\frac{27055878477817}{67108864} h^{16}-\frac{46326596042412903}{67108864} h-\frac{5822578274983}{67108864} h^{17}\\
\phantom{\tilde{f}_2:=}
+\frac{101711920551979}{8388608} h^{15}) \kappa^4+\frac{1}{33554432} h (24786556740 h^{17}+584913432490 h^{16}-2342497331229 h^{15}\\
\phantom{\tilde{f}_2:=}
-197939029933495 h^{14}-2432717286208029 h^{13}-12721660505193191 h^{12}\\
\phantom{\tilde{f}_2:=}
-20154534523490017 h^{11}+79244101754791277 h^{10}+341135256535815927 h^9\\
\phantom{\tilde{f}_2:=}
+110110170694944625 h^8-980452616030700471 h^7-218952263689161509 h^6\\
\phantom{\tilde{f}_2:=}
+2038228848786011337 h^5+179759392129858427 h^4-1564177676844557451 h^3\\
\phantom{\tilde{f}_2:=}
+602709015627097519 h^2-58657218181534119 h+600423209256465) \kappa^3-\frac{1}{134217728} h \\
\phantom{\tilde{f}_2:=}
\cdot(15814746165 h^{17}+243182479141 h^{16}-6727511541654 h^{15}-214827403401246 h^{14}\\
\phantom{\tilde{f}_2:=}
-2369157028586874 h^{13}-12499828242091954 h^{12}-25557423142145934 h^{11}\\
\phantom{\tilde{f}_2:=}
+40740153361625706 h^{10}+267353481272472560 h^9+189670062636941912 h^8\\
\phantom{\tilde{f}_2:=}
-699544766664324290 h^7-540367923019593626 h^6+1366928930536656346 h^5\\
\phantom{\tilde{f}_2:=}
+757965148352724594 h^4-1056146317659794202 h^3+217056976339137102 h^2\\
\phantom{\tilde{f}_2:=}
-6307841010550149 h+34495848583875) \kappa^2+\frac{1}{67108864} h^2 (714213105 h^{16}+1970792740 h^{15}\\
\phantom{\tilde{f}_2:=}
-667227521572 h^{14}-15852386236908 h^{13}-163368504986864 h^{12}-859564456329388 h^{11}\\
\phantom{\tilde{f}_2:=}
-1943968636814692 h^{10}+1607925129586244 h^9+16599297481440266 h^8+17615966503801900 h^7\\
\phantom{\tilde{f}_2:=}
-38258848104755532 h^6-56255239628904068 h^5+59913343377313336 h^4+77458565024388700 h^3\\
\phantom{\tilde{f}_2:=}
-40865206572185772 h^2+2404177561720620 h-34495848583875) \kappa-\frac{1}{134217728} h^3 (54073725 h^{15}\\
\phantom{\tilde{f}_2:=}
-1011489335 h^{14}-97189416455 h^{13}-1986147874307 h^{12}-19483918887719 h^{11}\\
\phantom{\tilde{f}_2:=}
-101131586545787 h^{10}-233547778816235 h^9+163691764889881 h^8+2031963839835047 h^7\\
\phantom{\tilde{f}_2:=}
+2740202277870123 h^6-3829178058687173 h^5-9114195900415897 h^4+3080665291005891 h^3\\
\phantom{\tilde{f}_2:=}
+11146213901963847 h^2-866187180306441 h+34495848583875).\\
\end{array}
\end{eqnarray*}
The resultant of $\tilde{f}_1$ and $\tilde{f}_2$ by eliminating variable $\kappa$ denoted as $\tilde{\mathcal{R}}_1$,
where
\begin{eqnarray*}
\begin{array}{l}
\tilde{\mathcal{R}}_1:=\mbox{res}(\tilde{f}_1, \tilde{f}_2, \kappa)
=-\frac{1}{72057594037927936} h^4 (h-1)^2 (h+23)^2 (2 h-1)^2 (3 h^3+16 h^2+15 h+14)^2\\
\phantom{\mathcal{R}_1=}
 \cdot(h^2+2 h-7)^4 (h+3)^6 (4 h^3+6 h^2-15 h+21) (69324 h^{25}+3784256 h^{24}+88925415 h^{23}\\
\phantom{\mathcal{R}_1=}
+1210510049 h^{22}+9321912097 h^{21}+31150056402 h^{20}-48628636343 h^{19}-621513307905 h^{18}\\
\phantom{\mathcal{R}_1=}
-449862701721 h^{17}+6309033684670 h^{16}+9292970454326 h^{15}-38554634019094 h^{14}\\
\phantom{\mathcal{R}_1=}
-64709338716470 h^{13}+130467748580116 h^{12}+188895548702946 h^{11}-213048078842242 h^{10}\\
\phantom{\mathcal{R}_1=}
-167516634645510 h^9+80910436459484 h^8-138237240652045 h^7+193470794940005 h^6\\
\phantom{\mathcal{R}_1=}
+135195659870885 h^5-204505372471814 h^4+105439411571397 h^3-32169849655725 h^2\\
\phantom{\mathcal{R}_1=}
+3724806231363 h-132182987706)\\
\end{array}
\end{eqnarray*}
Solving $\tilde{\mathcal{R}}_1=0$ yields five zeros $x_1$ in $\mathcal{P}$, i.e., $0.0658$, $ 0.1123$, $0.4822$, $0.5$ and $0.9627$.
Applying the successive pseudo-division to the equations $\tilde{f}_1=0$ and $\tilde{f}_2=0$ we get a solution $\kappa=\frac{\kappa_n}{\kappa_d}$, where
\begin{eqnarray*}
\begin{array}{l}
\kappa_n:=-h (104424 h^{36}+15806636 h^{35}+1060235958 h^{34}+41360310493 h^{33}+1033237688608 h^{32}\\
\phantom{\kappa_n:=}
+17093052680880 h^{31}+185234034631326 h^{30}+1214782758103459 h^{29}+3410711864703916 h^{28}\\
\phantom{\kappa_n:=}
-10263855243151606 h^{27}-100244026439753366 h^{26}-87721813509906287 h^{25}\\
\phantom{\kappa_n:=}
+1171764771207211552 h^{24}+2240927152940620996 h^{23}-9772808464958064482 h^{22}\\
\phantom{\kappa_n:=}
-22749592143215129389 h^{21}+57638777043031711252 h^{20}+121537499193080672062 h^{19}\\
\phantom{\kappa_n:=}
-246262891452996520490 h^{18}-282788105181178447125 h^{17}+766714881123795937480 h^{16}\\
\phantom{\kappa_n:=}
-109543664161325840664 h^{15}-1226489872897021980470 h^{14}+1630300109228497230369 h^{13}\\
\phantom{\kappa_n:=}
-806552374484503612940 h^{12}-1077797869120164446890 h^{11}+4106839101181189996638 h^{10}\\
\phantom{\kappa_n:=}
-3848854638191389868333 h^9+845991926265373897616 h^8+2063653729639496063332 h^7\\
\phantom{\kappa_n:=}
-4536379008677843744262 h^6+4748641501899313027689 h^5-3172898299865472357660 h^4\\
\phantom{\kappa_n:=}
+1353005890097332372470 h^3-365758729089763817556 h^2+58630913767660198356 h\\
\phantom{\kappa_n:=}
-2104881200090531352),\\
\kappa_d:=
54480 h^{37}-25280 h^{36}-429890512 h^{35}-29584573272 h^{34}-994950290565 h^{33}-20178839743180 h^{32}\\
\phantom{\kappa_n:=}
-258930573495916 h^{31}-2018820796155987 h^{30}-7678410965617803 h^{29}+6781705163337138 h^{28}\\
\phantom{\kappa_n:=}
+174416098225939254 h^{27}+344695353057908799 h^{26}-1704274219197667845 h^{25}\\
\phantom{\kappa_n:=}
-5694732570594415444 h^{24}+12128022268874697944 h^{23}+49503457543617590261 h^{22}\\
\phantom{\kappa_n:=}
-66857275387628213603 h^{21}-246372492634932965506 h^{20}+309415046153421965826 h^{19}\\
\phantom{\kappa_n:=}
+585374819884899902511 h^{18}-1139312137822506612103 h^{17}+11079676686133307176 h^{16}\\
\phantom{\kappa_n:=}
+2115483758750366890388 h^{15}-2796158967212325979553 h^{14}+1112696763869892676919 h^{13}\\
\phantom{\kappa_n:=}
+2422679013780879493454 h^{12}-7011710911652753395142 h^{11}+6151475587496528241397 h^{10}\\
\phantom{\kappa_n:=}
-1447591488200152741791 h^9-3831489956406538933284 h^8+8118838419345304712896 h^7\\
\phantom{\kappa_n:=}
-8146426987531933537793 h^6+5606696051148996416871 h^5-2476126155714773703294 h^4\\
\phantom{\kappa_n:=}
+691746159028922895774 h^3-118440927041560529931 h^2+3814475091599076912 h\\
\phantom{\kappa_n:=}
-2631715218685476.
\end{array}
\end{eqnarray*}
Among these five positive zeros of $\tilde{\mathcal{R}}_1=0$, only the third one $x_1\doteq0.4822$ is a feasible zero for $(\kappa,h,x_1)\in\mathcal{P}$, i.e.,
if $x_1\doteq0.4822$, then $\kappa\doteq493.2405>0.6516$ and $\tilde{f}_1=\tilde{f}_2=0$ and $\tilde{f}_3=-1.2279\times 10^{31}<0$.
Hence, it is possible that $V(f_1, f_2, f_3)=\emptyset$ but $V(f_1, f_2)\neq\emptyset$ in $\mathcal{P}$.
The computation of algebraic varieties $V(f_1, f_2)$ and $V(f_1, f_2, f_3)$ in the proof of Theorem \ref{thm4} is completed.

\end{document}